# Surreal substructures

Vincent Bagayoko

Fachbereich Mathematik und Statistik
Universität Konstanz
D-78457 Konstanz
Germany

Joris van der Hoeven

CNRS, LIX
Campus de l'École polytechnique
1, rue Honoré d'Estienne d'Orves
Bâtiment Alan Turing, CS35003
91120 Palaiseau, France

*May 3, 2023*

### Abstract

Conway's field **No** of surreal numbers comes both with a natural total order and an additional "simplicity relation" which is also a partial order. Considering **No** as a doubly ordered structure for these two orderings, an isomorphic copy of **No** into itself is called a *surreal substructure*. It turns out that many natural subclasses of **No** are actually of this type. In this paper, we study various constructions that give rise to surreal substructures and analyze important examples in greater detail.

## 1 Introduction

### 1.1 Surreal numbers

The class **No** of *surreal numbers* was discovered by Conway and studied in his well-known monograph *On Numbers and Games* [14]. Conway's original definition is somewhat informal and goes at follows:

> "If $L$ and $R$ are any two sets of (surreal) numbers, and no member of $L$ is $\geqslant$ any member of $R$, then there is a (surreal) number $\{L|R\}$. All (surreal) numbers are constructed in this way."

The magic of surreal numbers lies in the fact that many traditional operations on integers and real numbers can be defined in a very simple way on surreal numbers. Yet, the class **No** turns out to admit a surprisingly rich algebraic structure under these operations. For instance, the sum of two surreal numbers $x = \{x_L | x_R\}$ and $y = \{y_L | y_R\}$ is defined recursively by

$$x + y = \{x_L + y, x + y_L | x_R + y, x + y_R\}. \qquad (1.1)$$

In section 3 below, we recall similar definitions for subtraction and multiplication. Despite the fact that the basic arithmetic operations can be defined in such an "effortless" way, Conway showed that **No** actually forms a real-closed field that contains $\mathbb{R}$. Strictly speaking, some care is required here, since the surreal numbers **No** form a proper class. In particular, it contains all ordinal numbers $\alpha = \{\alpha_L | \emptyset\}$. We refer to appendix B for ways to deal with this kind of set-theoretic issues.

One convenient way to rigorously introduce surreal numbers $x$ is to regard them as "sign sequences" $x = (x[\beta])_{\beta < \alpha} \in \{-1, +1\}^\alpha$ indexed by the elements $\beta < \alpha$ of an ordinal number $\alpha = \ell(x)$, called the *length* of $x$: see section 2.1 below for details. Every ordinal $\alpha$ itself is represented as $\alpha = (\alpha[\beta])_{\beta < \alpha}$ with $\alpha[\beta] = 1$ for all $\beta < \alpha$. The number $1/2$ is represented by the sign sequence $+1, -1$ of length $2$. The ordering $\leqslant$ on **No** corresponds to the lexicographical ordering on sign sequences, modulo zero padding when comparing two surreal numbers of different lengths. The sign sequence representation also induces the important notion of *simplicity*: given $x, y \in \mathbf{No}$, we say that $x$ is *simpler* as $y$, and write $x \sqsubseteq y$, if the sign sequence of $x$ is a truncation of the sign sequence of $y$. The simplicity relation is denoted by $\leq_s$ in some previous works [11, 30, 3].

The sign sequence representation was introduced and studied systematically in Gonshor's book [24]. As we will see in section 3, it also allows for a natural extension of ordinal arithmetic to the surreal numbers. In order to avoid confusion, we will systematically use the notations $\alpha \dotplus \beta$ and $\alpha \dottimes \beta$ for ordinal sums and products and $\dot{\omega}^\beta$ for ordinal exponentiation. For instance, in **No**, we have $\omega \dotplus 1 = \omega + 1 = 1 + \omega \neq 1 \dotplus \omega = \omega$. Given an ordinal $\alpha$, it is also natural to define the set $\mathbf{No}(\alpha)$ of all surreal numbers $x$ of length $\ell(x) < \alpha$. It turns out that $\mathbf{No}(\alpha)$ is a real-closed subfield of **No** if and only if $\alpha$ is an $\varepsilon$-number, i.e. $\dot{\omega}^\alpha = \alpha$ [15, Proposition 4.7 and Corollary 4.9].

## 1.2 Exponentiation, derivation, and hyperseries

Quite some work has been dedicated to the extension of basic calculus to the surreal numbers and to the study of various operations in terms of sign sequences. In his book [24], Gonshor shows how to extend the real exponential function to **No**. This exponential function actually admits the same first order properties as the usual exponential function: the class **No** is elementarily equivalent to $\mathbb{R}$ as an exponential field. In fact, they are even elementarily equivalent as real exponential ordered fields equipped with restricted analytic functions [15, Theorem 2.1]. Here we recall that a restricted real analytic function is a power series $f \in \mathbb{R}[[x]]$ at the origin that converges on a small closed ball $[-r, r]$ with $r > 0$. Then it can be shown that the definition of $f(x)$ extends to surreal numbers $x$ with $-r \leqslant x \leqslant r$.

Another important question concerns the possibility to define a natural derivation $\partial$ on the surreal numbers, which is non-trivial in the sense that $\partial \omega = 1$. Such a derivation was first constructed by Berarducci and Mantova [11], while making use of earlier work by van der Hoeven and his student Schmeling [38]. It was shown in [3] that this "Italian" derivation $\partial_{\mathrm{BM}}$ has "similarly good properties" as the exponential function in the sense that **No** is elementary equivalent to the field of transseries as an H-field. Here transseries are a generalization of formal power series. They form an ordered exponential field $\mathbb{T}$ that comes with a derivation. The notion of an H-field captures the algebraic properties of this field $\mathbb{T}$ as well as those of so-called Hardy fields. We refer to [1] for more details.

The above results on the exponential function and the Italian derivation $\partial_{\mathrm{BM}}$ on **No** rely on yet another representation of surreal numbers as generalized power series $x = \sum_{\mathfrak{m} \in \mathbf{Mo}} x_\mathfrak{m} \mathfrak{m}$ with real coefficients and *monomials* $\mathfrak{m} \in \mathbf{Mo}$ such that $\mathfrak{m}$ is simpler than any other $0 < x \in \mathbf{No}$ with the same valuation as $\mathfrak{m}$: see section 2.3 for details. Indeed, ordinary power series and Laurent series in $\omega^{-1}$ can be regarded as *functions* in $\omega$, so they come with a natural derivation. More generally, the exponential function on **No** makes it possible to interpret any transseries in $\omega$ as a surreal number, which makes it again possible to derive such surreal numbers in a natural way.

Unfortunately, not all surreal numbers are transseries in $\omega$. For instance, the surreal number $\{\omega, e^\omega, e^{e^\omega}, \ldots | \emptyset\}$ is larger than any transseries in $\omega$. In order to be able to intepret all surreal numbers as functions in $\omega$, two ingredients are missing: on the one hand, we need to introduce ordinal "iterators" $E_\alpha$ of the exponential function that grow faster than finite iterates. For instance, we have $E_\omega(\omega) = \{\omega, e^\omega, e^{e^\omega}, \ldots | \emptyset\}$. On the other hand, we need to be able to represent so-called nested transseries such as

$$\sqrt{\omega} + e^{\sqrt{\log \omega} + e^{\sqrt{\log\log \omega} + e^{\cdot^{\cdot^{\cdot}}}}}. \tag{1.2}$$

In [29, 2], it was conjectured that any surreal number can be represented as a generalized "hyperseries" in $\omega$, when taking these above observations into account in a suitable manner. For early progress on the "series side", we refer to [26, 38, 29, 16]. In a series of preprints [8, 7, 5, 6], we recently managed to prove the conjecture. This work is quite technical and critically relies on the more accessible results in this paper.

The derivation $\partial_{BM}$ cannot be compatible with a composition law on **No** [12, Theorem 8.4]. More specifically, it was noted in [2] that the Italian derivation fails to satisfy $\partial_{BM}(E_\omega(x)) = (\partial_{BM} x) E'_\omega(x)$ for all $x$. Ultimately, the ability to represent surreal numbers as hyperseries evaluated at $\omega$ should lead to compatible definitions of a derivation and a composition on **No**.

## 1.3 Surreal substructures

In the course of our project to construct an isomorphism between **No** and a suitable class of hyperseries, we frequently encountered subclasses **S** of **No** that are naturally parameterized by **No** itself. For instance, Conway's generalized ordinal exponentiation $x \in \mathbf{No} \longmapsto \dot{\omega}^x \in \mathbf{Mo}$ is bijective, which leads to a natural parameterization of the class **Mo** of monomials by **No** (see Theorems 5.2 and 5.11). Similarly, nested expressions such as (1.2) do not give rise to a single surreal number, but rather to a class **Ne** of surreal numbers that is naturally parameterized by **No** (see Theorem 8.8). Yet another example is the class $\mathbf{La} = \bigcap_{n \in \mathbb{N}} \{(\exp \circ \overset{n \times}{\cdots} \circ \exp)(\mathfrak{m}) : \mathfrak{m} \in \mathbf{Mo}, \mathfrak{m} > \mathbb{R}\}$ of *log-atomic* surreal numbers that occurs crucially in the construction of derivations on **No** [11, Section 5.2].

In these three examples, the parameterizations turn out to be more than mere bijective maps: they actually preserve both the ordering $\leqslant$ and the simplicity relation $\sqsubseteq$. This leads to the definition of a *surreal substructure* of **No** as being an isomorphic copy of $(\mathbf{No}, \leqslant, \sqsubseteq)$ inside itself. Surreal substructures such as **Mo**, **Ne**, and **La** behave similarly as the surreal numbers themselves **No** in many regards. In our project, we have started to exploit this property for the definition and study of new functions on **No** such as hyperlogarithms and nested transseries.

The main goal of the present paper is to develop the basic theory of surreal substructures for its own sake and as a new tool to study surreal numbers. We hope to convey the sense that surreal substructures are at the same time very general and very rigid subclasses of **No** and that several problems regarding the enriched structure of **No** (highlighted in particular in the work of Gonshor [24], Lemire [31, 32, 33], Ehrlich [19, 18, 20], Kuhlmann–Matusinski [30], Berarducci–Mantova [11], and Aschenbrenner–van den Dries–van der Hoeven [3]) crucially involve surreal substructures. Even for very basic subclasses of **No**, such as $\mathbf{No}^{>} = \{x \in \mathbf{No} : x > 0\}$, we suggest that it deserves our attention when they form surreal substructures.

A substantial part of our paper (namely, sections 4, 5, and 6) is therefore devoted to basic but fundamental results. Some of these general facts (or variants of them) were known and rediscovered in different contexts [34, 19]. However, they mainly appeared as auxiliary tools in these works. In this paper, we aim at covering the most noteworthy facts in a self-contained and organized way. In the course of our exposition, we identify which properties of surreal substructures are systematic and which ones are proper to specific structures. We also include a wide range of examples. This effort culminates in the last two sections 7 and 8, where we present the examples that motivated our paper and that are important for our program to construct an isomorphism between **No** and the class of hyperseries. In Appendix A, we also compiled a small atlas for the most prominent examples of surreal substructures.

## 1.4 Summary of our contributions

Let us briefly outline the structure of the paper. In section 2, we recall the three main representations of surreal numbers. In section 3, we recall the definitions of basic arithmetic operations on surreal numbers. We also show how to extend the ordinal sum $\dotplus$ and the ordinal product $\dottimes$ to **No**.

In section 4, we introduce *surreal substructures*, our main object of study, as isomorphic copies of $(\mathbf{No}, \leqslant, \sqsubseteq)$ inside itself. Any surreal substructure **S** comes with a *defining isomorphism* $\Xi_\mathbf{S}: (\mathbf{No}, \leqslant, \sqsubseteq) \longrightarrow (\mathbf{S}, \leqslant, \sqsubseteq)$ that is unique and that we consider as a parameterization of the elements in **S** by **No**. Defining isomorphisms $\Xi_\mathbf{S}$ and $\Xi_\mathbf{T}$ can be composed to form the defining isomorphism $\Xi_\mathbf{U} = \Xi_\mathbf{S} \circ \Xi_\mathbf{T}$ of a new surreal structure $\mathbf{U} = \mathbf{S} \prec \mathbf{T}$ that we call the *imbrication* of **T** inside **S**. More generally, we will often switch between the study of surreal substructures and that of their parameterizations. A consequent part of section 4.1 is a reformulation of notions and arguments found in [34, 19, 20]; see Remark 4.8.

In section 5, we investigate the existence of *fixed points* for the defining isomorphism $\Xi_\mathbf{S}$ of a given surreal substructure **S**. More precisely, we give conditions on **S** under which the class $\mathbf{Fix_S}$ of such fixed points is itself a surreal substructure. Determining the class $\mathbf{Fix_S}$ allows us in some cases to compare the defining isomorphisms of two surreal substructures. This task leads us to study surreal substructures **S** which are closed under non-empty, set-sized suprema in $(\mathbf{No}, \sqsubseteq)$ of chains in $(\mathbf{S}, \sqsubseteq)$. Such a surreal substructure **S** is said to be **No**-*closed*, and has the following properties:

- Corollary 5.14: for an **No**-closed surreal substructure **S**, the class $\mathbf{Fix_S}$ is a surreal substructure, and it coincides with $\bigcap_{n \in \mathbb{N}} \Xi_\mathbf{S}^n(\mathbf{No})$, where $\Xi_\mathbf{S}^n$ denotes the $n$-fold composition of $\Xi_\mathbf{S}$ with itself. A similar result was first proved by Lurie [34, Theorem 8.2]; see Remark 5.15.

- Proposition 5.18: for an **No**-closed surreal substructure **S**, there is a decreasing sequence $(\mathbf{S}^{\prec \alpha})_{\alpha \in \mathbf{On}}$ of surreal substructures such that for ordinals $\alpha, \beta$, we have

    a) $\mathbf{S}^{\prec 0} = \mathbf{No}$ and $\mathbf{S}^{\prec 1} = \mathbf{S}$,

    b) $\mathbf{S}^{\prec(\alpha \dotplus \beta)} = \mathbf{S}^{\prec \alpha} \prec \mathbf{S}^{\prec \beta}$,

    c) $\mathbf{S}^{\prec(\alpha \dottimes \beta)} = (\mathbf{S}^{\prec \alpha})^{\prec \beta}$,

    d) $\mathbf{S}^{\prec \alpha} = \bigcap_{\gamma < \alpha} \mathbf{S}^{\prec \gamma}$ if $0 < \alpha$ is limit,

In fact any well-ordered sequence of **No**-closed surreal substructures can be similarly "imbricated", and thus **No**-closed surreal substructures can be seen as words in a rich language that captures at the same time the notions of fixed points, imbrications and intersections of surreal substructures. One direct application is a new proof of a theorem by Lemire [32]; see Remark 5.17.

In section 6, we study subclasses **Smp**$_\Pi$ whose elements are the simplest representatives of members in a convex partition $\Pi$ of a surreal substructure **S**. Under a set-theoretic condition on $\Pi$, we prove that this class forms a surreal substructure of **S** (Theorem 6.7) whose parameterization admits a short recursive definition. A weaker version of this theorem was first proved by Lurie [34]; see Remark 6.8. A particularly important special case is when the convex partition is induced by a group action (see section 6.3). We also introduce the notion of a *sharp* convex partition $\Pi$ of a surreal substructure **S** which makes **Smp**$_\Pi$ closed within **S** (Theorem 6.14).

Our final sections 7 and 8 concern the application of our results to some prominent examples of surreal substructures. This includes the structure **No**$_>$ of *purely infinite surreal numbers* of [24], the structure **Mo** of *monomials* of [14], the structure **La** of *log-atomic numbers* of [11], the structure **K** of *κ-numbers* of [30], and various structures of *nested monomials*, including **Ne**. Our results about nested monomials in section 8 are analogous to Lemire's work on continued exponential expressions [33], when replacing ordinal exponentiation by traditional exponentiation. The appendix A contains a short overview of the surreal substructures encountered in this paper.

## 1.5 Notations

We will systematically use a bold type face to denote classes such as **No** that may not be sets. Given a partially ordered class $(\mathbf{X}, <_\mathbf{X})$ and subclasses $\mathbf{A}, \mathbf{B}$ of $\mathbf{X}$, we write $\mathbf{A} <_\mathbf{X} \mathbf{B}$ if $a <_\mathbf{X} b$ for all $a \in \mathbf{A}$ and $b \in \mathbf{B}$. This holds in particular whenever $\mathbf{A} = \emptyset$ or $\mathbf{B} = \emptyset$. For elements $x_1, ..., x_n, y_1, ..., y_n$ of $\mathbf{X}$, we write $x_1, ..., x_m <_\mathbf{X} \mathbf{B}$ and $\mathbf{A} <_\mathbf{X} y_1, ..., y_n$ instead of $\{x_1, ..., x_m\} <_\mathbf{X} \mathbf{B}$ and $\mathbf{A} <_\mathbf{X} \{y_1, ..., y_n\}$. Given more than two subclasses $\mathbf{A}_1, ..., \mathbf{A}_n$ of $\mathbf{X}$, we also write $\mathbf{A}_1 <_\mathbf{X} \cdots <_\mathbf{X} \mathbf{A}_n$ whenever $\mathbf{A}_i <_\mathbf{X} \mathbf{A}_j$ for all $i < j$.

If $x \in \mathbf{X}$, we let $\mathbf{X}^{>x}$ denote the class of elements $y \in \mathbf{X}$ with $y > x$. In the special case when $(\mathbf{X}, e, \cdot, <_\mathbf{X})$ is an ordered monoid, we simply write $\mathbf{X}^> = \mathbf{X}^{>e}$ and $\mathbf{X}^< = \mathbf{X}^{<e}$.

We use similar notations for non-strict orders $\leqslant_\mathbf{X}$.

## 2 Different presentations of surreal numbers

Surreal numbers can be represented in three main ways: as sign sequences, as generalized Dedekind cuts, and as generalized power series over $\mathbb{R}$. In this section, we briefly recall how this works, and review the specific advantages of each representation. We refer to [14, 24, 19, 18, 35] for more details.

### 2.1 Surreal numbers as sign sequences

The sign sequence representation is most convenient for the rigourous development of the basic theory of surreal numbers. It was introduced by Gonshor [24, page 3] and we will actually use it to formally define surreal numbers as follows:

**Definition 2.1.** *A surreal number is a map* $x: \ell(x) \longrightarrow \{-1,1\}; \alpha \longmapsto x[\alpha]$, *where* $\ell(x) \in \mathbf{On}$ *is an ordinal number. We call* $\ell(x)$ *the* **length** *of x and the map* $x: \ell(x) \longrightarrow \{-1,1\}$ *the* **sign sequence** *of x. We write* **No** *for the class of surreal numbers.*

It follows from this definition that **No** is a proper class. Given a surreal number $x \in \mathbf{No}$, it is convenient to extend its sign sequence with zeros to a map $\mathbf{On} \longrightarrow \{-1,0,1\}$ and still denote this extension by $x$. In other words, we take $x[\alpha] = 0$ for all $\alpha \geqslant \ell(x)$. Given $x \in \mathbf{No}$ and $\alpha \in \mathbf{On}$, we also introduce its *restriction* $y = x \upharpoonright \alpha \in \mathbf{No}$ to $\alpha$ as being the zero padded restriction of the map $x$ to $\alpha$: we set $y[\beta] = x[\beta]$ for $\beta < \alpha$ and $y[\beta] = 0$ for $\beta \geqslant \alpha$.

The first main relation on **No** is its *ordering* $\leqslant$. We define it to be the restriction of the lexicographical ordering on the set of all maps from **On** to $\{-1,0,1\}$. More precisely, given distinct elements $x, y \in \mathbf{No}$, there exists a smallest ordinal $\alpha$ with $x[\alpha] \neq y[\alpha]$. Then we define $x < y$ if and only if $x[\alpha] < y[\alpha]$.

The second main relation on **No** is the simplicity relation $\sqsubseteq$: given numbers $x, y \in \mathbf{No}$, we say that $x$ is *simpler* than $y$, and write $x \sqsubseteq y$, if $x = y \upharpoonright \ell(x)$. We write $x_\sqsubset = \{a \in \mathbf{No} : a \sqsubset x\}$ for the set of surreal numbers that are strictly simpler than $x$. The partially ordered class $(\mathbf{No}, \sqsubseteq)$ is well-founded, and $(x_\sqsubset, \sqsubseteq)$ is well-ordered with order type $\mathrm{ot}(x_\sqsubset, \sqsubseteq) = \ell(x)$.

Every linearly ordered—and thus well-ordered—subset $X$ of $(\mathbf{No}, \sqsubseteq)$ has a *supremum* $s = \sup_\sqsubseteq X$ in $(\mathbf{No}, \sqsubseteq)$: for any $x \in X$ and $\alpha < \ell(x)$, one has $s[\alpha] = x[\alpha]$; for any $\alpha \in \mathbf{On}$ with $\alpha \geqslant \ell(x)$ all $x \in X$, one has $s[\alpha] = 0$. We will only consider suprema in $(\mathbf{No}, \sqsubseteq)$ and never in $(\mathbf{No}, \leqslant)$. Numbers $x$ that are equal to $\sup_\sqsubseteq x_\sqsubset$ are called *limit numbers*; other numbers are called *successor numbers*. Limit numbers are exactly the numbers whose length is a limit ordinal.

## 2.2 Surreal numbers as simplest elements in cuts

If $L, R$ are *sets* of surreal numbers satisfying $L < R$, then there is a simplest surreal number, written $\{L | R\}$, which satisfies $L < \{L | R\} < R$ [24, Theorem 2.1]. We call $\{|\}$ the *Conway bracket*. Notice that $\{L | R\}$ is *the* simplest such number in the strong sense that for all $x \in \mathbf{No}$ with $L < x < R$, we have $\{L | R\} \sqsubseteq x$. The converse implication $\forall x \in \mathbf{No}, \{L | R\} \sqsubseteq x \Longrightarrow L < x < R$ may fail: see Remark 4.21 below.

Now consider two more sets $L', R'$ of surreal numbers with $L' < R'$. If $L$ has no strict upper bound in $L'$ and $R$ has no strict lower bound in $R'$, then we say that $(L, R)$ is *cofinal* with respect to $(L', R')$. We say that $(L, R)$ and $(L', R')$ are *mutually cofinal* if they are cofinal with respect to one another, in which case it follows that $\{L | R\} = \{L' | R'\}$. These definitions naturally extend to pairs $(\mathbf{L}, \mathbf{R})$ of classes with $\mathbf{L} < \mathbf{R}$. Note however that $\{\mathbf{L} | \mathbf{R}\}$ is not necessarily defined for such classes. Indeed, there may be no number $x$ with $\mathbf{L} < x < \mathbf{R}$ (e.g. for $\mathbf{L} = \mathbf{No}$ and $\mathbf{R} = \emptyset$).

We call a pair $(L, R)$ of sets with $L < R$ a *cut representation* of $\{L | R\}$. Such representations are not unique; in particular, we may replace $(L, R)$ by any mutually cofinal pair $(L', R')$. For every surreal number $x$, we denote

$$x_L = \{a \in \mathbf{No} : a < x, a \sqsubseteq x\}$$
$$x_R = \{a \in \mathbf{No} : a > x, a \sqsubseteq x\},$$

which are *sets* of surreal numbers. We call $x_L$ and $x_R$ the sets of *left* and *right options* for $x$. By [24, Theorem 2.8], one has $x = \{x_L | x_R\}$ and the pair $(x_L, x_R)$ is called the *canonical representation* of $x$.

This identity $x = \{x_L | x_R\}$ is the fundamental intuition behind Conway's definition of surreal numbers precisely as the simplest numbers lying in the "cut" defined by sets $L < R$ of simpler and previously defined surreal numbers. Of course, this is a highly recursive representation that implicitly relies on transfinite induction.

Conway's cut representation is attractive because it allows for the recursive definition of functions using by well-founded induction on $(\mathbf{No}, \sqsubseteq)$ or its powers. For instance, there is a unique bivariate function $f$ such that for all $x, y \in \mathbf{No}$, we have

$$f(x,y) = \{f(x_L, y), f(x, y_L) | f(x_R, y), f(x, y_R)\}. \tag{2.1}$$

Here we understand that $f(x_L, y), f(x, y_L)$ denotes the set $\{f(x', y) : x' \in x_L\} \cup \{f(x, y') : y' \in y_L\}$ and similarly for $f(x_R, y), f(x, y_R)$. This recursive definition is justified by the fact that the elements of the sets $x_L \times \{y\}, \{x\} \times y_L, x_R \times \{y\}$, and $\{x\} \times y_R$ are all strictly simpler than $(x, y)$ for the product order on $(\mathbf{No}, \sqsubseteq) \times (\mathbf{No}, \sqsubseteq)$. This precise equation is actually the one that Conway used to define the addition $+ = f$ on $\mathbf{No}$. We will recall similar definitions of a few other arithmetic operations in section 3 below.

## 2.3 Surreal numbers as well-based series

Let $C$ be a field and let $\mathfrak{M}$ be a totally ordered multiplicative group for the ordering $\preccurlyeq$. A subset $\mathfrak{S} \subseteq \mathfrak{M}$ is said to be *well-based* if it is well-ordered for the opposite ordering of $\preccurlyeq$ (i.e. there are no infinite chains $\mathfrak{m}_1 \prec \mathfrak{m}_2 \prec \cdots$ in $\mathfrak{M}$). A *well-based series* in $\mathfrak{M}$ and over $C$ is a map $f : \mathfrak{M} \longrightarrow C$ whose *support* $\operatorname{supp} f := \{\mathfrak{m} \in \mathfrak{M} : f(\mathfrak{m}) \neq 0\}$ is a well-based subset of $\mathfrak{M}$. Such a series is usually written as $f = \sum_{\mathfrak{m} \in \mathfrak{M}} f_\mathfrak{m} \mathfrak{m}$, where $f_\mathfrak{m} = f(\mathfrak{m})$ and the set of all such series is denoted by $C[[\mathfrak{M}]]$. Elements in $C$ and $\mathfrak{M}$ are respectively called *coefficients* and *monomials*. We call $\mathfrak{M}$ the *monomial group*. The support of any non-zero element $f \in C[[\mathfrak{M}]]$ admits a largest element for $\preccurlyeq$, which is called the *dominant monomial* of $f$ and denoted by $\mathfrak{d}_f$.

It was shown by Hahn [25] that $C[[\mathfrak{M}]]$ forms a field for the natural sum and the usual Cauchy convolution product

$$f + g := \sum_{\mathfrak{m} \in \mathfrak{M}} (f_\mathfrak{m} + g_\mathfrak{m}) \mathfrak{m}, \qquad fg := \sum_{\mathfrak{m} \in \mathfrak{M}} \left( \sum_{\mathfrak{v}\mathfrak{w} = \mathfrak{m}} f_\mathfrak{v} g_\mathfrak{w} \right) \mathfrak{m}.$$

In $C[[\mathfrak{M}]]$, there is also a natural notion of infinite sums: if $I$ is a set and $(f_i)_{i \in I}$ is a family of well-based series in $C[[\mathfrak{M}]]$, then we say that it is *summable* if $\bigcup_{i \in I} \operatorname{supp} f_i$ is well-based and $\{i \in I : f_{i, \mathfrak{m}} \neq 0\}$ is finite for every $\mathfrak{m} \in \mathfrak{M}$. In that case, we define the *sum* $f = \sum_{i \in I} f_i \in C[[\mathfrak{M}]]$ of this family by

$$f := \sum_{\mathfrak{m} \in \mathfrak{M}} \left( \sum_{i \in I} f_{i, \mathfrak{m}} \right) \mathfrak{m}.$$

Consider a second monomial group $\mathfrak{N}$ and a map $\varphi : C[[\mathfrak{M}]] \longrightarrow C[[\mathfrak{N}]]$. We say that $\varphi$ is *strongly linear* if it is $C$-linear and for every summable family $(f_i)_{i \in I}$ in $C[[\mathfrak{M}]]$, the family $(\varphi(f_i))_{i \in I}$ is summable in $C[[\mathfrak{N}]]_\lambda$ with $\varphi(\sum_{i \in I} f_i) = \sum_{i \in I} \varphi(f_i)$. By [28, Proposition 10], in order to show that a linear map $\varphi$ is strongly linear, it suffices to prove that the above condition holds for families of scalar multiples of monomials. So $\varphi$ is strongly linear if and only if for all $f \in C[[\mathfrak{M}]]$, the family $(f_\mathfrak{m} \varphi(\mathfrak{m}))_{\mathfrak{m} \in \operatorname{supp} f}$ is summable, with

$$\varphi(f) = \sum_{\mathfrak{m} \in \operatorname{supp} f} f_\mathfrak{m} \varphi(\mathfrak{m}).$$

Since the support of any $f \in C[[\mathfrak{M}]]$ is well-based, the order type $\operatorname{ot}(f) = \operatorname{ot}(\operatorname{supp} f)$ of $\operatorname{supp} f$ for the opposite order of $\leqslant$ is an ordinal. Now consider an $\varepsilon$-number $\lambda$. We recall that this means that $\dot{\omega}^\lambda = \lambda$, where $\dot{\omega}^\lambda$ stands for Cantor's $\lambda$-th ordinal power of $\omega$. It is known [23, Corollary 6.4] that the series $f \in C[[\mathfrak{M}]]$ with $\operatorname{ot}(f) < \lambda$ form a subfield $C[[\mathfrak{M}]]_\lambda$ of $C[[\mathfrak{M}]]$.

The ordering on **No** induces a natural valuation $v$ on **No** whose residue field is $\mathbb{R}$. The *Archimedean class* of a non-zero surreal number $x$ is the class $A_x$ of all $y \in \mathbf{No}$ with the same valuation as $x$. One of the discoveries of Conway was that $A_x \cap \mathbf{No}^>$ admits a simplest element that we will denote by $\mathfrak{d}_x$. Let $\mathbf{Mo} := \{\mathfrak{d}_x : x \in \mathbf{No}^{\neq}\}$ be the class of all $\mathfrak{d}_x$ that we may obtain in this way. Conway also constructed an order preserving bijection $\dot{\omega} : \mathbf{No} \longrightarrow \mathbf{Mo}; x \longmapsto \dot{\omega}^x$ that extends Cantor's ordinal exponentiation.

Through this $\dot{\omega}$-map and the so-called Conway normal form [14, Chapter 5], it turns out that the field **No** is naturally isomorphic to a field of well-based series $\mathbb{R}[[\mathbf{Mo}]]_{\mathbf{On}}$, for which **Mo** becomes the monomial group. For this series representation, any number $x \in \mathbf{No}$ has a set-sized support $\operatorname{supp} x$. The *Conway normal form* of $x$ coincides with its expression as a series $x = \sum_{\mathfrak{m} \in \mathbf{Mo}} x_\mathfrak{m} \mathfrak{m}$. For $x, y \in \mathbf{No}$ we sometimes write $x +\!\!\!+ y$ instead of $x + y$ in order to indicate that we have $\operatorname{supp} y \prec \operatorname{supp} x$, and thus that $x$ is a truncation of $x + y$ as a series.

# 3 Arithmetic on surreal numbers

In the sequel of this paper, by "number", we will always mean "surreal number".

## 3.1 Surreal arithmetic

We already explained the usefulness of Conway's cut representation for the recursive definition of functions on **No** and mentioned the addition (2.1) as an example. In fact, one may define all basic ring operations in a similar way:

$$
\begin{align}
0 &= \{\,|\,\} \tag{3.1}\\
1 &= \{0\,|\,\} \tag{3.2}\\
-x &= \{-x_R\,|\,-x_L\} \tag{3.3}\\
x+y &= \{x_L+y, x+y_L \,|\, x_R+y, x+y_R\} \tag{3.4}\\
xy &= \{x'y+xy'-x'y', x''y+xy''-x''y'' \,|\, x'y+xy''-x'y'', x''y+xy'-x''y'\} \\
&\quad (x' \in x_L,\ x'' \in x_R,\ y' \in y_L,\ y'' \in y_R). \tag{3.5}
\end{align}
$$

One major discovery of Conway was that the surreal numbers **No** actually form a real closed field for these operations and the ordering $\leqslant$. As an ordered field, it naturally contains the dyadic numbers, which are the numbers with finite length, and the real numbers, which are the numbers of length $\ell(r) \leqslant \omega$ whose sign sequence does not end with infinitely many consecutive identical signs.

The class **On** of ordinals is also naturally embedded into $(\mathbf{No}, \leqslant)$ by identifying an ordinal $\alpha$ with the constant sequence of length $\alpha$ with $\alpha[\beta] = 1$ for all $\beta < \alpha$. Thus, in **No**, expressions such as

$$\pi\sqrt{\omega_1} - \tfrac{2}{\omega}, \frac{\omega^{3/4}+1}{1-\omega^2}, \ldots$$

make sense and are amenable to various computations and comparisons. See [14, Chapter 1] for more details on the field operations on **No**. See [24, Chapters 1, 2 and 3] for more details on those operation in the framework of sign sequences and on the correspondence between cuts and sign sequences.

Using hints from Kruskal, Gonshor also defined an exponential function on **No**, which we denote by exp [24, Page 145]. This function extends the usual exponential function on $\mathbb{R}$. In fact, it turns out that **No** is an elementary extension of $\mathbb{R}$ as an ordered exponential field [15, Corollary 5.5]. In other words, the usual exponential function and its extended version to **No** satisfy the same first order properties over $\mathbb{R}$.

In order to define $\exp x$ for $x \in \mathbf{No}$ using a recursive equation, one needs to find an appropriate characterization of the cut formed by $\exp x$ inside the field generated by $x$, $x_\sqsubset$, and $\exp x_\sqsubset$. In exponential fields, the natural inequalities satisfied by such cuts involve truncated Taylor series expansions. Given $n \in \mathbb{N}$ and $a \in \mathbf{No}$, let

$$[a]_n = \sum_{k \leqslant n} \frac{a^k}{k!}.$$

If $x \in \mathbf{No}$ and $x' \in x_L$ is such that $\exp(x')$ is already defined, then for $n \in \mathbb{N}$, we should have

$$\exp(x) = \exp(x') \exp(x - x') > \exp(x') [x - x']_n$$

and one expects that such inequalities give sharp approximations of $\exp x$. Following this line of thought, Gonshor defined

$$\exp x = \left\{ 0, [x-x']_\mathbb{N} \exp x', [x-x'']_{2\mathbb{N}+1} \exp x'' \;\middle|\; \frac{\exp x''}{[x-x'']_{2\mathbb{N}+1}}, \frac{\exp x'}{[x'-x]_\mathbb{N}} \right\} \\ (x' \in x_L, x'' \in x_R). \tag{3.6}$$

The reciprocal of exp, defined on $\mathbf{No}^>$, is denoted log. This also leads to a natural powering operation: given $x \in \mathbf{No}^>$ and $y \in \mathbf{No}$, we define $x^y = \exp(y \log(x))$. Given $r \in \mathbb{R}$, we have $\dot{\omega}^r = \omega^r$, but for more general elements $x \in \mathbf{No}$, one does not necessarily have $\dot{\omega}^x = \omega^x$. (see [9] for more details).

## 3.2 Extending ordinal arithmetic

We write $\mathbf{On}^>$ and $\mathbf{On}_{\lim}$ for the classes of non-zero and limit ordinal numbers, respectively. The class of ordinal numbers is equipped with two distinct sets of operations: Cantor's (non-commutative) ordinal arithmetic and Hessenberg's (commutative) arithmetic. For ordinals $\alpha, \beta$, we will denote their ordinal sum, product, and exponentiation by $\alpha \dotplus \beta$, $\alpha \dot{\times} \beta$ and $\dot{\alpha}^\beta$. Their Hessenberg sum and product coincide with their sum and product when seen as surreal numbers [24, Theorems 4.5 and 4.6]; accordingly, we denote them by $\alpha + \beta$ and $\alpha \beta$. We assume that the reader is familiar with elementary computations in ordinal arithmetic. In this section, we define operations on surreal numbers which extend ordinal arithmetic.

For numbers $x, y$, we let $x \dotplus y$ denote the number, called the *concatenation sum* of $x$ and $y$, whose sign sequence is the concatenation of that of $y$ at the end of that of $x$. So $x \dotplus y$ is the number of length $\ell(x \dotplus y) = \ell(x) \dotplus \ell(y)$, which satisfies

$$\begin{aligned}(x \dotplus y)[\alpha] &= x[\alpha] & (\alpha < \ell(x)) \\ (x \dotplus y)[\ell(x) \dotplus \beta] &= y[\beta] & (\beta < \ell(y))\end{aligned}$$

It is easy to check that this extends the definition of ordinal sums. Moreover, the concatenation sum is associative and satisfies $\sup_{\sqsubseteq} (x \dotplus y_{\sqsubset}) = x \dotplus y$ whenever $x \in \mathbf{No}$ and $y \in \mathbf{No}$ is a limit number.

We let $x \dottimes y$ denote the number of length $\ell(x) \dottimes \ell(y)$, called the *concatenation product* of $x$ and $y$, whose sign sequence is defined by

$$(x \dottimes y)[\ell(x) \dottimes \alpha \dotplus \beta] = y[\alpha] x[\beta] \qquad (\alpha < \ell(y), \beta < \ell(x)).$$

Here we consider $y[\alpha] x[\beta]$ as a product in $\{-1, +1\}$. Informally speaking, given $x \in \mathbf{No}$ and $\alpha \in \mathbf{On}$, the number $x \dottimes \alpha$ is the $\alpha$-fold right-concatenation of $x$ with itself, whereas $\alpha \dottimes x$ is the number obtained from $x$ by replacing each sign $\alpha$ times by itself. We note that $\dottimes$ extends Cantor's ordinal product.

The operations $\dotplus$ and $\dottimes$ will be useful in what follows for the construction of simple yet interesting examples of surreal substructures. The remainder of this section is devoted to the collection of basic properties of these operations. The proofs can be skipped at a first reading, but we included them here for completeness and because we could not find them in the literature. We refer to [14, First Part] for a different extension of the ordinal product to the class of games (which properly contains $\mathbf{No}$).

**Lemma 3.1.** *For $x, y, z \in \mathbf{No}$, we have*

    *a)* $x \dottimes (y \dottimes z) = (x \dottimes y) \dottimes z$.

    *b)* $x \dottimes 1 = x$ *and* $x \dottimes (-1) = -x$.

    *c)* $x \dottimes (y \dotplus z) = (x \dottimes y) \dotplus (x \dottimes z)$.

    *d)* $x \dottimes y = \sup_{\sqsubseteq} (x \dottimes y_{\sqsubset})$ *if $y$ is limit.*

**Proof.** *a)* Both $x \dottimes (y \dottimes z)$ and $(x \dottimes y) \dottimes z$ have length $\ell(x) \dottimes \ell(y) \dottimes \ell(z)$. Let $\alpha < \ell(y \dottimes z)$ and $\delta < \ell(x)$. Write $\alpha = \ell(y) \dottimes \beta \dotplus \gamma$ where $\beta < \ell(z)$ and $\gamma < \ell(y)$. Then

$$\begin{aligned}(x \dottimes (y \dottimes z))[\ell(x) \dottimes \alpha \dotplus \delta] &= (y \dottimes z)[\alpha] x[\delta] \\ &= z[\beta] y[\gamma] x[\delta] \\ &= z[\beta] (x \dottimes y)[\ell(x) \dottimes \gamma \dotplus \delta] \\ &= ((x \dottimes y) \dottimes z)[\ell(x) \dottimes \ell(y) \dottimes \beta \dotplus \ell(x) \dottimes \gamma \dotplus \delta] \\ &= ((x \dottimes y) \dottimes z)[\ell(x) \dottimes \alpha \dotplus \delta].\end{aligned}$$

*b)* The numbers $x \dottimes 1$ and $x \dottimes (-1)$ have length $\ell(x) \dottimes 1 = \ell(x)$. For $\beta < \ell(x)$, we have $(x \dottimes 1)[\beta] = 1[0] x[\beta] = x[\beta]$ and $(x \dottimes (-1))[\beta] = (-1)[0] x[\beta] = -x[\beta]$.

*c)* The number $x \dottimes (y \dotplus z)$ has length

$$\begin{aligned}\ell(x) \dottimes \ell(y \dotplus z) &= \ell(x) \dottimes (\ell(y) \dotplus \ell(z)) \\ &= \ell(x) \dottimes \ell(y) \dotplus \ell(x) \dottimes \ell(z) \\ &= \ell(x \dottimes y) \dotplus \ell(x \dottimes z) \\ &= \ell((x \dottimes y) \dotplus (x \dottimes z)).\end{aligned}$$

Let $\beta < \ell(x)$ and $\alpha < \ell(y \dotplus z)$. If $\alpha < \ell(y)$, then

$$\begin{aligned}(x \dottimes (y \dotplus z))(\ell(x) \dottimes \alpha \dotplus \beta) &= (y \dotplus z)[\alpha] x[\beta] \\ &= y[\alpha] x[\beta] \\ &= (x \dottimes y)[\ell(x) \dottimes \alpha \dotplus \beta] \\ &= ((x \dottimes y) \dotplus (x \dottimes z))[\ell(x) \dottimes \alpha \dotplus \beta].\end{aligned}$$

Otherwise, there is $\eta < \ell(z)$ such that $\alpha = \ell(y) \dotplus \eta$ and then

$$\begin{aligned} x \dottimes (y \dotplus z)[\ell(x) \dottimes \alpha \dotplus \beta] &= (y \dotplus z)[\alpha] x[\beta] \\ &= z[\eta] x[\beta] \\ &= (x \dottimes z)[\ell(x) \dottimes \eta \dotplus \beta] \\ &= ((x \dottimes y) \dotplus (x \dottimes z))[\ell(x) \dottimes \alpha \dotplus \beta]. \end{aligned}$$

*d*) The previous identities imply in particular that $x \dottimes y_\sqsubset$ is linearly ordered by simplicity, which means that the supremum $\sup_\sqsubseteq (x \dottimes y_\sqsubset)$ is well defined in $(\mathbf{No}, \sqsubseteq)$. Assume $y$ is limit. If $y = 0$, then we have $x \dottimes y = 0 = \sup_\sqsubseteq x \dottimes 0_\sqsubset$. Assume $y \neq 0$. Notice that we have $\ell(y) = \sup_\sqsubseteq \ell(y_\sqsubset)$, so

$$\ell(x \dottimes y) = \ell(x) \dottimes \sup_\sqsubseteq \ell(y_\sqsubset) = \sup_\sqsubseteq (\ell(x) \dottimes \ell(y_\sqsubset)) = \sup_\sqsubseteq \ell(x \dottimes y_\sqsubset).$$

Let $\beta < \ell(x)$ and $\alpha < \ell(y)$. Since $y$ is a limit number, there is $u \in y_\sqsubset$ such that $\alpha < \ell(u)$. Then

$$(x \dottimes y)[\ell(x) \dottimes \alpha \dotplus \beta] = y[\alpha] x[\beta] = u[\alpha] x[\beta] = (x \dottimes u)[\ell(x) \dottimes \alpha \dotplus \beta]. \qquad \square$$

**Remark 3.2.** The previous lemma can be regarded as an alternative way to define the concatenation product. Yet another way is through the equation

$$\forall x > 0, \forall y, \; x \dottimes y \;=\; \{x \dottimes y_L + x_L, x \dottimes y_R \dotplus (-x_R) \,|\, x \dottimes y_L \dotplus x_R, x \dottimes y_R \dotplus (-x_L)\}. \qquad (3.7)$$

Likewise, the contatenation sum has the following equation [18, Proposition 2]:

$$\forall x, \forall y, \; x \dotplus y \;=\; \{x_L, x \dotplus y_L \,|\, x \dotplus y_R, x_R\}. \qquad (3.8)$$

Note that these two equations are *not* uniform in the sense of Definition 4.29 below.

**Proposition 3.3.** *Let $x, y, z \in \mathbf{No}$.*

  *a*) *If $x \neq 0$, then $y \sqsubseteq z$ if and only if $x \dottimes y \sqsubseteq x \dottimes z$.*
  *b*) *If $0 < x$, then $y < z$ if and only if $x \dottimes y < x \dottimes z$.*

**Proof.** *a*) If $y \sqsubseteq z$, then for $a \in \mathbf{No}$ with $z = y \dotplus a$, Lemma 3.1(c) implies that

$$x \dottimes y \sqsubseteq (x \dottimes y) \dotplus (x \dottimes a) = x \dottimes z.$$

Conversely, if $x \dottimes y \sqsubseteq x \dottimes z$, then since $x \neq 0$, we may compute, for $\alpha < \ell(y)$, the sign $y[\alpha] x[0] = (x \dottimes y)[\ell(x) \dottimes \alpha] = (x \dottimes z)[\ell(x) \dottimes \alpha] = z[\alpha] x[0]$. We deduce that $y[\alpha] = z[\alpha]$, so $y \sqsubseteq z$.

*b*) If $y < z$, then given the maximal common initial segment $u$ of $y$ and $z$, we have $(x \dottimes u) \sqsubseteq (x \dottimes y), (x \dottimes z)$, with $\ell(x \dottimes u) = \ell(x) \dottimes \ell(u)$. Thus $(x \dottimes y)[\ell(x) \dottimes \ell(u)] = y[\ell(u)] x[0] = y[\ell(u)]$ is strictly smaller than $z[\ell(u)] = z[\ell(u)] x[0] = (x \dottimes z)[\ell(x) \dottimes \ell(u)]$, which means that $x \dottimes y < x \dottimes z$. Since the order $\leqslant$ is linear, this suffices to prove the result. $\square$

# 4 Surreal substructures

## 4.1 Surreal substructures and their parameterizations

Let $\mathbf{X}$ be a subclass of $\mathbf{No}$ and let $\mathcal{R} = (\leqslant_i)_{i \in I}$ be a family of ordering relations on $\mathbf{No}$. Then we say that a function $f: \mathbf{X} \longrightarrow \mathbf{No}$ is $\mathcal{R}$-*increasing* if $f$ is increasing for each $\leqslant_i$ with $i \in I$. If $f$ is also injective, then we say that it is *strictly $\mathcal{R}$-increasing*. If we have $x \leqslant_i y \Longleftrightarrow f(x) \leqslant_i f(y)$ for all $x, y \in \mathbf{X}$ and $i \in I$, then we call $f$ an $\mathcal{R}$-*embedding* of $(\mathbf{X}, (\leqslant_i)_{i \in I})$ into $(\mathbf{No}, (\leqslant_i)_{i \in I})$. We simply say that $f$ is an *embedding* if $f$ is a $(\leqslant, \sqsubseteq)$-embedding.

**Definition 4.1.** *A* **surreal substructure** *is the image of an embedding of* **No** *into itself.*

**Example 4.2.** Given $a \in \mathbf{No}$, the map $x \longmapsto a \dotplus x$ is an embedding of $(\mathbf{No}, \leqslant, \sqsubseteq)$ into itself. If $a > 0$, then so is the map $x \longmapsto a \dottimes x$, by Proposition 3.3. Consequently:

- For $a \in \mathbf{No}$, the map $x \longmapsto a \dotplus x$ gives rise to the surreal substructure $a \dotplus \mathbf{No}$ of numbers whose sign sequences begin with the sign sequence of $a$.

- For $0 < a \in \mathbf{No}$, the map $x \longmapsto a \dottimes x$ induces the surreal substructure $a \dottimes \mathbf{No}$ of numbers whose sign sequences are (possibly empty or transfinite) concatenations of the sign sequences of $a$ and $-a$.

**Example 4.3.** Let $\varphi$ be an embedding of **No** into itself with image **S**. Then the map $\psi: x \longmapsto -\varphi(-x)$ defines another embedding of **No** into itself with image $-\mathbf{S} = \{-x : x \in \mathbf{S}\}$. In other words, if **S** is a surreal substructure, then so is $-\mathbf{S}$.

We claim that any strictly $(\leqslant, \sqsubseteq)$-increasing map $f: \mathbf{No} \longrightarrow \mathbf{No}$ is automatically an embedding. We first need a lemma.

**Lemma 4.4.** *If $x, y, z$ are numbers such that $x \sqsubseteq y$ and $x \not\sqsubseteq z$, then we have $x < z$ if and only if $y < z$, and $z < x$ if and only if $z < y$.*

**Proof.** Since $x \not\sqsubseteq z$, we have $x < z$ if and only if there is $\eta_x < \ell(x)$ with $x \upharpoonright \eta_x = z \upharpoonright \eta_x$ and $x[\eta_x] < z[\eta_x]$. Now $x \sqsubseteq y$ so $y \not\sqsubseteq z$ and likewise $y < z$ holds if and only if there is $\eta_y < \ell(y)$ with $y \upharpoonright \eta_y = z \upharpoonright \eta_y$ and $y[\eta_y] < z[\eta_y]$. Notice that $y \upharpoonright \eta_y = z \upharpoonright \eta_y$ and $y \sqsupseteq x \not\sqsubseteq z$ imply that $\eta_y < \ell(x)$. In both cases, since $x \sqsubseteq y$, we have $x[\eta_x] = y[\eta_x]$ and $x[\eta_y] = y[\eta_y]$. Therefore the existence of $\eta_x$ yields that of $\eta_y = \eta_x$ and *vice versa*. The other equivalence follows by symmetry. □

**Lemma 4.5.** *Assume that* **X** *is a convex subclass of* $(\mathbf{No}, \leqslant)$. *Then every strictly $(\leqslant, \sqsubseteq)$-increasing function $\varphi: \mathbf{X} \longrightarrow \mathbf{No}$ is an embedding $(\mathbf{X}, \leqslant, \sqsubseteq) \longrightarrow (\mathbf{No}, \leqslant, \sqsubseteq)$.*

**Proof.** Since $(\mathbf{No}, \leqslant)$ is a linear order, the function $\varphi$ is automatically an embedding for $\leqslant$, so we need only prove that it is an embedding for $\sqsubseteq$. Assume for contradiction that there are elements $x < y$ of **X** such that $x \not\sqsubseteq y$ and $\varphi(x) \sqsubseteq \varphi(y)$. Let $z$ be the $\sqsubseteq$-maximal common initial segment of $x$ and $y$. We have $x < z \leqslant y$, so $z \in \mathbf{X}$. Since $\varphi$ is strictly $(\leqslant, \sqsubseteq)$-increasing, we have $\varphi(x) < \varphi(z) \leqslant \varphi(y)$ and $\varphi(x) \not\sqsubseteq \varphi(z)$, which given our assumption $\varphi(x) \sqsubseteq \varphi(y)$ contradicts the previous lemma. Hence $\varphi(x) \not\sqsubseteq \varphi(y)$, which concludes the proof. □

Since a surreal substructure **S** is an isomorphic copy of **No** into itself, it should induce a natural Conway bracket $\{\}_\mathbf{S}$ on **S**. This actually leads to an equivalent definition of surreal substructures. Let us investigate this in more detail.

Let **S** be an arbitrary subclass of **No**. We say that **S** is *rooted* if it admits a simplest element, called its *root*, and which we denote by $\mathbf{S}^\bullet$. Given subclasses $\mathbf{L} < \mathbf{R}$ of **S**, we let $(\mathbf{L} | \mathbf{R})_\mathbf{S}$ denote the class of elements $x \in \mathbf{S}$ such that $\mathbf{L} < x < \mathbf{R}$. If $(\mathbf{L} | \mathbf{R})_\mathbf{S}$ is rooted, then we let $\{\mathbf{L} | \mathbf{R}\}_\mathbf{S}$ denote its root. If $L = \mathbf{L}$ and $R = \mathbf{R}$ are sets, then we call $(L | R)_\mathbf{S}$ the *cut* in **S** defined by $L$ and $R$. If for any subsets $L < R$ of **S** the class $(L | R)_\mathbf{S}$ is rooted, then we say that **S** *admits an induced Conway bracket*.

**Proposition 4.6.** *Let* **S** *admit an induced Conway bracket. Then the map* $\Xi_\mathbf{S}\colon \mathbf{No}\longrightarrow \mathbf{S}$ *defined by*

$$\forall x\in \mathbf{No}, \Xi_\mathbf{S}\, x = \{\Xi_\mathbf{S}\, x_L\,|\,\Xi_\mathbf{S}\, x_R\}_\mathbf{S}$$

*is an isomorphism* $(\mathbf{No},\leqslant,\sqsubseteq)\longrightarrow (\mathbf{S},\leqslant,\sqsubseteq)$.

**Proof.** We first justify that $\Xi_\mathbf{S}$ is well defined. Let $x\in \mathbf{No}$ be such that $\Xi_\mathbf{S}$ is well-defined and strictly $\leqslant$-increasing on $x_\sqsubset$, with values in **S**. We have $\Xi_\mathbf{S}\, x_L < \Xi_\mathbf{S}\, x_R$ where those sets are in **S** so $\Xi_\mathbf{S}\, x$ is a well-defined element of $(\Xi_\mathbf{S}\, x_L\,|\,\Xi_\mathbf{S}\, x_R)_\mathbf{S}$, and $\Xi_\mathbf{S}$ is strictly $\leqslant$-increasing on $\{x\}\cup x_L\cup x_R$. By induction, $\Xi_\mathbf{S}$ is a strictly increasing map $\mathbf{No}\longrightarrow \mathbf{S}$. Let $y\in \mathbf{No}$ with $x\sqsubseteq y$, so that $x_L < y < x_R$. By definition, the number $\Xi_\mathbf{S}\, x$ is the simplest element $u\in \mathbf{S}$ with $\Xi_\mathbf{S}\, x_L < u < \Xi_\mathbf{S}\, x_R$. Since $\Xi_\mathbf{S}\, y\in \mathbf{S}$ and $\Xi_\mathbf{S}\, x_L < \Xi_\mathbf{S}\, y < \Xi_\mathbf{S}\, y_L$, it follows that $\Xi_\mathbf{S}\, x\sqsubseteq \Xi_\mathbf{S}\, y$. We deduce from Lemma 4.5 that $\Xi_\mathbf{S}$ is an embedding of $(\mathbf{No},\leqslant,\sqsubseteq)$ into itself.

We now prove that $\mathbf{S}=\Xi_\mathbf{S}\,\mathbf{No}$ by induction on $y\in \mathbf{S}$ for $\sqsubseteq$. Let $y\in \mathbf{S}$ be such that $y_\sqsubset\cap \mathbf{S}$ is a subset of $\Xi_\mathbf{S}\,\mathbf{No}$. Let $\Xi_\mathbf{S}\, L' = L = y_L\cap \mathbf{S}$ and $R = y_R\cap \mathbf{S} = \Xi_\mathbf{S}\, R'$ where since $\Xi_\mathbf{S}$ is strictly $\leqslant$-increasing and thus injective, the sets $L', R'$ are uniquely determined and satisfy $L' < R'$. Since **S** admits an induced Conway bracket, the cut $(L\,|\,R)_\mathbf{S}$ is rooted and contains $y$, so $\{L\,|\,R\}_\mathbf{S}\sqsubseteq y$. Since $\{L\,|\,R\}_\mathbf{S}\notin L\cup R$, we necessarily have $y=\{L\,|\,R\}_\mathbf{S} = \Xi_\mathbf{S}\,\{L'\,|\,R'\}$. By induction, we conclude that $\mathbf{S} = \Xi_\mathbf{S}\,\mathbf{No}$. □

**Proposition 4.7.** *Let* **S** *be a subclass of* **No**. *Then* **S** *is a surreal substructure if and only if it admits an induced Conway bracket.*

**Proof.** Assume that **S** admits an induced Conway bracket. By the previous proposition, **S** is the range of the strictly $(\leqslant,\sqsubseteq)$-increasing function $\Xi_\mathbf{S}\colon \mathbf{No}\longrightarrow \mathbf{No}$, whence **S** is a surreal substructure. Conversely, consider an embedding $\varphi$ of **No** into itself with image **S**. Let $L < R$ be subsets of **S** and define $(L', R') = (\varphi^{-1}(L), \varphi^{-1}(R))$. The function $\varphi$ is strictly $\leqslant$-increasing so $L' < R'$, and we may consider the number $x = \{L'\,|\,R'\}$. Now let $y\in (L\,|\,R)_\mathbf{S}$. We have $\varphi^{-1}(y)\in (L'\,|\,R')$, so $x\sqsubseteq \varphi^{-1}(y)$. Since $\varphi$ is $\sqsubseteq$-increasing, this implies $\varphi(x)\sqsubseteq y$, which proves that $\varphi(x) = \{L\,|\,R\}_\mathbf{S}$, so **S** admits an induced Conway bracket. □

**Remark 4.8.** More generally, one may discard the existence condition for the Conway bracket and consider subclasses **X** of **No** that satisfy the following condition:

**IN.** For all subsets $L, R$ of **X** with $L < R$, the class $(L\,|\,R)_\mathbf{X}$ is either empty or rooted.

A subclass $\mathbf{X}\subseteq \mathbf{No}$ satisfies **IN** if and only if there is a (unique) $\sqsubseteq$-initial subclass $\mathbf{I_S}$ of **No** and a (unique) isomorphism $(\mathbf{I_S},\leqslant,\sqsubseteq)\longrightarrow (\mathbf{S},\leqslant,\sqsubseteq)$. This is in particular the case for the classes $\mathbf{Smp}_\Pi$ described in Section 6 below. For more details on this more general kind of subclasses, we refer to [19].

In this paper, we focus on surreal substructures. The characterizations given in Proposition 4.7 and Proposition 4.13 are known results. The second one was first proved (for more general types of ordinal sequences) by Lurie [34, Theorem 8.3], and both of them were proved by Ehrlich [19, Theorems 1 and 4].

**Proposition 4.9.** *Let* **S** *be a surreal substructure. The function* $\Xi_\mathbf{S}$ *is the unique surjective strictly* $(\leqslant,\sqsubseteq)$-*increasing function* $\mathbf{No}\longrightarrow \mathbf{S}$.

**Proof.** Let $\varphi$ be a strictly $(\leqslant,\sqsubseteq)$-increasing function $\mathbf{No} \longrightarrow \mathbf{S}$ with image $\mathbf{S}$. By Lemma 4.5, it is an embedding. Given $x \in \mathbf{No}$ such that $\varphi$ and $\Xi_\mathbf{S}$ coincide on $x_\sqsubset$, the numbers $\varphi(x)$ and $\Xi_\mathbf{S} x$ of $\mathbf{S}$ are both the simplest element of $(\Xi_\mathbf{S} x_L | \Xi_\mathbf{S} x_R)_\mathbf{S}$ and are thus equal. It follows by induction that $\varphi = \Xi_\mathbf{S}$. □

**Lemma 4.10.** *Let $\mathbf{S}$ be a surreal substructure. For $x \in \mathbf{No}$, we have $\ell(x) \leqslant \ell(\Xi_\mathbf{S} x)$.*

**Proof.** By Proposition 4.6, the map $\Xi_\mathbf{S}$ realizes an embedding of $(x_\sqsubset, \sqsubseteq)$ into $((\Xi_\mathbf{S} x)_\sqsubset, \sqsubseteq)$, so the order type $\ell(x)$ of the former is smaller than that of the latter, namely $\ell(\Xi_\mathbf{S} x)$. □

Given a surreal substructure $\mathbf{S}$, we call $\Xi_\mathbf{S}$ the *defining surreal isomorphism* of *parametrization* of $\mathbf{S}$. The above uniqueness property is fundamental; it allows us in particular to perform constructions on surreal substructures *via* their defining surreal isomorphisms and *vice versa*.

## 4.2 Cut representations

Let $\mathbf{S}$ be a surreal substructure. Given an element $x \in \mathbf{S}$ and subsets $L, R$ of $\mathbf{S}$ with $L < R$, we say that $(L, R)$ is a *cut representation* of $x$ in $\mathbf{S}$ if $x = \{L | R\}_\mathbf{S}$. We refer to elements in $L$ and $R$ as *left* and *right options* of the representation. For $x \in \mathbf{S}$, we write

$$(x_L^\mathbf{S}, x_R^\mathbf{S}) := (x_L \cap \mathbf{S}, x_R \cap \mathbf{S})$$

and call this pair the *canonical representation* of $x$ in $\mathbf{S}$. We also write $x_\sqsubset^\mathbf{S}$ for the set $x_\sqsubset \cap \mathbf{S}$.

A $\sqsubseteq$-*final substructure* of $\mathbf{S}$ is a rooted final segment $\mathbf{T}$ of $\mathbf{S}$ for $\sqsubseteq$ (and thereby necessarily a substructure). It is easy to see that this is the case if and only if $\mathbf{T}$ is rooted and $\mathbf{T}$ is the class $\mathbf{S}^{\sqsupseteq \mathbf{T}^\bullet}$ of elements $x \in \mathbf{S}$ such that $\mathbf{T}^\bullet \sqsubseteq x$.

**Proposition 4.11.** *Let $\mathbf{S}$ be a surreal substructure and let $(L, R)$ and $(L', R')$ be cut representations in $\mathbf{S}$. For $x \in \mathbf{S}$, we have*

*a)* $\{L | R\}_\mathbf{S} \leqslant \{L' | R'\}_\mathbf{S}$ *if and only if* $\{L | R\}_\mathbf{S} < R'$ *and* $L < \{L' | R'\}_\mathbf{S}$.

*b)* $(x_L^\mathbf{S}, x_R^\mathbf{S})$ *is a cut representation of $x$ in $\mathbf{S}$ with respect to which any other cut representation of $x$ in $\mathbf{S}$ is cofinal.*

*c)* $\mathbf{S}^{\sqsupseteq x} = (x_L^\mathbf{S} | x_R^\mathbf{S})_\mathbf{S}$.

**Proof.** The assertions *a)* and *b)* are true when $\mathbf{S} = \mathbf{No}$ by [24, Theorems 2.5 and 2.9]. By Proposition 4.6, the function $\Xi_\mathbf{S}$ is an isomorphism $(\mathbf{No}, \leqslant, \sqsubseteq) \longrightarrow (\mathbf{S}, \leqslant, \sqsubseteq)$, satisfying the relation $\forall a \in \mathbf{No}, (\Xi_\mathbf{S} a_L, \Xi_\mathbf{S} a_R) = ((\Xi_\mathbf{S} a)_L^\mathbf{S}, (\Xi_\mathbf{S} a)_R^\mathbf{S})$, so *a)* and *b)* hold in general. We have $\mathbf{S}^{\sqsupseteq x} \supseteq (x_L^\mathbf{S} | x_R^\mathbf{S})_\mathbf{S}$, since $x = (x_L^\mathbf{S} | x_R^\mathbf{S})_\mathbf{S}^\bullet$. Conversely, for $y \in \mathbf{S}^{\sqsupseteq x}$ and $x' \in x_\sqsubset^\mathbf{S}$, we have $x' \sqsubset y$ and $y[\ell(x')] = x[\ell(x')] \in \{-1, 1\}$, so $y - x'$ and $x - x'$ have the same sign. We conclude that $x_L^\mathbf{S} < y < x_R^\mathbf{S}$, which completes the proof of *c)*. □

## 4.3 Imbrications

Let $\mathbf{S}, \mathbf{T}$ be two surreal substructures. Then there is a unique $(\leqslant, \sqsubseteq)$-isomorphism $\Xi_\mathbf{T}^\mathbf{S} := \Xi_\mathbf{T} \Xi_\mathbf{S}^{-1} : \mathbf{S} \longrightarrow \mathbf{T}$ that we call the *surreal isomorphism* between $\mathbf{S}$ and $\mathbf{T}$. The composition $\Xi_\mathbf{S} \circ \Xi_\mathbf{T}$ is also an embedding, so its image $\mathbf{S} \prec \mathbf{T} := \Xi_\mathbf{S} \mathbf{T}$ is again a surreal substructure that we call the *imbrication* of $\mathbf{T}$ into $\mathbf{S}$. We say that $\mathbf{T}$ is a *left factor* (resp. *right factor*) of $\mathbf{S}$ if there is a surreal substructure $\mathbf{U}$ such that $\mathbf{S} = \mathbf{T} \prec \mathbf{U}$ (resp. $\mathbf{S} = \mathbf{U} \prec \mathbf{T}$).

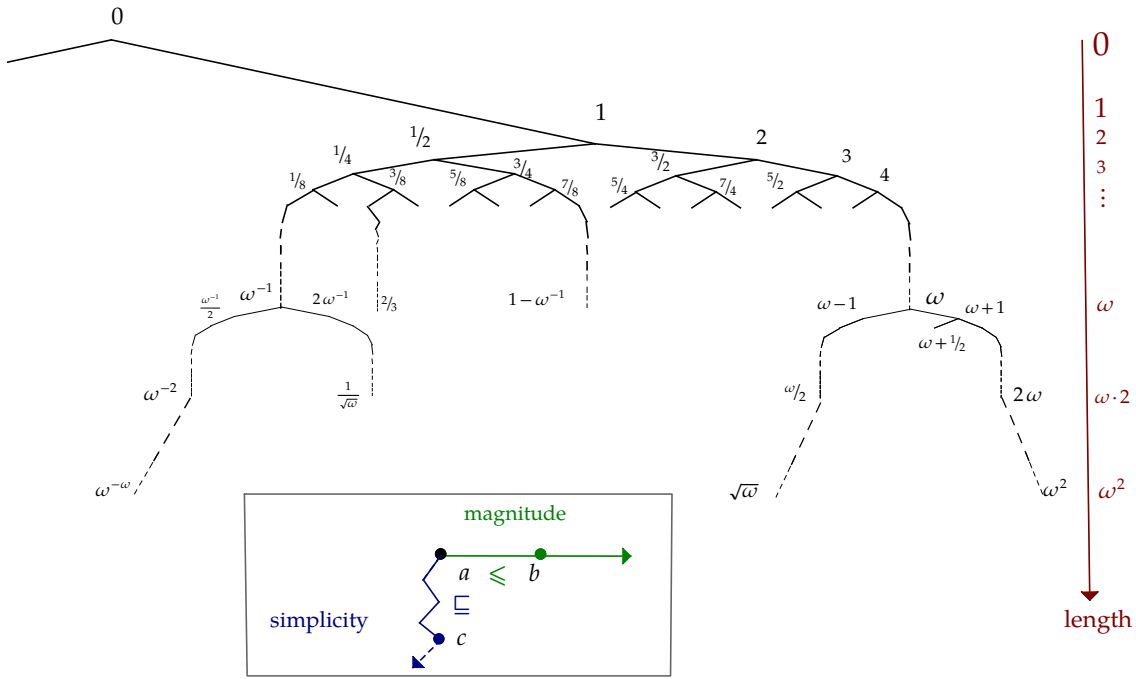

**Figure 4.1.** The class of positive surreal numbers as a tree. For clarity, only a few numbers up to the length $\omega^2$ are represented. Negative numbers are obtained through symmetry w.r.t. the $y$-axis.

By the associativity of the composition of functions, the imbrication of surreal substructures is associative. Right factors are determined by the two other substructures. More precisely, since $\Xi_\mathbf{T}$ is injective, the relation $\mathbf{S} = \mathbf{T} \prec \mathbf{U} = \Xi_\mathbf{T} \mathbf{U}$ yields $\mathbf{U} = \Xi_\mathbf{T}^{-1}(\mathbf{S})$. The same does not hold for left factors:

$$(1 \,\dotplus\, \mathbf{No}) \,\dotplus\, (\omega \,\dotplus\, \mathbf{No}) = \mathbf{No} \prec (\omega \,\dotplus\, \mathbf{No}) = \omega \,\dotplus\, \mathbf{No}.$$

**Proposition 4.12.** *If $\mathbf{S}, \mathbf{T}$ are surreal substructures, then $\mathbf{T}$ is a left factor of $\mathbf{S}$ if and only if $\mathbf{S} \subseteq \mathbf{T}$.*

**Proof.** If $\mathbf{S} = \mathbf{T} \prec \mathbf{U}$, then $\mathbf{S} = \Xi_\mathbf{T} \mathbf{S} \subseteq \mathbf{T}$. Assume that $\mathbf{S} \subseteq \mathbf{T}$ and let $\mathbf{U} = \Xi_\mathbf{T}^{-1}(\mathbf{S})$. We have $\mathbf{U} = (\Xi_\mathbf{T}^{-1} \!\restriction\! \mathbf{S})\, \Xi_\mathbf{S}\, \mathbf{No}$ where $\Xi_\mathbf{T}^{-1} \!\restriction\! \mathbf{S}$ and $\Xi_\mathbf{S}$ are respectively embeddings $(\mathbf{S}, \leqslant, \sqsubseteq) \longrightarrow (\mathbf{No}, \leqslant, \sqsubseteq)$ and $(\mathbf{No}, \leqslant, \sqsubseteq) \longrightarrow (\mathbf{S}, \leqslant, \sqsubseteq)$ so $(\Xi_\mathbf{T}^{-1} \!\restriction\! \mathbf{S})\, \Xi_\mathbf{S}$ is an embedding $(\mathbf{No}, \leqslant, \sqsubseteq) \longrightarrow (\mathbf{No}, \leqslant, \sqsubseteq)$. Hence $\mathbf{U}$ is a surreal substructure with $\Xi_\mathbf{T} \mathbf{U} = \mathbf{S}$, which means that $\mathbf{T} \prec \mathbf{U} = \mathbf{S}$. □

## 4.4 Surreal substructures as trees

Through the identification $\mathbf{No} \approx \{-1, 1\}^{<\mathbf{On}}$, the class of surreal numbers can naturally be represented by a full binary tree of uniform depth $\mathbf{On}$, as illustrated in Figure 4.1.

For each ordinal $\alpha$, we let $\mathbf{No}(\alpha)$ denote the subtree of $\mathbf{No}$ of nodes of depth $<\alpha$, that is, the set of numbers $x$ with $\ell(x) < \alpha$. This can be represented as the subtree obtained by cropping the picture at depth $\alpha$. In order to characterize surreal substructures in tree-theoretic terms, we need to investigate chains for $\sqsubseteq$: given a subclass $\mathbf{X} \subseteq \mathbf{No}$, a $\sqsubseteq$-*chain* in $\mathbf{X}$ is a linearly ordered (and thus well-ordered) *subset* $C$ of $(\mathbf{X}, \sqsubseteq)$. If a $\sqsubseteq$-chain $C$ in $(\mathbf{X}, \sqsubseteq)$ admits a supremum in $(\mathbf{X}, \sqsubseteq)$, we denote it $\sup_{\mathbf{X}, \sqsubseteq} C$. Note that the empty set has a supremum in $(\mathbf{X}, \sqsubseteq)$ if and only if $\mathbf{X}$ has a root, in which case $\sup_{\mathbf{X}, \sqsubseteq} \emptyset = \mathbf{X}^\bullet$. We say that $y \in \mathbf{X}$ is the *left successor* of $x \in \mathbf{X}$ if $y < x$ and $z \sqsupseteq y$ for every $z < x$ in $\mathbf{X}$. Right successors are defined similarly.

**Proposition 4.13.** *Let* **S** *be a class of surreal numbers. Then the following assertions are equivalent:*

*a)* **S** *is a surreal substructure.*

*b)* *Every element of* **S** *has a left and a right successor in* **S** *and every $\sqsubseteq$-chain in* **S** *has a supremum in* $(\mathbf{S}, \sqsubseteq)$.

**Proof.** Let **S** be a surreal substructure. In **No**, any element $x$ clearly admits a left successor $\{x_L | x\}$ and a right successor $\{x | x_R\}$, and every $\sqsubseteq$-chain clearly admits a supremum. Since these properties are preserved by the isomorphism $\Xi_\mathbf{S}$, we deduce *b)*.

Assume now that *b)* holds. We derive *a)* by inductively defining an isomorphism $\Xi: (\mathbf{No}, \sqsubseteq, \leqslant) \longrightarrow (\mathbf{S}, \sqsubseteq, \leqslant)$. Applying *b)* to the empty chain, we note that the supremum of $\emptyset$ in $(\mathbf{S}, \sqsubseteq)$ is the minimum of **S** for $\sqsubseteq$. So **S** is rooted and we may define $\Xi 0 = \mathbf{S}^\bullet$. Let $0 < \alpha$ be an ordinal such that $\Xi$ is defined and strictly $(\leqslant, \sqsubseteq)$-increasing on $\mathbf{No}(\alpha)$. We distinguish two cases:

- If $\alpha$ is limit, then let $x$ be a surreal number with length $\alpha$. Thus $x$ is a limit number and $\Xi x_\sqsubset$ is a $\sqsubseteq$-chain in **S**. We define $\Xi x = \sup_{\mathbf{X}, \sqsubseteq} \Xi x_\sqsubset$.

- Assume now that $\alpha$ is successor, let $x$ be a number with length $\alpha$, and write $x = u \dotplus \sigma$ where $\sigma \in \{-1, 1\}$. Let $u_{-1}$ and $u_1$ be the left and right successors of $\Xi u$. Then we define $\Xi x = u_\sigma$.

In both cases, this defines $\Xi$ on $\mathbf{No}(\alpha + 1)$ and the extension is clearly strictly $\sqsubseteq$-increasing and strictly $\leqslant$-increasing on every set $x_\sqsubseteq := \{x\} \cup x_\sqsubset$ for $x \in \mathbf{No}(\alpha + 1)$.

It remains to be shown that $\Xi$ is strictly $\leqslant$-increasing on $\mathbf{No}(\alpha + 1)$. Given $a < b$ in $\mathbf{No}(\alpha + 1)$, let $c \in \mathbf{No}(\alpha)$ be their $\sqsubseteq$-maximal common initial segment. We either have $a \leqslant c < b$ and thus $\Xi a \leqslant \Xi c < \Xi b$, or $a < c \leqslant b$ and thus $\Xi a < \Xi c \leqslant \Xi b$. So $\Xi$ is strictly $\leqslant$-increasing on $\mathbf{No}(\alpha + 1)$.

By induction, the function $\Xi$ is defined and $(\leqslant, \sqsubseteq)$-increasing on $\mathbf{No} = \bigcup_{\alpha \in \mathbf{On}} \mathbf{No}(\alpha)$. Note that $(\mathbf{S}, \sqsubseteq)$ is well-founded since $(\mathbf{No}, \sqsubseteq)$ is well-founded and $\mathbf{S} \subseteq \mathbf{No}$. By induction over $y \in \mathbf{S}$, let us show that $y$ lies in the range of $\Xi$. If $y$ is the left or right successor of an element $v \in \mathbf{S}$, then the induction hypothesis implies the existence of some $u \in \mathbf{No}$ with $v = \Xi u$, and we get $y = \Xi(u \dotpm 1)$. Otherwise, we have $y = \sup_\sqsubseteq y_\sqsubset^\mathbf{S} = \Xi \sup_\sqsubseteq C$ where $C = \{x \in \mathbf{No} : \Xi x \sqsubset y\}$. We conclude that $\Xi$ is an isomorphism. □

**Example 4.14.** Consider the class **Inc** defined by $\Xi_{\mathbf{Inc}} 0 := 1$, $\Xi_{\mathbf{Inc}}(u \dotplus \sigma) = (\Xi_{\mathbf{Inc}} u) \dotplus \sigma \dotplus 1$, for all $u \in \mathbf{No}$ and $\sigma \in \{-1, 1\}$ and $\Xi_{\mathbf{Inc}} \sup_\sqsubseteq C = (\sup_\sqsubseteq \Xi_{\mathbf{Inc}} C) \dotplus 1$ for every non-empty $\sqsubseteq$-chain $C$ without maximum in $(\mathbf{No}, \sqsubseteq)$. It is easy to check that we have $\ell(\Xi_{\mathbf{Inc}} x) > \ell(x)$ for every surreal number $x$.

**Example 4.15.** Let $\mathbf{S} = \mathbf{No}^\geqslant \setminus \{1\}$. Then $(\mathbf{S}, \sqsubseteq)$ is isomorphic to $(\mathbf{No}, \sqsubseteq)$, but **S** is not a surreal structure. In other words, the condition *b)* cannot be replaced by the weaker condition that $(\mathbf{S}, \sqsubseteq)$ and $(\mathbf{No}, \sqsubseteq)$ be isomorphic.

The characterization *b)* gives us some freedom in constructing a surreal substructure: one only has to provide a mechanism for chosing left and right successors of already constructed elements, as well as least upper bounds for already constructed branches (i.e. $\sqsubseteq$-chains). Intuitively speaking, this corresponds to a way to "draw" **S** as a full binary tree inside the binary tree that represents **No**: see Figure 4.2.

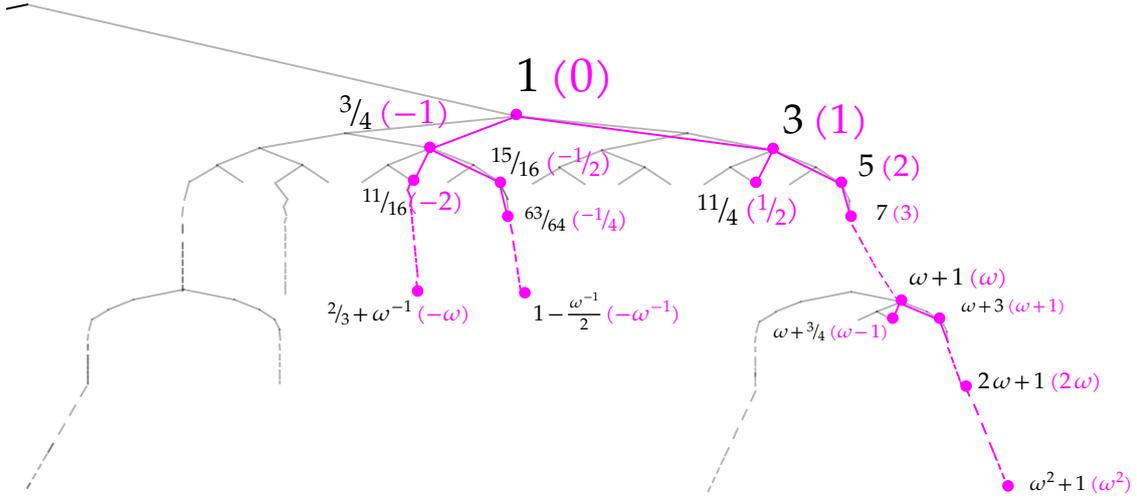

**Figure 4.2.** The (sub)tree representation of the surreal substructure **Inc** (purple) from Example 4.14 within **No** (grey). The labels have the form $\Xi_{\mathbf{Inc}}\, x\, (x)$. For instance $\Xi_{\mathbf{Inc}}(-2) = {}^{11}\!/_{16}$.

## 4.5 Convex subclasses

If $\mathbf{X} \subseteq \mathbf{Y}$ are subclasses of **No**, recall that **X** is *convex in* **Y** if

$$\forall x, z \in \mathbf{X}, \forall y \in \mathbf{Y}, (x \leqslant y \leqslant z \Longrightarrow y \in \mathbf{X}),$$

and **X** is $\sqsubseteq$-*convex in* **Y** if

$$\forall x, z \in \mathbf{X}, \forall y \in \mathbf{Y}, (x \sqsubseteq y \sqsubseteq z \Longrightarrow y \in \mathbf{X}).$$

We simply say that **X** is convex (resp. $\sqsubseteq$-convex) if it is convex (resp. $\sqsubseteq$-convex) in **No**. We let $\mathbf{Hull}_\mathbf{Y}(\mathbf{X})$ denote the convex hull of **X** in **Y**, that is, for every number $y$, we have $y \in \mathbf{Hull}_\mathbf{Y}(\mathbf{X})$ if and only if $y \in \mathbf{Y}$ and there are elements $x, z$ of **X** such that $x \leqslant y \leqslant z$. The convex hull of **X** in **Y** is the smallest convex subclass of **Y** containing **X**.

**Lemma 4.16.** *Assume that* **S** *is a surreal substructure. Then every non-empty convex subclass of* **S** *is rooted.*

**Proof.** In view of Propositions 4.6 and 4.7, it suffices to prove the lemma for $\mathbf{S} = \mathbf{No}$. Let **C** be a non-empty convex subclass of **No**. Assume for contradiction that $u, v \in \mathbf{C}$ are two simplest elements with $u < v$. Let $\alpha$ be the smallest ordinal such that $u[\alpha] < v[\alpha]$. Since $u \not\sqsubseteq v$ and $v \not\sqsubseteq u$, we must have $u[\alpha] = -1$ and $v[\alpha] = 1$. Now consider the number $w$ whose sign sequence is $u\!\upharpoonright\!\alpha = v\!\upharpoonright\!\alpha$. Then $u < w < v$, whence $w \in \mathbf{C}$, but also $w \sqsubseteq u$; a contradiction. $\square$

**Lemma 4.17.** *If* **C** *is a non-empty final segment of* **No**, *then* $\mathbf{C}^\bullet$ *is the smallest ordinal in* **C**.

**Proof.** Given $x \in \mathbf{C}$, we have $x \leqslant \ell(x) \in \mathbf{C}$, so **C** contains an ordinal. Let $\iota$ denote the smallest ordinal in **C**. Given another ordinal $\eta < \iota$, we have $\eta \notin \mathbf{C}$ by minimality of $\iota$. Since **C** is a final segment of **No**, it follows that $\eta < \mathbf{C}$. For any $x \in \mathbf{C}$, we deduce that $x$ lies in the cut $(\iota_L | \emptyset)$, whence $\iota = \{\iota_L | \emptyset\} \sqsubseteq x$. This shows that $\iota = \mathbf{C}^\bullet$. $\square$

**Proposition 4.18.** *Let **S** be a surreal substructure.*

a) *A convex subclass **C** of **S** is a surreal substructure if and only if it has no cofinal or coinitial subset.*

b) *For subsets $L < R$ of **S**, the cut $(L|R)_{\mathbf{S}}$ is a surreal substructure.*

c) *If $\mathbf{T} \subseteq \mathbf{S}$ is a surreal substructure, then $\mathbf{Hull}_{\mathbf{S}}(\mathbf{T})$ is a surreal substructure.*

d) *If **T** is a surreal substructure, $(L|R)_{\mathbf{S}}$ is a cut in **S** and $f: \mathbf{T} \longrightarrow \mathbf{S}$ is strictly monotonic and surjective, then $f^{-1}((L|R)_{\mathbf{S}})$ is a surreal substructure.*

e) *The intersection of any set-sized decreasing family of surreal substructures that are convex in **S** is a surreal substructure.*

**Proof.** *a*) Assume that **C** has no cofinal or coinitial subset and let $L < R$ be subsets of **C**.

- If both $L$ and $R$ are empty, then $L < c < R$ for any $c \in \mathbf{C}$. Notice that $\mathbf{C} \neq \emptyset$, since $\emptyset$ is not cofinal in **C**.
- If $L = \emptyset$ and $R \neq \emptyset$, then there exists an $x \in \mathbf{C}$ with $x < R$, since $R$ is not coinitial in **C**. Let $y = \{x | R\}_{\mathbf{S}}$ and $r \in R$. Then $x < y < r$, so $y \in \mathbf{C}$, and $y \in (L|R)_{\mathbf{C}}$.
- Similarly, if $L \neq \emptyset$ and $R = \emptyset$, then $\{L|y\}_{\mathbf{S}} \in (L|R)_{\mathbf{C}}$ for some $y > L$ in **C**.
- If $L \neq \emptyset$ and $R \neq \emptyset$, then $\{L|R\}_{\mathbf{S}} \in \mathbf{C}$, by convexity.

In each of the above cases, we have shown that $(L|R)_{\mathbf{C}}$ is a non-empty convex subclass of **S**. By Lemma 4.16, it is rooted. By Proposition 4.7, it follows that **C** is a surreal substructure. Conversely, if **C** is a surreal substructure, then given a subset $X$ of **C**, we have

$$\mathbf{C} \ni \{\emptyset | X\}_{\mathbf{C}} < X < \{X|\emptyset\}_{\mathbf{C}} \in \mathbf{C},$$

so $X$ is neither cofinal nor coinitial in **C**.

*b*) This is a direct consequence of the previous point: the cut $(L|R)_{\mathbf{S}}$ is by definition a convex subclass of **S**, and given a subset $X$ of $(L|R)_{\mathbf{S}}$ we have

$$(L|R)_{\mathbf{S}} \ni \{L|X\}_{\mathbf{S}} < X < \{X|R\}_{\mathbf{S}} \in (L|R)_{\mathbf{S}}.$$

By Proposition 4.7, it follows that $(L|R)_{\mathbf{S}}$ is a surreal substructure.

*c*) Since **T** is a surreal substructure, it has no cofinal or coinitial subset. It follows that the same holds for $\mathbf{Hull}_{\mathbf{S}}(\mathbf{T})$, which is thus a surreal substructure.

*d*) We have $f^{-1}((L|R)_{\mathbf{S}}) = (f^{-1}(L) | f^{-1}(R))_{\mathbf{T}}$ is $f$ is increasing and $f^{-1}((L|R)_{\mathbf{S}}) = (f^{-1}(R) | f^{-1}(L))_{\mathbf{T}}$ if $f$ is decreasing. In both cases, $f^{-1}((L|R)_{\mathbf{S}})$ is a cut in **T**, hence a surreal substructure by *c*).

*e*) Let $(I, <)$ be a linearly ordered set and let $(\mathbf{C}_i)_{i \in I}$ be decreasing for $\subseteq$. Its intersection $\mathbf{C} := \bigcap_{i \in I} \mathbf{C}_i$ is convex. Let $X$ be a subset of **C**. For $i \in I$, we have $X \subseteq \mathbf{C}_i$ whence $l_i < X < r_i$ where $l_i = (\emptyset | X)_{\mathbf{C}_i}^{\bullet}$ and $r_i = (X | \emptyset)_{\mathbf{C}_i}^{\bullet}$. Writing $l = \{l_i : i \in I | X\}_{\mathbf{S}}$ and $r = \{X | r_i : i \in I\}_{\mathbf{S}}$, we have $l < X < r$. Moreover, for $i \in I$, we have $l_i < l < r < r_i$ so $l, r \in \mathbf{C}_i$ by convexity. This proves that $l, r \in \mathbf{C}$ and consequently that $X$ is neither cofinal nor coinitial in **C**. Therefore **C** is a surreal substucture by *a*). □

**Example 4.19.** Cuts $(L|R)_{\mathbf{S}}$ where $L < R$ are subsets of **S** include $\sqsubseteq$-final substructures of **S** and non-empty open intervals of **S**, which are therefore convex surreal substructures. Note that non-empty convex classes of **No** which are open in the order topology may fail to be surreal substructures. One counterexample is the class $\mathbf{No}^{\preccurlyeq} := \mathbf{Hull}(\mathbb{Z})$ of *finite surreal numbers*, since it admits the cofinal subset $\mathbb{N}$.

**Example 4.20.** Here are some further examples and counterexamples of convex surreal substructures that we will consider later on.

- The class $\mathbf{No}^> := (\{0\}|\emptyset)$ of strictly positive surreal numbers is a convex surreal substructure, and it is in fact the $\sqsubseteq$-final substructure $\mathbf{No}^{\sqsupseteq 1}$ of $\mathbf{No}$.

- Likewise, the class $\mathbf{No}^{>,>} := (\mathbb{N}|\emptyset) = \mathbf{No}^{\sqsupseteq \omega}$ of positive *infinite surreal numbers* is a convex surreal substructure.

- The class $\mathbf{No}^{\prec} := (\mathbb{R}^{<0}|\mathbb{R}^{>0})$ of *infinitesimals* forms a surreal substructure which can be split as the union of $\{0\}$ and the two $\sqsubseteq$-final substructures $\mathbf{No}^{\sqsupseteq -\omega^{-1}}$, $\mathbf{No}^{\sqsupseteq \omega^{-1}}$.

- Although every interval $(-n-1, n+1)$ for $n \in \mathbb{N}$ is a convex surreal substructure, their increasing union $\mathbf{No}^{\lessapprox}$ is not a surreal substructure.

**Remark 4.21.** For subsets $L < R$ of $\mathbf{S}$, the cut $(L|R)_\mathbf{S}$ may fail to be a $\sqsubseteq$-final substructure of $\mathbf{S}$. In fact, by Proposition 4.11(c), it is a $\sqsubseteq$-final substructure of $\mathbf{S}$ if and only if the canonical representation of $\{L|R\}_\mathbf{S}$ in $\mathbf{S}$ is cofinal with respect to $(L,R)$, in which case we have $(L|R)_\mathbf{S} = \mathbf{S}^{\sqsupseteq\{L|R\}_\mathbf{S}}$.

Any convex subclass $\mathbf{C}$ of $\mathbf{S}$ is a generalized cut $\mathbf{C} = (\mathbf{L}|\mathbf{R})_\mathbf{S}$ in $\mathbf{S}$ where $\mathbf{L}$ is the class of strict lower bounds of $\mathbf{C}$ in $\mathbf{S}$ and $\mathbf{R}$ is the class of its strict upper bounds. However, those classes may not always be replaced by sets. In fact, the class $\mathbf{C}$ is a cut $\mathbf{C} = (L|R)_\mathbf{S}$ with subsets $L < R$ of $\mathbf{S}$ if and only if such sets can be found that are mutually cofinal with $(\mathbf{L}, \mathbf{R})$. The existence thus amounts to $\mathrm{cof}(\mathbf{L},<), \mathrm{cof}(\mathbf{B},>) \in \mathbf{On}$ since cofinality is invariant under mutual cofinality (see the end of Appendix B for notes about cofinal well-ordered subsets).

**Example 4.22.** Recall that $\omega/2 = \omega \dotplus (-\omega)$. Let $x_\alpha = \omega/2 \dottimes \alpha$ for each $\alpha \in \mathbf{On}$ and consider the class $\mathbf{C} = \{y \in \mathbf{No} : \forall \alpha \in \mathbf{On}, y > x_\alpha\}$. Then $\mathbf{C}$ is a convex surreal substructure of $\mathbf{No}$. Indeed, the sequence $(y_\alpha)_{\alpha \in \mathbf{On}}$ with $y_\alpha = \omega \dotplus (\omega/2 \dotttimes (-\alpha))$ is strictly decreasing and coinitial in $\mathbf{C}$. This shows that $\mathbf{C}$ does not admit a coinitial subset. As a non-empty final segment of $\mathbf{No}$, the class $\mathbf{C}$ also admits no cofinal subset. Proposition 4.18 thus implies that $\mathbf{C}$ is a surreal substructure. We have $\mathrm{cof}(\{x_\alpha : \alpha \in \mathbf{On}\}, <) = \mathbf{On}$, so $\mathbf{C}$ is not a cut in $\mathbf{No}$.

## 4.6 Cut equations

We already noted that the Conway bracket allows for elegant recursive definitions of functions on $\mathbf{No}$. Let us now study such definitions in more detail and examine how they generalize to arbitrary surreal substructures.

**Definition 4.23.** *Let $\mathbf{S}, \mathbf{T}$ be surreal substructures. Let $\lambda, \rho$ be functions defined for cut representations in $\mathbf{S}$ and such that $\lambda(L,R), \rho(L,R)$ are subsets of $\mathbf{T}$ whenever $(L,R)$ is a cut representation in $\mathbf{S}$. We say that a function $F: \mathbf{S} \longrightarrow \mathbf{T}$ has* **cut equation** $\{\lambda|\rho\}_\mathbf{T}$ *if for all $x \in \mathbf{S}$, we have*

$$\lambda(x_L^\mathbf{S}, x_R^\mathbf{S}) < \rho(x_L^\mathbf{S}, x_R^\mathbf{S}) \quad \text{and}$$
$$F(x) = \{\lambda(x_L^\mathbf{S}, x_R^\mathbf{S}) | \rho(x_L^\mathbf{S}, x_R^\mathbf{S})\}_\mathbf{T}.$$

*We say that the cut equation is **extensive** if it satisfies*

$$\forall x, y \in \mathbf{S}, (x \sqsubseteq y \Longrightarrow (\lambda(x_L^{\mathbf{S}}, x_R^{\mathbf{S}}) \subseteq \lambda(y_L^{\mathbf{S}}, y_R^{\mathbf{S}}) \wedge \rho(x_L^{\mathbf{S}}, x_R^{\mathbf{S}}) \subseteq \rho(y_L^{\mathbf{S}}, y_R^{\mathbf{S}}))).$$

**Note.** We will see in the proof of Proposition 4.27 below that extensive cut equations preserve simplicity.

**Example 4.24.** A simple example of a cut equation is (3.3): $\forall x \in \mathbf{No}, -x = \{-x_R \mid -x_L\}$. Here we have $\mathbf{S} = \mathbf{T} = \mathbf{No}$ and we can take $\lambda(x_L, x_R) = -x_R$ and $\rho(x_L, x_R) := -x_L$. Note that this cut equation is extensive.

Taking $\mathbf{S} = \mathbf{No}$ and $\mathbf{T} = \mathbf{No}^{>}$, $\lambda(x_L, x_R) = x_L \cap \mathbf{No}^{>}$ and $\rho(x_L, x_R) = x_R \cap \mathbf{No}^{>}$, we obtain the function $F$ with $F(x) = 0$ for all $x \leqslant 0$ and $F(x) = x$ for all $x > 0$.

See Example 4.32 below for more examples.

**Remark 4.25.** Our notion of cut equation is not restrictive on the function, since any function $F: \mathbf{S} \longrightarrow \mathbf{T}$ has cut equation $(\lambda, \rho)$ with $\lambda(L, R) := F(\{L \mid R\}_{\mathbf{S}})_L^{\mathbf{T}}$ and $\rho(L, R) := F(\{L \mid R\}_{\mathbf{S}})_R^{\mathbf{T}}$. Thus it should not be confused with the notions of *recursive definition* in [22] and *genetic definition* in [37].

**Example 4.26.** Given sets $\Lambda, P$ of functions $\mathbf{S} \longrightarrow \mathbf{T}$, cut equations of the form $(\lambda, \rho)$ with

$$\begin{aligned} \lambda(x_L^{\mathbf{S}}, x_R^{\mathbf{S}}) &= \{\xi(l) : \xi \in \Lambda, l \in x_L^{\mathbf{S}}\} \\ \rho(x_L^{\mathbf{S}}, x_R^{\mathbf{S}}) &= \{\psi(r) : \psi \in P, r \in x_R^{\mathbf{S}}\} \end{aligned}$$

are extensive. We will write $\{\lambda(x_L^{\mathbf{S}}, x_R^{\mathbf{S}}) \mid \rho(x_L^{\mathbf{S}}, x_R^{\mathbf{S}})\}_{\mathbf{T}} = \{\Lambda(x_L^{\mathbf{S}}) \mid P(x_R^{\mathbf{S}})\}_{\mathbf{T}}$ in this case. Note that it is common to consider well-defined cut equations of the form

$$F(x) = \{\Lambda(x_L^{\mathbf{S}}) \mid P(x_R^{\mathbf{S}})\}_{\mathbf{T}},$$

where $F$ itself belongs to $\Lambda$ and $P$.

**Proposition 4.27.** *Let $\mathbf{S}, \mathbf{T}$ be surreal substructures. Let $F: \mathbf{S} \longrightarrow \mathbf{T}$ be strictly $\leqslant$-increasing with extensive cut equation $\{\lambda \mid \rho\}_{\mathbf{T}}$. Then $F(\mathbf{S})$ is a surreal substructure, and we have $F = \Xi_{F(\mathbf{S})}^{\mathbf{S}}$.*

**Proof.** We claim that $F$ is $\sqsubseteq$-increasing. Indeed, let $x, y \in \mathbf{S}$ with $x \sqsubseteq y$. We have $x_L^{\mathbf{S}} < y < x_R^{\mathbf{S}}$, so $x_L^{\mathbf{S}} \subseteq y_L^{\mathbf{S}}$ and $x_R^{\mathbf{S}} \subseteq y_R^{\mathbf{S}}$. We deduce by extensivity of $(\lambda, \rho)$ that $\lambda(x_L^{\mathbf{S}}, x_R^{\mathbf{S}}) \subseteq \lambda(y_L^{\mathbf{S}}, y_R^{\mathbf{S}})$ and $\rho(x_L^{\mathbf{S}}, x_R^{\mathbf{S}}) \subseteq \rho(y_L^{\mathbf{S}}, y_R^{\mathbf{S}})$, and thus $\lambda(x_L^{\mathbf{S}}, x_R^{\mathbf{S}}) < F(y) < \rho(x_L^{\mathbf{S}}, x_R^{\mathbf{S}})$. This implies that $F(x) \sqsubseteq F(y)$. Thus $F$ is strictly $(\leqslant, \sqsubseteq)$-increasing. So the composition $F \circ \Xi_{\mathbf{S}}: \mathbf{No} \longrightarrow F(\mathbf{S})$ is strictly $(\leqslant, \sqsubseteq)$-increasing. The function $\Xi_{\mathbf{S}}: (\mathbf{No}, \leqslant, \sqsubseteq) \longrightarrow (\mathbf{S}, \leqslant, \sqsubseteq)$ is an embedding by Proposition 4.6, so $F$ embeds $\mathbf{S}$ into $\mathbf{T}$. In particular, $F(\mathbf{S})$ is a surreal substructure. By Proposition 4.9, we conclude that $F = \Xi_{F(\mathbf{S})}^{\mathbf{S}}$. □

As an application, we get the following well-known result (see [11, Proposition 4.22]).

**Proposition 4.28.** *Let $\varphi$ be a number, and let $\mathbf{No}^{<\operatorname{supp}\varphi}$ denote the class of numbers $x$ with $x \prec \operatorname{supp}\varphi$. Then $\mathbf{No}^{<\operatorname{supp}\varphi}$ and $\varphi + \mathbf{No}^{<\operatorname{supp}\varphi}$ are surreal substructures with*

$$\forall x \in \mathbf{No}, \Xi_{\varphi + \mathbf{No}^{<\operatorname{supp}\varphi}} x = \varphi + \Xi_{\mathbf{No}^{<\operatorname{supp}\varphi}} x.$$

**Proof.** We have $\mathbf{No}^{<\mathrm{supp}\,\varphi} = (-\mathbb{R}^{>}\mathrm{supp}\,\varphi \,|\, \mathbb{R}^{>}\mathrm{supp}\,\varphi)$. By Proposition 4.18(b), this is a surreal substructure. Recall that for $x \in \mathbf{No}$, we have $\varphi + x = \{\varphi_L + x, \varphi + x_L \,|\, \varphi + x_R, \varphi_R + x\}$. If $x \in \mathbf{No}^{<\mathrm{supp}\,\varphi}$, then we have $\varphi_L + x < \varphi + \mathbf{No}^{<\mathrm{supp}\,\varphi} < \varphi_R + x$ so we may write

$$\varphi + x = \{\varphi + x_L \,|\, \varphi + x_R\}_{\varphi + \mathbf{No}^{<\mathrm{supp}\,\varphi}}$$
$$= \{\varphi + x_L^{\mathbf{No}^{<\mathrm{supp}\,\varphi}} \,|\, \varphi + x_R^{\mathbf{No}^{<\mathrm{supp}\,\varphi}}\}_{\varphi + \mathbf{No}^{<\mathrm{supp}\,\varphi}}.$$

Seen as a cut equation in $x$, this is an extensive cut equation, so by Proposition 4.27, we see that $\varphi + \mathbf{No}^{<\mathrm{supp}\,\varphi}$ is a surreal substructure and that $x \mapsto \varphi + x$ realizes the isomorphism $\mathbf{No}^{<\mathrm{supp}\,\varphi} \longrightarrow \varphi + \mathbf{No}^{<\mathrm{supp}\,\varphi}$. □

**Definition 4.29.** *Let $F$ be a function $\mathbf{S} \longrightarrow \mathbf{T}$ with cut equation $(\lambda, \rho)$. We say that $(\lambda, \rho)$ is **uniform** at $x \in \mathbf{S}$ if we have*

$$\lambda(L, R) < \rho(L, R) \quad \text{and}$$
$$F(x) = \{\lambda(L, R) \,|\, \rho(L, R)\}$$

*whenever $(L, R)$ is a cut representation of $x$ in $\mathbf{S}$. We say that $(\lambda, \rho)$ is **uniform** if it is uniform at every $x \in \mathbf{S}$.*

**Example 4.30.** Let $a \in \mathbf{No}$. The following cut equation for the function $y \mapsto a \dotplus y : \mathbf{No} \longrightarrow 1 \dotplus \mathbf{No}$ obtained from (3.8)

$$\forall x \in \mathbf{No}, a \dotplus y = \{a_L, a \dotplus y_L \,|\, a \dotplus y_R, a_R\},$$

is uniform. On the contrary, the following cut equation for $x \mapsto x \dotplus 1$ is not uniform:

$$\forall x \in \mathbf{No}, x \dotplus 1 = \{x, x_L \,|\, x_R\}.$$

Indeed, we have $0 = \{\emptyset \,|\, 1\}$ and $0 \dotplus 1 = 1$, but $\{0, \emptyset \,|\, 1\} = \{0 \,|\, 1\} = {}^1\!/_2$.

**Example 4.31.** Let $b \in \mathbf{No}^{>}$. By (3.7), the function $y \mapsto b \dot\times y : \mathbf{No} \longrightarrow b \dot\times \mathbf{No}$ has the following cut equation

$$\forall y \in \mathbf{No}, b \dot\times y = \{b \dot\times y_L \dotplus b_L, b \dot\times y_R \dotplus (-b_R) \,|\, b \dot\times y_L \dotplus b_R, b \dot\times y_R \dotplus (-b_L)\},$$

which is uniform. On the contrary, the cut equation for $x \mapsto x \dot\times {}^1\!/_2$ is not uniform:

$$\forall x \in \mathbf{No}, x \dot\times {}^1\!/_2 = \{x_L, x \dotplus (-x_R) \,|\, x_R, x \dotplus (-x_L)\}.$$

Indeed, if we were to apply this cut equation to the cut presentation $(\{{}^1\!/_2\}, \emptyset)$ of $1$, then we would have ${}^1\!/_2$ as a left option and $1 \dotplus (-{}^1\!/_2) \leqslant {}^1\!/_2$ as a right option, which cannot be.

**Example 4.32.** Most common definitions of unary functions $\mathbf{No} \longrightarrow \mathbf{No}$ have known simple cut equations, and many of them are uniform, in particular throughout the work of H. Gonshor in [24]. For instance, the classical cut equations (3.3) and (3.6) for the functions $x \mapsto -x$ and $x \mapsto \exp x$ are uniform, so for $x \in \mathbf{No}$ and for any cut representation $(L_x, R_x)$ of $x$ in $\mathbf{No}$, we have

$$-x = \{-R_x \,|\, -L_x\}, \text{ and}$$
$$\exp x = \left\{ 0,\, [x-l]_{\mathbb{N}} \exp l,\, [x-r]_{2\mathbb{N}+1} \exp r \,\bigg|\, \frac{\exp r}{[x-r]_{2\mathbb{N}+1}},\, \frac{\exp l}{[l-x]_{\mathbb{N}}} \right\} \qquad (l \in L_x, r \in R_x).$$

**Example 4.33.** We will also need an extension of the notion of uniform cut equation to functions $f: \mathbf{No} \times \mathbf{No} \longrightarrow \mathbf{No}$. Specifically, by [24, Theorem 3.2], the classical cut equation (3.4) for the sum of two numbers $x, y$ is uniform in the sense that, given cut representations $(L_x, R_x)$ and $(L_y, R_y)$ of $x, y$ in $\mathbf{No}$, we have

$$x + y = \{L_x + y, x + L_y | x + R_y, R_x + y\}. \tag{4.1}$$

Similarly for the multiplication, we have

$$x + y = \{x'y + xy' - x'y', x''y + xy'' - x''y'' | x'y + xy'' - x'y'', x''y + xy' - x''y'\},$$

where $x'$, $x''$, $y'$ and $y''$ range in $L_x$, $R_x$, $L_y$ and $R_y$ respectively.

Uniform cut equations have the interesting property that they can be composed.

**Lemma 4.34.** *Let* $\mathbf{S}_0, \mathbf{S}_1, \mathbf{S}_2$ *be surreal substructures. Let* $F_1: \mathbf{S}_0 \longrightarrow \mathbf{S}_1$ *and* $F_2: \mathbf{S}_1 \longrightarrow \mathbf{S}_2$ *be functions with uniform cut equations*

$$F_1 \equiv \{\lambda_1 | \rho_1\}_{\mathbf{S}_1}$$
$$F_2 \equiv \{\lambda_2 | \rho_2\}_{\mathbf{S}_2}.$$

*Then* $F_2 \circ F_1$ *has the uniform cut equation* $(\lambda_{12}, \rho_{12})$ *where for every cut representation* $(L, R)$ *in* $\mathbf{S}_0$, *we have* $\lambda_{12}(L, R) = \lambda_2(\lambda_1(L, R), \rho_1(L, R))$ *and* $\rho_{12}(L, R) = \rho_2(\lambda_1(L, R), \rho_1(L, R))$.

**Proof.** Let $x \in \mathbf{S}_0$, let $(L, R)$ be a cut representation of $x$ in $\mathbf{S}_0$. By uniformity of the cut equation of $F_1$ at $x$, we have

$$F_1(x) = \{\lambda_1(L, R) | \rho_1(L, R)\}_{\mathbf{S}_1}.$$

By uniformity of the cut equation of $F_2$ at $F_1(x)$, we have

$$F_2(F_1(x)) = \{\lambda_2(\lambda_1(L, R), \rho_1(L, R)) | \rho_2(\lambda_1(L, R), \rho_1(L, R))\},$$

whence the result. □

Recall that a class $\mathbf{X} \subseteq \mathbf{No}$ is *cofinal* (resp. *coinitial*) *with respect to* a class $\mathbf{Y} \subseteq \mathbf{No}$ if every element of $\mathbf{Y}$ has an upper bound (resp. lower bound) in $\mathbf{X}$. If $\mathbf{X} \subseteq \mathbf{Y}$, then we simply say that $\mathbf{X}$ is cofinal (resp. coinitial) *in* $\mathbf{Y}$.

**Lemma 4.35.** *When* $\mathbf{S}, \mathbf{T}$ *are surreal substructures, the cut equation* $\Xi_{\mathbf{T}}^{\mathbf{S}} x \equiv \{\Xi_{\mathbf{T}}^{\mathbf{S}} x_L^{\mathbf{S}} | \Xi_{\mathbf{T}}^{\mathbf{S}} x_R^{\mathbf{S}}\}_{\mathbf{T}}$ *is uniform and extensive.*

**Proof.** Let us first prove uniformity in the case when $\mathbf{S} = \mathbf{No}$. Let $L < R$ be sets of surreal numbers and let $x = \{L|R\}$. Since $\Xi_{\mathbf{T}}$ is strictly increasing and ranges in $\mathbf{T}$, the number $y = \{\Xi_{\mathbf{T}} L | \Xi_{\mathbf{T}} R\}_{\mathbf{T}}$ is well defined and $\Xi_{\mathbf{T}} L < \Xi_{\mathbf{T}} x < \Xi_{\mathbf{T}} R$, which yields $y \sqsubseteq \Xi_{\mathbf{T}} x$. Moreover, the set $L$ is cofinal in $x_L$ whereas $R$ is coinitial in $x_R$, so $\Xi_{\mathbf{T}} x_L < y < \Xi_{\mathbf{T}} x_R$. Hence $\Xi_{\mathbf{T}} x \sqsubseteq y$ and $\Xi_{\mathbf{T}} x = y$, which shows that the cut equation $\Xi_{\mathbf{T}} x \equiv \{\Xi_{\mathbf{T}} x_L | \Xi_{\mathbf{T}} x_R\}_{\mathbf{T}}$ is uniform.

Now consider the general case and let $\Xi_{\mathbf{S}} A = L < R = \Xi_{\mathbf{S}} B$ be subsets of $\mathbf{S}$. Setting $z := \{A|B\}$ and $x := \{L|R\}_{\mathbf{S}}$, we have $x = \Xi_{\mathbf{S}} z$ by uniformity of the cut equation for $\Xi_{\mathbf{S}}$. Furthermore,

$$\{\Xi_{\mathbf{T}}^{\mathbf{S}} L | \Xi_{\mathbf{T}}^{\mathbf{S}} R\}_{\mathbf{T}} = \{\Xi_{\mathbf{T}} A | \Xi_{\mathbf{T}} B\}_{\mathbf{T}}$$
$$= \Xi_{\mathbf{T}} z,$$

by uniformity of the cut equation for $\Xi_{\mathbf{T}}$. Hence $\{\Xi_{\mathbf{T}}^{\mathbf{S}} L \,|\, \Xi_{\mathbf{T}}^{\mathbf{S}} R\}_{\mathbf{T}} = \Xi_{\mathbf{T}} \Xi_{\mathbf{S}}^{-1} x = \Xi_{\mathbf{T}}^{\mathbf{S}} z$, which proves that $\Xi_{\mathbf{T}}^{\mathbf{S}} \equiv \{\Xi_{\mathbf{T}}^{\mathbf{S}} L \,|\, \Xi_{\mathbf{T}}^{\mathbf{S}} R\}_{\mathbf{T}}$ is uniform. This cut equation has the form $\Xi_{\mathbf{T}}^{\mathbf{S}} z = \{\Lambda(z_L^{\mathbf{S}}) \,|\, \mathrm{P}(z_R^{\mathbf{S}})\}_{\mathbf{T}}$ where $\Lambda = \mathrm{P} = \{\Xi_{\mathbf{T}}^{\mathbf{S}}\}$ are sets of functions, so it is extensive. □

The above proposition shows that surreal isomorphisms satisfy natural extensive cut equations. Inversily, Proposition 4.27 shows that extensive cut equations give rise to surreal isomorphisms. As an application, if we admit that the operation

$$\forall x \in \mathbf{No},\ \dot{\omega}^x := \{0, \mathbb{N}\,\dot{\omega}^{x_L} \,|\, 2^{-\mathbb{N}}\dot{\omega}^{x_R}\}$$

is well defined, then we see that it defines a surreal isomorphism. This is the parametrization of the class $\mathbf{Mo}$ of *monomials*, that is, Conway's $\omega$-map. This cut equation is also uniform (see [24, corollary of Theorem 5.2]), and we can for instance compute, for every number $x$, the number

$$\begin{aligned}
\dot{\omega}^{\dot{\omega}^x} &= \dot{\omega}^{\{0, \mathbb{N}\dot{\omega}^{x_L} | 2^{-\mathbb{N}}\dot{\omega}^{x_R}\}} \\
&= \{0, \mathbb{N}\,\dot{\omega}^0, \mathbb{N}\,\omega^{\mathbb{N}\dot{\omega}^{x_L}} \,|\, 2^{-\mathbb{N}}\dot{\omega}^{2^{-\mathbb{N}}\dot{\omega}^{x_R}}\} \\
&= \{\mathbb{N}, \omega^{\mathbb{N}\dot{\omega}^{x_L}} \,|\, \dot{\omega}^{2^{-\mathbb{N}}\dot{\omega}^{x_R}}\}.
\end{aligned}$$

Whenever they exist, this shows the usefulness of extensive cut equations. Unfortunately, many common surreal functions such as the exponential do not admit extensive cut equations. The next proposition describes a more general type of cut equation that is sometimes useful.

**Proposition 4.36.** *Let $\mathbf{S}, \mathbf{T}$ be surreal substructures. Let $\Lambda$ be a function from $\mathbf{S}$ to the class of subsets of $\mathbf{T}$ such that for $x, y \in \mathbf{S}$ with $x < y$, the set $\Lambda(y)$ is cofinal with respect to $\Lambda(x)$. For $x \in \mathbf{S}$, let $\Lambda[x]$ denote the class of elements $u$ of $\mathbf{S}$ such that $\Lambda(x)$ and $\Lambda(u)$ are mutually cofinal. Let $\{\lambda \,|\, \rho\}_{\mathbf{T}}$ be an extensive cut equation on $\mathbf{S}$. Let $F: \mathbf{S} \longrightarrow \mathbf{T}$ be strictly increasing with cut equation*

$$\forall x \in \mathbf{S}, F(x) = \{\Lambda(x), \lambda(x_L^{\mathbf{S}}, x_R^{\mathbf{S}}) \,|\, \rho(x_L^{\mathbf{S}}, x_R^{\mathbf{S}})\}_{\mathbf{T}}$$

*Then $F$ induces an embedding $(\Lambda[x], \leqslant, \sqsubseteq) \longrightarrow (\mathbf{T}, \leqslant, \sqsubseteq)$ for each element $x$ of $\mathbf{S}$.*

**Proof.** Let $x \in \mathbf{S}$. If $u, w \in \Lambda[x]$ and $v \in \mathbf{S}$ satisfies $u \leqslant v \leqslant w$, then $\Lambda(v)$ is cofinal with respect to $\Lambda(u)$ and hence to $\Lambda(x)$, and $\Lambda(x)$ is cofinal with respect to $\Lambda(w)$ and hence to $\Lambda(v)$, so $v \in \Lambda[x]$. Therefore $\Lambda[x]$ is a non-empty convex subclass of $\mathbf{S}$. Note that for $u \in \Lambda[x]$, we have

$$F(u) = \{\Lambda(x), \lambda(u_L^{\mathbf{S}}, u_R^{\mathbf{S}}) \,|\, \rho(u_L^{\mathbf{S}}, u_R^{\mathbf{S}})\}_{\mathbf{T}}.$$

For numbers $u, v$ lying in $\Lambda[x]$ with $u \sqsubseteq v$, we have

$$\Lambda(x) \cup \lambda(u_L^{\mathbf{S}}, u_R^{\mathbf{S}}) \subseteq \Lambda(x) \cup \lambda(v_L^{\mathbf{S}}, v_R^{\mathbf{S}}) < F(v) < \rho(v_L^{\mathbf{S}}, v_R^{\mathbf{S}}) \supseteq \rho(u_L^{\mathbf{S}}, u_R^{\mathbf{S}}),$$

which implies that $F(u) \sqsubseteq F(v)$. Since $\Lambda[x]$ is a non-empty convex subclass of $\mathbf{S}$ and $\Xi_{\mathbf{S}}: \mathbf{No} \longrightarrow \mathbf{S}$ is increasing and bijective, the class $\mathbf{C} := \Xi_{\mathbf{S}}^{-1}(\Lambda[x])$ is a non-empty convex subclass of $\mathbf{No}$ on which $F \circ \Xi_{\mathbf{S}}$ is strictly $(\leqslant, \sqsubseteq)$-increasing. By Lemma 4.5, the function $F \circ \Xi_{\mathbf{S}}$ induces an embedding $(\mathbf{C}, \leqslant, \sqsubseteq) \longrightarrow (\mathbf{T}, \leqslant, \sqsubseteq)$ and thus $F$ induces an embedding $(\Lambda[x], \leqslant, \sqsubseteq) \longrightarrow (\mathbf{T}, \leqslant, \sqsubseteq)$. □

**Example 4.37.** A typical example is the following cut equation of [11, Theorem 3.8(1)] for the exponential function on the class $\mathbf{Mo}^{\succ} := \{\mathfrak{m} \in \mathbf{Mo} : \mathbb{R} \prec \mathfrak{m}\}$ of infinite monomials:

$$\forall \mathfrak{m} \in \mathbf{Mo}, \exp \mathfrak{m} = \{\mathfrak{m}^{\mathbb{N}}, (\exp \mathfrak{m}_L^{\mathbf{Mo}})^{\mathbb{N}} \mid (\exp \mathfrak{m}_R^{\mathbf{Mo}})^{\mathbb{N}}\}.$$

Here we have $\Lambda(\mathfrak{m}) = \mathfrak{m}^{\mathbb{N}}$ and $\Lambda[\mathfrak{m}] = \{\mathfrak{n} \in \mathbf{Mo}^{\succ} : \exists p, q \in \mathbb{N}, \mathfrak{m}^{1/p} \prec \mathfrak{n} \prec \mathfrak{m}^p\}$.

# 5 Fixed points

After introducing the $\omega$-map as a way to parameterize the class $\mathbf{Mo}$ of monomials, Conway remarks that for any ordinal $\alpha$, the number $\dot{\omega}^\alpha$ coincides with Cantor's $\alpha$-th ordinal power of $\omega$. He then goes on with the definition of *generalized $\varepsilon$-numbers* as surreal numbers $a$ such that $\dot{\omega}^a = a$. It turns out that the class of generalized $\varepsilon$-numbers can be parameterized as well and actually forms a surreal substructure: see Conway's informal discussion [14, p 34–35] and Gonshor's formal proof [24, Theorem 9.1 and Corollary 9.2]. Gonshor gives further conditions for the class of fixed points of a surreal function to be a surreal substructure [24, Theorem 9.4].

In this section, we consider the more general problem of deciding, given a surreal substructure $\mathbf{S}$, whether $\Xi_{\mathbf{S}}$ admits fixed points, and possibly a whole surreal substructure of fixed points. A related fixed point theorem was obtained by Lurie [34, Theorem 8.2] in a somewhat different context.

## 5.1 Fixed points and iterations of defining isomorphisms

For operators $\Omega: \mathbf{X} \longrightarrow \mathbf{Y}$ where $\mathbf{Y} \subseteq \mathbf{X}$ are subclasses of $\mathbf{No}$ and $n \in \mathbb{N}$, it will be convenient to write $\Omega^n$ for the $n$-fold composition of $\Omega$ with itself. In particular, $\Omega^0 = \mathrm{id}_{\mathbf{X}}$.

**Definition 5.1.** *Let $\mathbf{S}$ be a surreal substructure. We say that a number $x$ is $\mathbf{S}$-fixed if $\Xi_{\mathbf{S}} x = x$. We let $\mathbf{Fix}_{\mathbf{S}}$ denote the class of $\mathbf{S}$-fixed numbers. Notice that $\mathbf{Fix}_{\mathbf{S}}$ is a subclass of $\mathbf{S}$.*

If $\mathbf{U}, \mathbf{V}, \mathbf{W}$ are surreal substructures with $\mathbf{U} = \mathbf{V} \prec \mathbf{W}$, then for every number $x$, we have $\Xi_{\mathbf{U}} x \geqslant \Xi_{\mathbf{V}} x$ if and only if $\Xi_{\mathbf{W}} x \geqslant x$, and $\Xi_{\mathbf{U}} x \sqsupseteq \Xi_{\mathbf{V}} x$ if and only if $\Xi_{\mathbf{W}} x \sqsupseteq x$. In particular, the parametrizations of $\mathbf{U}$ and $\mathbf{V}$ coincide exactly on $\mathbf{Fix}_{\mathbf{W}}$.

**Proposition 5.2.** *If $\mathbf{S}$ is a surreal substructure, then $\mathbf{Fix}_{\mathbf{S}} = \bigcap_{n \in \mathbb{N}} \Xi_{\mathbf{S}}^n \mathbf{No}$.*

**Proof.** Let $\mathbf{S}^{\prec \omega} = \bigcap_{n \in \mathbb{N}} \Xi_{\mathbf{S}}^n \mathbf{No}$. For $n \in \mathbb{N}$, we have $\mathbf{Fix}_{\mathbf{S}} = \Xi_{\mathbf{S}}^n \mathbf{Fix}_{\mathbf{S}} \subseteq \Xi_{\mathbf{S}}^n \mathbf{No}$, so $\mathbf{Fix}_{\mathbf{S}} \subseteq \mathbf{S}^{\prec \omega}$. Assume for contradiction that $\mathbf{Fix}_{\mathbf{S}}$ is a proper subclass of $\mathbf{S}^{\prec \omega}$, and consider $x \in \mathbf{S}^{\prec \omega} \setminus \mathbf{Fix}_{\mathbf{S}}$ with minimal length. For $n \in \mathbb{N}^{>}$, let $x_n \in \mathbf{No}$ with $x = \Xi_{\mathbf{S}}^n x_n$. For all $n \in \mathbb{N}^{>}$, we have $x_n \in \mathbf{S}^{\prec \omega}$, so by our minimality assumption and Lemma 4.10, we have $\forall n \in \mathbb{N}, \ell(x) = \ell(x_n)$.

Recall that $x$ is not $\mathbf{S}$-fixed, so $x_0 \neq x_1$. By symmetry, we may assume without loss of generality that $x_0 < x_1$, which implies that $x_n < x_{n+1}$ for all $n \in \mathbb{N}$. For $n \in \mathbb{N}$, let $u_n$ be the $\sqsubseteq$-maximal element of $\mathbf{S}$ with $u_n \sqsubseteq x_n, x_{n+1}$. This element is well-defined since $\mathbf{S}$ is a surreal substructure and $x_n, x_{n+1} \in \mathbf{S}$. The number $\Xi_{\mathbf{S}}^{-1} u_n$ is $\sqsubseteq$-maximal in $\mathbf{No}$ with $\Xi_{\mathbf{S}}^{-1} u_n \sqsubseteq x_{n+1}, x_{n+2}$, whence $u_{n+1} \sqsubseteq \Xi_{\mathbf{S}}^{-1} u_n$, so $\Xi_{\mathbf{S}} u_{n+1} \sqsubseteq u_n$.

Since $\ell(x_n) = \ell(x_{n+1})$ and $x_{n+1} \neq x_n$, we have $x_n \not\sqsubseteq x_{n+1}$ and $x_{n+1} \not\sqsubseteq x_n$. We deduce that $u_n \sqsubset x_n, x_{n+1}$ and that $x_n < u_n < x_{n+1}$. In particular, we have $x_{n+1} = \Xi_{\mathbf{S}}^{-1} x_n < \Xi_{\mathbf{S}}^{-1} u_n$ so $u_n < \Xi_{\mathbf{S}}^{-1} u_n$, so $u_n$ is not $\mathbf{S}$-fixed, and we have $\ell(u_n) < \ell(x_n) = \ell(x)$.

Since $\Xi_\mathbf{S} u_{n+1} \sqsubseteq u_n$ for each $n \in \mathbb{N}$, Lemma 4.10 implies $\ell(u_0) \geqslant \ell(u_1) \geqslant \cdots$. The latter decreasing sequence of ordinals is necessarily stationary; let $n_0 \in \mathbb{N}$ be such that $\ell(u_n) = \ell(u_{n_0})$ for all $n \geqslant n_0$. By Lemma 4.10, it follows that $\Xi_\mathbf{S} u_{n+1} = u_n$ for all $n \geqslant n_0$, whence $u_{n_0} \in \mathbf{S}^{\prec \omega} \setminus \mathbf{Fix}_\mathbf{S}$. But $\ell(u_{n_0}) < \ell(x)$, which contradicts the minimality of $\ell(x)$. This absurdity completes our proof. □

**Example 5.3.** Here are some examples of structures of fixed points where $\bigcap_{n \in \mathbb{N}} \Xi_\mathbf{S}^n \mathbf{No}$ is a surreal substructure:

- If $\mathbf{S}$ is the $\sqsubseteq$-final substructure $a \dotplus \mathbf{No} = \mathbf{No}^{\sqsupseteq a}$, then for any surreal number $x$, the sign sequence of $\Xi_\mathbf{S} x = a \dotplus x$ is obtained through concatenation of the sign sequences of $a$ and $x$. Thus $\mathbf{S}$-fixed numbers are numbers whose sign sequences start with $\omega$ copies of the sign sequence of $a$, that is $\mathbf{Fix}_{\mathbf{No}^{\sqsupseteq a}} = \mathbf{No}^{\sqsupseteq a \dot\times \omega}$.

- Consider $\mathbf{S} = a \dot\times \mathbf{No}$ where $a$ is a strictly positive number. Let $a_0 = 1$ and $a_n = a \dot\times a_n$ for $n \in \mathbb{N}$. We claim that $\mathbf{Fix}_\mathbf{S} = a_\omega \dot\times \mathbf{No}$ where $a_\omega = \sup_\sqsubseteq \{a_n : n \in \mathbb{N}\}$.

  Indeed, since $1 \sqsubseteq a_1$ and $a \dot\times \cdot$ is a surreal isomorphism, we have $a_n \sqsubseteq a_{n+1}$ for every $n \in \mathbb{N}$, so $a_\omega$ is well defined. We have $a \dot\times a_\omega = \sup_\sqsubseteq a \dot\times a_\mathbb{N} = \sup_\sqsubseteq a_{\mathbb{N}+1} = a_\omega$. For every number $x = a_\omega \dot\times x'$ where $x' \in \mathbf{No}$, we have $a \dot\times x = (a \dot\times a_\omega) \dot\times x' = a_\omega \dot\times x' = x$, so $a_\omega \dot\times \mathbf{No} \subseteq \mathbf{Fix}_\mathbf{S}$. Conversely if $x \in \mathbf{Fix}_{a \dot\times \mathbf{No}}$, then $a \dot\times x \sqsubseteq x$ so $\ell(a) \dot\times \ell(x) \leqslant \ell(x)$, so $\ell(x)$ is equal to $\ell(a)^\omega \dot\times \alpha$ for some ordinal $\alpha$. For $n \in \mathbb{N}$, $\ell(a_n) = \ell(a)^n$, so $\ell(a_\omega) = \sup_{n \in \mathbb{N}} \ell(a_n) = \ell(a)^\omega$. Let $x'$ denote the number of length $\alpha$ defined at the level of sign sequences by

  $$\forall \beta < \alpha, x'[\beta] = x[\ell(a)^\omega \dot\times \beta].$$

We claim that $x = a_\omega \dot\times x'$. Indeed, for $\beta < \alpha$ and $\gamma < \ell(a)^\omega$, there is $n \in \mathbb{N}$ such that $\gamma < \ell(a)^n$, and we have

$$\begin{aligned}(a_\omega \dot\times x')[\ell(a)^\omega \dot\times \beta \dotplus \gamma] &= x'[\beta] a_\omega[\gamma] \\ &= x[\ell(a)^\omega \dot\times \beta] a_n[\gamma] \\ &= (a_n \dot\times x)[\ell(a)^n \dot\times (\ell(a)^\omega \dot\times \beta) \dotplus \gamma] \\ &= (a_n \dot\times x)[(\ell(a)^\omega \dot\times \beta) \dotplus \gamma] \\ &= x[\ell(a)^\omega \dot\times \beta \dotplus \gamma].\end{aligned}$$

  Thus $a_\omega \dot\times x' = x$, so $\mathbf{Fix}_\mathbf{S} = a_\omega \dot\times \mathbf{No}$.

  We let $\mathbf{No}_>$ denote the surreal substructure $\mathbf{Fix}_{2 \dot\times \mathbf{No}} = \omega \dot\times \mathbf{No}$ which is the class of surreal numbers, whose sign sequence contains no consecutive distinct signs. Elements in $f \in \mathbf{No}_>$ are called *purely infinite numbers*, since their supports $\operatorname{supp} f$ as series $f = \sum_{\mathfrak{m} \in \mathbf{Mo}} f_\mathfrak{m} \mathfrak{m}$ contains only infinitely large monomials: see Proposition 7.4 below.

- As mentioned at the beginning of this section, if $\mathbf{S} = \mathbf{Mo}$ is the class of monomials, then $\Xi_\mathbf{S}$ is the $\omega$-map $x \mapsto \dot\omega^x$, and its fixed points are called generalized $\varepsilon$-numbers. For $x \in \mathbf{No}$, the number $\Xi_{\mathbf{Fix}_{\mathbf{Mo}}} x$ is usually denoted $\varepsilon_x$, and the $\varepsilon$-map $x \mapsto \varepsilon_x$ extends the parametrization of $\varepsilon$-numbers in $\mathbf{On}$. We refer to [24, Chapter 9] for a detailed study.

- If $\mathbf{S} = 1 \dotplus \mathbf{Mo}^\prec$ (where $\mathbf{Mo}^\prec = \mathbf{Mo} \cap \mathbf{No}^\prec$), then for $x \in \mathbf{No}$, we have

$$\Xi_\mathbf{S} x = 1 + \dot\omega^{(-1) \dotplus x}.$$

Consider the function $\Phi: x \longmapsto 1 + \dot{\omega}^{x-3/2}: \mathbf{No} \longrightarrow \mathbf{No}$. For all $y \in \mathbf{No}^{\prec}$ and $r \in \mathbb{R}$, we have $r + y = r \dotplus y$ by [24, Theorem 5.12]. Recall that $-1/2 = (-1) \dotplus 1$. Thus for $x \in 1 + \mathbf{No}^{\prec}$, we have

$$x - 3/2 = (x-1) - 1/2 = ((-1) \dotplus 1) + (x-1) = (-1) \dotplus (1 \dotplus (x-1)) = (-1) \dotplus x.$$

So $\Xi_\mathbf{S}$ and $\Phi$ coincide on $1 + \mathbf{No}^{\prec}$. Since $\mathbf{S}$ and the class of fixed points of $\Phi$ are contained in $1 + \mathbf{No}^{\prec}$, we deduce that $\mathbf{Fix_S}$ is the class of fixed points of $\Phi$.

Now, informally speaking, we would like to consider the expression

$$1 \dotplus \dot{\omega}^{-1/2 \dotplus \dot{\omega}^{-1/2 + \dot{\omega}^{\cdot^{\cdot^{\cdot}}}}}$$

as a notation for "the" fixed point of $\Phi$. However, this expression is inherently ambiguous, since $\mathbf{Fix_S}$ actually contains many elements. The map $\Xi_{\mathbf{Fix_S}}$ can be regarded as a notation to provide an unambiguous expression for each fixed point $x$, using a single surreal parameter $u$ with $x = \Xi_{\mathbf{Fix_S}}(u)$. In a similar manner, one may regard the notation $\varepsilon_u$ as a way to disambiguate

$$\dot{\omega}^{\dot{\omega}^{\dot{\omega}^{\cdot^{\cdot^{\cdot}}}}}.$$

- If $\mathbf{S}$ is the interval $(-\omega, \omega)$, then we can see that $\Xi_\mathbf{S}$ fixes $\mathbf{No}^{\preccurlyeq}$ pointwise and replaces the initial segment $\omega$ (resp. $-\omega$) in the sign sequence of a positive (resp. negative) infinite number with $\omega - 1$ (resp. $1 - \omega$). Since $\omega/2 = \omega \dotplus (-\omega)$, we deduce that the defining isomorphism $\Xi_\mathbf{S}$ fixes $\mathbf{No}^{\preccurlyeq}$, $\mathbf{No}^{\sqsupseteq -\omega/2}$, and $\mathbf{No}^{\sqsupseteq \omega/2}$ pointwise. One can check that the class

$$\mathbf{Fix_S} = \mathbf{No}^{\sqsupseteq -\omega/2} \sqcup \mathbf{No}^{\preccurlyeq} \sqcup \mathbf{No}^{\sqsupseteq \omega/2}$$

is a surreal substructure.

In general, the class $\mathbf{Fix_S}$ may not be a surreal substructure. For instance, the class $\mathbf{Inc}$ defined in Example-4.14 satisfies $\forall x \in \mathbf{No}, \ell(\Xi_{\mathbf{Inc}} x) > \ell(x)$, and consequently has no fixed point. This raises the question of finding a condition on $\mathbf{S}$ that will ensure $\mathbf{Fix_S}$ to be a surreal substructure. One obvious first idea is to investigate when decreasing intersections of surreal substructures are surreal substructures.

## 5.2 Closed subclasses

We introduce a notion of closed subclasses $\mathbf{X}$ of an ambient surreal substructure $\mathbf{S} \supseteq \mathbf{X}$. In the case when $\mathbf{X}$ is a surreal substructure, we characterize its closedness in terms of its defining surreal isomorphism.

**Definition 5.4.** *Let $\mathbf{S}$ be a surreal substructure. Let $\mathbf{X}$ be a subclass of $\mathbf{S}$. We say that $\mathbf{X}$ is $\mathbf{S}$-closed, if the supremum in $(\mathbf{S}, \sqsubseteq)$ of any non-empty $\sqsubseteq$-chain in $\mathbf{X}$ lies in $\mathbf{X}$.*

**Example 5.5.**

- The intervals $(-\omega-1, \omega+1), (0,7)$ and $(0, \omega^2+1)$ are $\mathbf{No}$-closed convex surreal substructures. The interval $(-\omega, \omega)$ is a surreal substructure which is not $\mathbf{No}$-closed, since $\sup_\sqsubseteq \mathbb{N} = \omega \notin (-\omega, \omega)$.

- The structure $\mathbf{No}_>$ introduced in Example 5.3 is a non-convex $\mathbf{No}$-closed surreal substructure since having no different consecutive signs in one's sign sequence is preserved by taking suprema in $\mathbf{No}$.

- Likewise, the structure $2 \mathbin{\dot\times} \mathbf{No}$ is $\mathbf{No}$-closed.

- If $\mathbf{T}$ is a surreal substructure defined by the tree construction (see Proposition 4.13), then it is $\mathbf{No}$-closed if and only if for each non-empty $\sqsubseteq$-chain $X$ in $\mathbf{T}$, the element $\Xi_\mathbf{T} \sup_\sqsubseteq X$ of $\mathbf{T}$ is defined as $\sup_\sqsubseteq \Xi_\mathbf{T} X$. In particular, the surreal substructure $\mathbf{Inc}$ from Example 4.14 is not $\mathbf{No}$-closed.

- The class $\bigsqcup_{\alpha \in \mathbf{On}} \mathbf{No}^{\exists \alpha - 1}$ is $\mathbf{No}$-closed but has a proper class of $\sqsubseteq$-minimal elements $\{\alpha - 1 : \alpha \in \mathbf{On}_{\mathrm{lim}}\}$ (in particular, it has no root).

The term "closed" suggests the existence of a topology. Indeed, we have:

**Proposition 5.6.** *Let $\mathbf{S}$ be a surreal substructure. Arbitrary intersections and finite unions of $\mathbf{S}$-closed subclasses of $\mathbf{S}$ are $\mathbf{S}$-closed.*

**Proof.** It is clear that $\emptyset$ and $\mathbf{S}$ are $\mathbf{S}$-closed. Let $\mathbf{X}_\mathbf{I}$ be the intersection of a (possibly proper class-sized) non-empty family $(\mathbf{X}_i)_{i \in \mathbf{I}}$ of $\mathbf{S}$-closed subclasses of $\mathbf{S}$. Let $C$ be a non-empty $\sqsubseteq$-chain in $\mathbf{X}_\mathbf{I}$. We have $\sup_{\mathbf{S},\sqsubseteq} C \in \mathbf{X}_i$ for all $i \in \mathbf{I}$, whence $\sup_{\mathbf{S},\sqsubseteq} C \in \mathbf{X}_\mathbf{I}$ and $\mathbf{X}_\mathbf{I}$ is $\mathbf{S}$-closed.

Let $\mathbf{X}_1, \mathbf{X}_2$ be $\mathbf{S}$-closed subclasses of $\mathbf{S}$ and let $C$ be a non-empty $\sqsubseteq$-chain in $\mathbf{X}_1 \cup \mathbf{X}_2$. If $C$ admits a $\sqsubseteq$-maximum, then $\sup_{\mathbf{S},\sqsubseteq} C = \max C \in \mathbf{X}_1 \cup \mathbf{X}_2$. Otherwise, let $i \in \{1, 2\}$ be such that $C \cap \mathbf{X}_i$ is $\sqsubseteq$-cofinal in $C$. Then $\sup_{\mathbf{S},\sqsubseteq} C = \sup_{\mathbf{S},\sqsubseteq} C \cap \mathbf{X}_i \in \mathbf{X}_i \subseteq \mathbf{X}_1 \cup \mathbf{X}_2$, so $\mathbf{X}_1 \cup \mathbf{X}_2$ is $\mathbf{S}$-closed. $\square$

**Lemma 5.7.** *If $\mathbf{S}$ is a surreal substructure and $\mathbf{T}$ is a $\sqsubseteq$-final substructure of $\mathbf{S}$, then $\mathbf{T}$ is $\mathbf{S}$-closed.*

**Proof.** The class $\mathbf{T}$ is $\sqsubseteq$-final in $\mathbf{S}$, thus suprema of non-empty $\sqsubseteq$-chains in $\mathbf{T}$ lie in $\mathbf{T}$. $\square$

It will sometimes be useful to comprehend closure in terms of projections.

**Proposition 5.8.** *Let $\mathbf{S}$ be a surreal substructure. A rooted subclass $\mathbf{X}$ of $\mathbf{S}$ is $\mathbf{S}$-closed if and only if every element $x$ of $\mathbf{S}^{\exists \mathbf{X}^\bullet}$ has a $\sqsubseteq$-maximal initial segment $\mu_\mathbf{X}^\mathbf{S}(x)$ lying in $\mathbf{X}$.*

**Proof.** Assume that $\mathbf{X}$ is $\mathbf{S}$-closed. Consider $x \in \mathbf{S}$ with $\mathbf{X}^\bullet \sqsubseteq x$. Then the set of initial segments of $x$ lying in $\mathbf{X}$ is non-empty and closed under taking suprema in $\mathbf{S}$. Consequently, $x$ indeed admits a $\sqsubseteq$-maximal initial segment $\mu_\mathbf{X}^\mathbf{S}(x)$ in $\mathbf{X}$. Inversely, assume that $\mu_\mathbf{X}^\mathbf{S}$ is well defined on $\mathbf{S}^{\exists \mathbf{X}^\bullet}$ and let $C$ be a non-empty $\sqsubseteq$-chain in $\mathbf{X}$. If $C$ has a $\sqsubseteq$-maximum, then $\sup_{\mathbf{S},\sqsubseteq} C = \max_\sqsubseteq C \in \mathbf{X}$. Otherwise, $\mu_\mathbf{X}^\mathbf{S}(\sup_{\mathbf{S},\sqsubseteq} C) \not\sqsubset \sup_{\mathbf{S},\sqsubseteq} C$, so $\mu_\mathbf{X}^\mathbf{S}(\sup_{\mathbf{S},\sqsubseteq} C) = \sup_{\mathbf{S},\sqsubseteq} C \in \mathbf{X}$. This shows that $\mathbf{X}$ is $\mathbf{S}$-closed. $\square$

**Definition 5.9.** *If $\mathbf{X} \subseteq \mathbf{S}$ is rooted and $\mathbf{S}$-closed, then we define $\mu_\mathbf{X}^\mathbf{S}$ to be the function $\mathbf{S}^{\exists \mathbf{X}^\bullet} \twoheadrightarrow \mathbf{X}$ that sends each element $x$ of $\mathbf{S}^{\exists \mathbf{X}^\bullet}$ to the $\sqsubseteq$-maximal initial segment of $x$ that lies in $\mathbf{X}$. It is by definition surjective, $\sqsubseteq$-increasing, and satisfies the relation $\mu_\mathbf{X}^\mathbf{S} \circ \mu_\mathbf{X}^\mathbf{S} = \mu_\mathbf{X}^\mathbf{S}$. We call it the **topological projection** $\mathbf{S}^{\exists \mathbf{X}^\bullet} \to \mathbf{X}$.*

Since $\mu_\mathbf{S}^\mathbf{X}$ is $\sqsubseteq$-increasing when it exists, its fibers are $\sqsubseteq$-convex in $\mathbf{S}^{\exists \mathbf{X}^\bullet}$.

**Lemma 5.10.** *Let $\mathbf{T} \subseteq \mathbf{S}$ be surreal substructures and let $\mathbf{X} \subseteq \mathbf{T}$ be rooted. If $\mathbf{X}$ is $\mathbf{T}$-closed and $\mathbf{T}$ is $\mathbf{S}$-closed, then $\mathbf{X}$ is $\mathbf{S}$-closed, and we have $\mu_\mathbf{X}^\mathbf{S} \equiv \mu_\mathbf{X}^\mathbf{T} \circ \mu_\mathbf{T}^\mathbf{S}$ on $\mathbf{S}_0^{\exists \mathbf{X}^\bullet}$.*

**Proof.** Let $x \in \mathbf{S}^{\sqsupseteq \mathbf{X}^\bullet}$. Since $\mathbf{T}^\bullet \sqsubseteq \mathbf{T} \supseteq \mathbf{X}$, we have $\mathbf{T}^\bullet \sqsubseteq \mathbf{X}^\bullet$, whence $x \in \mathbf{S}^{\sqsupseteq \mathbf{T}^\bullet}$. The class $\mathbf{T}$ is $\mathbf{S}$-closed so $x$ has a maximal initial segment $\mu_\mathbf{T}^\mathbf{S}(x)$ lying in $\mathbf{T}$. Now $\mathbf{X}^\bullet$ is an initial segment of $x$ lying in $\mathbf{T}$, whence $\mathbf{X}^\bullet \sqsubseteq \mu_\mathbf{T}^\mathbf{S}(x)$. We may thus consider the maximal initial segment $\mu_\mathbf{X}^\mathbf{T}(\mu_\mathbf{T}^\mathbf{S}(x))$ of $\mu_\mathbf{T}^\mathbf{S}(x)$ that lies in $\mathbf{X}$. If $z \in \mathbf{X}$ is simpler than $x$, then $z \sqsubseteq \mu_\mathbf{T}^\mathbf{S}(x)$, since $z \in \mathbf{T}$. Similarly, $z \sqsubseteq \mu_\mathbf{X}^\mathbf{T}(\mu_\mathbf{T}^\mathbf{S}(x))$, since $z \in \mathbf{X}$. This proves that $\mu_\mathbf{X}^\mathbf{T}(\mu_\mathbf{T}^\mathbf{S}(x))$ is the maximal initial segment of $x$ lying in $\mathbf{X}$. □

We will mostly consider closures of surreal substructures in other ones. In this situation, closure can be regarded as a property of the defining surreal isomorphism:

**Lemma 5.11.** *If $\mathbf{T} \subseteq \mathbf{S}$ are surreal substructures, then $\mathbf{T}$ is $\mathbf{S}$-closed if and only if for any non-empty $\sqsubseteq$-chain $X$ of $\mathbf{No}$, we have $\Xi_\mathbf{T} \sup_\sqsubseteq X = \sup_{\mathbf{S}, \sqsubseteq} \Xi_\mathbf{T} X$.*

**Proof.** Assume that the relation holds. Let $Y$ be a non-empty $\sqsubseteq$-chain in $\mathbf{T}$ and consider the set $X = \Xi_\mathbf{T}^{-1}(Y)$. Since $\Xi_\mathbf{T}$ is an $\sqsubseteq$-embedding, the set $X$ is a non-empty $\sqsubseteq$-chain in $\mathbf{No}$, whence $\Xi_\mathbf{T} \sup_\sqsubseteq X = \sup_{\mathbf{T}, \sqsubseteq} \Xi_\mathbf{T} X = \sup_{\mathbf{T}, \sqsubseteq} Y$ (see Proposition 4.13). Our assumption on $\Xi_\mathbf{T}$ gives $\Xi_\mathbf{T} \sup_\sqsubseteq X = \sup_{\mathbf{S}, \sqsubseteq} \Xi_\mathbf{T} X = \sup_{\mathbf{S}, \sqsubseteq} Y$, so $\sup_{\mathbf{S}, \sqsubseteq} Y = \sup_{\mathbf{T}, \sqsubseteq} Y \in \mathbf{T}$, and $\mathbf{T}$ is $\mathbf{S}$-closed. Conversely, assume $\mathbf{T}$ is $\mathbf{S}$-closed. Let $X \subset \mathbf{No}$ be a non-empty $\sqsubseteq$-chain. Since $\Xi_\mathbf{T}$ is $\sqsubseteq$-increasing, the set $\Xi_\mathbf{T} X$ is a non-empty $\sqsubseteq$-chain in $\mathbf{T}$, so $\sup_{\mathbf{S}, \sqsubseteq} \Xi_\mathbf{T} X \in \mathbf{T}$, whence $\sup_{\mathbf{S}, \sqsubseteq} \Xi_\mathbf{T} X = \sup_{\mathbf{T}, \sqsubseteq} \Xi_\mathbf{T} X = \Xi_\mathbf{T} \sup_\sqsubseteq X$, which is the desired equality. □

**Lemma 5.12.** *Let $\mathbf{U}, \mathbf{V}, \mathbf{W}$ be surreal substructures.*

*a)* *If $\mathbf{V} \subseteq \mathbf{U}$, then $\mathbf{V}$ is $\mathbf{U}$-closed if and only if $\Xi_\mathbf{V}$ sends $\mathbf{No}$-closed subclasses of $\mathbf{No}$ onto $\mathbf{U}$-closed subclasses of $\mathbf{U}$.*

*b)* *If $\mathbf{V}$ and $\mathbf{W}$ are $\mathbf{No}$-closed, then so is $\mathbf{V} \prec \mathbf{W}$.*

*c)* *If $\mathbf{V}$ and $\mathbf{V} \prec \mathbf{W}$ are $\mathbf{No}$-closed, then so is $\mathbf{W}$.*

**Proof.** *a)* Assume $\mathbf{V}$ is $\mathbf{U}$-closed and $\mathbf{X}$ is a closed subclass of $\mathbf{No}$. Let $Y$ be a non-empty $\sqsubseteq$-chain in $\Xi_\mathbf{V} \mathbf{X}$. The set $\Xi_\mathbf{V}^{-1}(Y)$ is a non-empty $\sqsubseteq$-chain in $\mathbf{X}$ so its supremum lies in $\mathbf{X}$, and $\Xi_\mathbf{V} \mathbf{X} \ni \Xi_\mathbf{V} \sup_\sqsubseteq \Xi_\mathbf{V}^{-1} Y = \sup_{\mathbf{U}, \sqsubseteq} \Xi_\mathbf{V} \Xi_\mathbf{V}^{-1}(Y) = \sup_{\mathbf{U}, \sqsubseteq} Y$, so $\Xi_\mathbf{V} \mathbf{X}$ is $\mathbf{U}$-closed. Conversely, if $\Xi_\mathbf{V}$ sends closed classes of surreal numbers onto $\mathbf{U}$-closed subclasses of $\mathbf{U}$, then in particular $\mathbf{V} = \Xi_\mathbf{V} \mathbf{No}$ is $\mathbf{U}$-closed.

*b)* This is a direct consequence of *a)*.

*c)* Assume that $\mathbf{V} \prec \mathbf{W}$ and $\mathbf{V}$ are $\mathbf{No}$-closed. Let $X$ be a non-empty $\sqsubseteq$-chain in $\mathbf{No}$. Then $\Xi_\mathbf{V} \Xi_\mathbf{W} \sup_\sqsubseteq X = \sup_\sqsubseteq \Xi_\mathbf{V} \Xi_\mathbf{W} X = \Xi_\mathbf{V} \sup_\sqsubseteq \Xi_\mathbf{W} X$, and since $\Xi_\mathbf{V}$ is injective, we get $\Xi_\mathbf{W} \sup_\sqsubseteq X = \sup_\sqsubseteq \Xi_\mathbf{W} X$, so $\mathbf{W}$ is $\mathbf{No}$-closed. □

We now come to the main interest of the notion of closure.

**Proposition 5.13.** *Let $0 < \alpha$ be a limit ordinal. Let $\mathbf{S}$ be a surreal substructure and let $(\mathbf{S}_\beta)_{\beta < \alpha}$ be a decreasing sequence of $\mathbf{S}$-closed surreal substructures of $\mathbf{S}$. Then its intersection $\bigcap_{\beta < \alpha} \mathbf{S}_\beta$ is an $\mathbf{S}$-closed surreal substructure.*

**Proof.** We use the characterization of surreal substructures given in Proposition 4.13. By Proposition 5.6, the class $\mathbf{S}_\alpha := \bigcap_{\beta < \alpha} \mathbf{S}_\beta$ is $\mathbf{S}$-closed. In particular, the class $\mathbf{S}_\alpha$ has suprema of non-empty $\sqsubseteq$-chains. We also have $\sup_{\sqsubseteq, \mathbf{S}_\alpha} \emptyset = \mathbf{S}_\alpha^\bullet = \sup_{\sqsubseteq, \mathbf{S}} \{\mathbf{S}_\beta^\bullet : \beta < \alpha\}$ which lies in $\mathbf{S}_\alpha$ by the $\mathbf{S}$-closure of each structure $\mathbf{S}_\beta$ for $\beta < \alpha$, so the empty $\sqsubseteq$-chain has a supremum as well. □

Let us now treat the case of left and right successors. Given $u \in \mathbf{S}_\alpha$, let $u_{\beta,-1} < u$ and $u_{\beta,1} > u$ be the left and right successors of $u$ in $\mathbf{S}_\beta$, for each ordinal $\beta < \alpha$. For $\beta < \gamma < \alpha$, we have $u_{\gamma,-1} \in \mathbf{S}_\beta$ and $u_{\gamma,-1} < u$, so $u_{\beta,-1} \sqsubseteq u_{\gamma,-1}$ by the definition of left successors. Similarly, we get $u_{\beta,1} \sqsubseteq u_{\gamma,1}$. Thus the sets $\{u_{\beta,-1} : \beta < \alpha\}$ and $\{u_{\beta,1} : \beta < \alpha\}$ are $\sqsubseteq$-chains whose suprema $u_{-1}, u_1$ in $\mathbf{S}$ satisfy $u_{-1} < u < u_1$. For $v \in \mathbf{S}_\alpha$ with $u < v$ and $\beta < \alpha$, we have $u, v \in \mathbf{S}_\beta$ so $u_{\beta,1} \sqsubseteq v$, whence $u_1 \sqsubseteq v$. This means that $u_1$ is the right successor of $u$ in $\mathbf{S}_\alpha$. Likewise, $u_{-1}$ is the left successor of $u$ in $\mathbf{S}_\alpha$. We conclude that $\mathbf{S}_\alpha$ is a surreal substructure. □

**Corollary 5.14.** *If the surreal substructure $\mathbf{S}$ is $\mathbf{No}$-closed, then $\mathbf{Fix_S}$ is an $\mathbf{No}$-closed surreal substructure.*

**Proof.** This is a direct consequence of Lemma 5.12, Proposition 5.13 and Proposition 5.2. □

**Remark 5.15.** Corollary 5.14 is similar to [34, Theorem 8.2]. Lurie's result is more general, but when applied to an $\mathbf{No}$-closed surreal substructure $\mathbf{S}$, it only concludes that $\mathbf{Fix_S}$ is a "good tree". Good trees need not be surreal substructures. For instance,

$$\mathbf{No}^{\exists -2} \sqcup \mathbf{No}^{\exists -1/2}\{0\} \sqcup \mathbf{No}^{\exists 1/2} \sqcup \mathbf{No}^{\exists 2}$$

is a good tree, but not a surreal substructure, since 0 has two right successors and two left successors.

## 5.3 Transfinite right-imbrications of surreal substructures

The class of $\mathbf{No}$-closed surreal substructures being closed under decreasing intersections, we are now in a position to define a notion of transfinite right-imbrications of $\mathbf{No}$-closed surreal substructures.

**Theorem 5.16.** *Let $\alpha$ be an ordinal. Let $\mathbf{U} = (\mathbf{U}_\beta)_{\beta < \alpha}$ be a sequence of $\mathbf{No}$-closed surreal substructures. We define a sequence $(\prec_{\gamma < \beta} \mathbf{U}_\gamma)_{\beta \leqslant \alpha}$ of $\mathbf{No}$-closed surreal substructures by the following rules:*

- $\prec_{\gamma < 0} \mathbf{U}_\gamma = \mathbf{No}$.
- $\prec_{\gamma < \beta+1} \mathbf{U}_\gamma = (\prec_{\gamma < \beta} \mathbf{U}_\gamma) \prec \mathbf{U}_\beta$ *if $\beta < \alpha$,*
- $\prec_{\gamma < \beta} \mathbf{U}_\gamma = \bigcap_{\beta' < \beta} \prec_{\gamma < \beta'} \mathbf{U}_\gamma$ *if $0 < \beta \leqslant \alpha$ is limit.*

*Then each class $\prec_{\gamma < \beta} \mathbf{U}_\gamma$ for $\beta \leqslant \alpha$ is an $\mathbf{No}$-closed surreal substructure, and if $\beta \dotplus \delta \leqslant \alpha$, then we have*

$$\underset{\gamma < \beta \dotplus \delta}{\prec} \mathbf{U}_\gamma = \underset{\gamma < \beta}{\prec} \mathbf{U}_\gamma \prec \underset{\beta \leqslant \gamma < \beta \dotplus \delta}{\prec} \mathbf{U}_\gamma. \tag{5.1}$$

**Proof.** We first need to prove that the definition is warranted. We do this by transfinite induction, while proving at the same time that the sequence $(\prec_{\gamma < \beta} \mathbf{U}_\gamma)_{\beta \leqslant \alpha}$ is decreasing, and that each term is an $\mathbf{No}$-closed surreal substructure. Let $\beta \leqslant \alpha$ be such that these assumptions hold strictly below $\beta$. If $\beta = \beta' + 1$ is a successor ordinal, then $\prec_{\gamma < \beta'} \mathbf{U}_\gamma$ and $\mathbf{U}_{\beta'}$ are $\mathbf{No}$-closed surreal substructures, whence $\prec_{\gamma < \beta} \mathbf{U}_\gamma := (\prec_{\gamma < \beta'} \mathbf{U}_\gamma) \prec \mathbf{U}_{\beta'}$ is well defined and $\mathbf{No}$-closed (by Lemma 5.12). The surreal substructure $\prec_{\gamma < \beta} \mathbf{U}_\gamma$ is a left factor of $\prec_{\gamma < \beta} \mathbf{U}_\gamma$, which implies that $\prec_{\gamma < \beta} \mathbf{U}_\gamma \subseteq \prec_{\gamma < \beta'} \mathbf{U}_\gamma$. If $\beta$ is limit, the intersection that defines $\prec_{\gamma < \beta} \mathbf{U}_\gamma$ is an $\mathbf{No}$-closed surreal substructure by Proposition 5.13, and $(\prec_{\gamma < \beta'} \mathbf{U}_\gamma)_{\beta' \leqslant \beta}$ is clearly decreasing.

We prove the identity (5.1) by induction on $\beta \dotplus \delta$. Let $\sigma$ be an ordinal such that (5.1) holds for any sequence $\mathbf{U}$ and $\beta, \delta$ with $\beta \dotplus \delta < \sigma$. Let $\beta, \delta$ be such that $\beta \dotplus \delta = \sigma$. If $\delta = \eta + 1$ for some ordinal $\eta$, then

$$\begin{aligned}
\prec_{\gamma < \beta \dotplus \delta} \mathbf{U}_\gamma &= \prec_{\gamma < \beta \dotplus \eta + 1} \mathbf{U}_\gamma \\
&= \prec_{\gamma < \beta \dotplus \eta} \mathbf{U}_\gamma \prec \mathbf{U}_{\beta \dotplus \eta} \\
&= \prec_{\gamma < \beta} \mathbf{U}_\gamma \prec \prec_{\beta \leqslant \gamma < \beta \dotplus \eta} \mathbf{U}_\gamma \prec \mathbf{U}_{\beta \dotplus \eta} \\
&= \prec_{\gamma < \beta} \mathbf{U}_\gamma \prec \prec_{\beta \leqslant \gamma < \beta \dotplus \eta + 1} \mathbf{U}_\gamma \\
&= \prec_{\gamma < \beta} \mathbf{U}_\gamma \prec \prec_{\beta \leqslant \gamma < \beta \dotplus \delta} \mathbf{U}_\gamma.
\end{aligned}$$

If $\delta$ is limit, then we have

$$\begin{aligned}
\prec_{\gamma < \beta \dotplus \delta} \mathbf{U}_\gamma &= \bigcap_{\eta < \delta} \prec_{\gamma < \beta \dotplus \eta} \mathbf{U}_\gamma \\
&= \bigcap_{\eta < \delta} \left( \prec_{\gamma < \beta} \mathbf{U}_\gamma \prec \prec_{\beta \leqslant \gamma < \beta \dotplus \eta} \mathbf{U}_\gamma \right) \\
&= \bigcap_{\eta < \delta} \Xi_{\prec_{\gamma < \beta} \mathbf{U}_\gamma} \left( \prec_{\beta \leqslant \gamma < \beta \dotplus \eta} \mathbf{U}_\gamma \right) \\
&= \Xi_{\prec_{\gamma < \beta} \mathbf{U}_\gamma} \left( \bigcap_{\eta < \delta} \prec_{\beta \leqslant \gamma < \beta \dotplus \eta} \mathbf{U}_\gamma \right) \\
&= \Xi_{\prec_{\gamma < \beta} \mathbf{U}_\gamma} \left( \prec_{\beta \leqslant \gamma < \beta \dotplus \delta} \mathbf{U}_\gamma \right) \\
&= \prec_{\gamma < \beta} \mathbf{U}_\gamma \prec \prec_{\beta \leqslant \gamma < \beta \dotplus \delta} \mathbf{U}_\gamma.
\end{aligned}$$

(The injectivity of $\Xi_{\prec_{\gamma < \beta} \mathbf{U}_\gamma}$ allowed us to move it through intersections). □

**Example 5.17.** In [14, p 34–35], Conway informally discussed continued exponential expressions of the form

$$x = u_0 \pm \omega^{u_1 \pm \omega^{u_2 \pm \omega^{\cdot^{\cdot^{\cdot}}}}}.$$

He outlined an approach for proving that the class of numbers that can be expressed in this way is order isomorphic to $\mathbf{No}$. Conway's ideas were rigorously worked out by Lemire [31, 32, 33]. He first proved the following result in the case when $u_0 = u_1 = \cdots = 0$: given $(z_i)_{i \in \mathbb{N}} \in \{-1, 1\}^{\mathbb{N}}$, let $\mathbf{E}_z$ be the class of numbers $x$ such that there exists a sequence $(x_i)_{i \in \mathbb{N}} \in \mathbf{No}^{\mathbb{N}}$ with

$$x = z_0 \omega^{z_1 \omega^{\cdot^{\cdot^{z_i x_i}}}}$$

for all $i \in \mathbb{N}$. Then $\mathbf{E}_z$ is order isomorphic to $\mathbf{No}$. Moreover, writing $\varphi_z : \mathbf{No} \longrightarrow \mathbf{E}_z$ for the isomorphism, $\varphi_z$ has fixed points of any order $\alpha \in \mathbf{On}$, and the class $\mathbf{E}_z^\alpha$ of such fixed points is also order isomorphic to $\mathbf{No}$. This result follows from Theorem 5.16 by taking $\mathbf{U}_{\omega \alpha + i} = z_i \mathbf{Mo}$ for all $\alpha \in \mathbf{On}$ and $i < \omega$. Then $\mathbf{E}_z^\alpha = \prec_{\beta < \omega^{1+\alpha}} \mathbf{U}_\beta$ for all $\alpha \in \mathbf{On}$.

A similar result was proved by Lemire for more general continued exponential expressions [32, Theorem 4]. This result is more involved and presents similarities with our results about nested expansions in section 8 below.

**Proposition 5.18.** *Let $\mathbf{S}$ be an $\mathbf{No}$-closed surreal substructure. For each ordinal $\alpha$, let*

$$\mathbf{S}^{\prec \alpha} := \prec_{\beta < \alpha} \mathbf{S}.$$

*Each $\mathbf{S}^{\prec \alpha}$ is an $\mathbf{No}$-closed surreal substructure, and for $\alpha, \beta \in \mathbf{On}$, we have:*

$$\mathbf{S}^{\prec(\alpha \dotplus \beta)} = \mathbf{S}^{\prec \alpha} \prec \mathbf{S}^{\prec \beta}. \tag{5.2}$$
$$\mathbf{S}^{\prec(\alpha \dottimes \beta)} = (\mathbf{S}^{\prec \alpha})^{\prec \beta}. \tag{5.3}$$

**Proof.** Most of this is a direct consequence of Theorem 5.16; we only need to prove the identity (5.3). Let $\pi \in \mathbf{On}$ be such that this identity holds for $\alpha \dottimes \beta < \pi$. Let $\alpha, \beta$ be ordinal numbers with $\alpha \dottimes \beta = \pi$. Corollary 5.14 justifies that the same construction can be applied to the structure $\mathbf{S}^{\prec \alpha}$. If $\beta = \eta + 1$ for $\eta \in \mathbf{On}$, then we have

$$\begin{aligned}
(\mathbf{S}^{\prec \alpha})^{\prec \beta} &= (\mathbf{S}^{\prec \alpha})^{\prec \eta} \prec \mathbf{S}^{\prec \alpha} \\
&= \mathbf{S}^{\prec(\alpha \dottimes \eta)} \prec \mathbf{S}^{\prec \alpha} \\
&= \mathbf{S}^{\prec(\alpha \dottimes \eta \dotplus \alpha)} \\
&= \mathbf{S}^{\prec(\alpha \dottimes \beta)},
\end{aligned}$$

where we used (5.2) as well as the inductive hypothesis. If $\beta$ is limit, then

$$\begin{aligned}
(\mathbf{S}^{\prec \alpha})^{\prec \beta} &= \bigcap_{\eta < \beta} (\mathbf{S}^{\prec \alpha})^{\prec \eta} \\
&= \bigcap_{\eta < \beta} \mathbf{S}^{\prec(\alpha \dottimes \eta)} \\
&= \bigcap_{\gamma < \alpha \dottimes \beta} \mathbf{S}^{\prec \gamma} \\
&= \mathbf{S}^{\prec(\alpha \dottimes \beta)}.
\end{aligned} \qquad \square$$

Note that for $n \in \mathbb{N}$, the structure $\mathbf{S}^{\prec n}$ is the $n$-fold imbrication of $\mathbf{S}$ into itself, and we have $\Xi_{\mathbf{S}^{\prec n}} = (\Xi_\mathbf{S})^n$. For $\alpha \in \mathbf{On}$, we have $\mathbf{S}^{\prec(\omega \dottimes \alpha)} = \mathbf{Fix}_\mathbf{S}^{\prec \alpha}$, by Proposition 5.2 and the identity (5.3). Thus transfinite right-imbrications of $\mathbf{S}$ with itself allow us to define higher order fixed points of $\Xi_\mathbf{S}$ as being elements of the $\mathbf{S}^{\prec \dot\omega^\alpha}$ with $0 < \alpha \in \mathbf{On}$. As we have seen, imbrication is left-distributive on decreasing intersections that form a surreal substructure. It is not right-distributive in general. For instance if $\mathbf{S}$ is a proper $\mathbf{No}$-closed surreal substructure of $\mathbf{No}$, then $\mathbf{S}^{\prec \omega} \prec \mathbf{S}$ is a proper subclass of $\mathbf{S}^{\prec \omega} = \bigcap_{n \in \mathbb{N}} (\mathbf{S}^{\prec n} \prec \mathbf{S})$.

**Example 5.19.** We will see in section 7.2 that the class $\mathbf{No}_{>}^{\prec \alpha}$ coincides with $\dot\omega^\alpha \dottimes \mathbf{No}$.

**Example 5.20.** The class $\mathbf{Mo}^{\prec \omega}$ of fixed points of the $\omega$-map was studied before in [14, 24, 34]; numbers in $\mathbf{Mo}^{\prec \omega}$ are called *generalized $\varepsilon$-numbers*. It also comes up in the study of the exponential function and the length of sign sequences [15, 30]. The class $\mathbf{Mo}^{\prec \dot\omega^\alpha}$ corresponds to a higher order fixed points of the $\omega$-map and we expect it to play a similar role as $\mathbf{Mo}^{\prec \omega}$ for the study of the $\alpha$-th hyperexponential function.

# 6 Convex partitions

Throughout this section, **S** stands for a surreal substructure.

## 6.1 Convex partitions

**Definition 6.1.** *Let $\Pi$ be a partition of **S** into convex subclasses. We say that $\Pi$ is a **convex partition** of **S**. For $x \in \mathbf{S}$ we let $\Pi[x]$ denote the member of $\Pi$ containing $x$ and recall that this class is rooted (by Lemma 4.16). We say that $x \in \mathbf{S}$ is $\Pi$-simple if $x = \Pi[x]^\bullet$, and we let $\mathbf{Smp}_\Pi$ denote the class of $\Pi$-simple elements of **S**. For $x, y \in \mathbf{S}$ we write:*

$$x =_\Pi y \quad \text{if } \Pi[x] = \Pi[y],$$
$$x <_\Pi y \quad \text{if } \Pi[x] < \Pi[y],$$
$$x \leqslant_\Pi y \quad \text{if } \Pi[x] = \Pi[y] \text{ or } \Pi[x] < \Pi[y].$$

**Remark 6.2.** Convex partitions are sometime called *condensations* [36, Definition 4.1].

We can obtain **S** as $\mathbf{Smp}_{\Pi_{\mathrm{disc}}}$ through the *discrete partition* $\Pi_{\mathrm{disc}}$ with $\Pi_{\mathrm{disc}}[x] = \{x\}$ for all $x \in \mathbf{S}$. Let $\pi_\Pi(x) := \Pi[x]^\bullet \in \mathbf{S}$ for all $x \in \mathbf{S}$. The map $\pi_\Pi \colon \mathbf{S} \longrightarrow \mathbf{Smp}_\Pi$ is a surjective, increasing projection. We refer to it as the $\Pi$-*simple projection*.

For the remainder of this subsection, let $\Pi$ be a convex partition of **S**. A *quasi-order* (or *preorder*) is a binary relation that is reflexive and transitive. The following lemma states basic facts on partitions of a linear order into convex subclasses.

**Lemma 6.3.** *The relation $\leqslant_\Pi$ is a linear quasi-order and restricts to a linear order on $\mathbf{Smp}_\Pi$. For $x, y \in \mathbf{S}$, we have $x \leqslant_\Pi y$ if and only if $\pi_\Pi(x) \leqslant \pi_\Pi(y)$.*

**Proof.** It is well known that the partition $\Pi$ corresponds to the equivalence relation $=_\Pi$ on **S**. The transitivity and irreflexivity of $<_\Pi$ follow from that of $<$ on subclasses of **No**. That its restriction to $\mathbf{Smp}_\Pi$ is a linear order is a direct consequence of the definition of $\mathbf{Smp}_\Pi$ and the equivalence stated above, which we now prove. If $\Pi$ has only one member, then the result is trivial. Else let $x, y \in \mathbf{S}$ with $x <_\Pi y$. We have $\pi_\Pi(x) \in \Pi[x] < \Pi[y] \ni \pi_\Pi(y)$ so $\pi_\Pi(x) < \pi_\Pi(y)$. Conversely, assume that $\pi_\Pi(x) < \pi_\Pi(y)$. Then $\Pi[x] \neq \Pi[y]$ which since $\Pi$ is a partition implies that $\Pi[x] \cap \Pi[y] = \emptyset$. For $x' \in \Pi[x]$, there may be no element $z$ of $\Pi[y]$ such that $z \leqslant x$ for this would imply $z \leqslant x \leqslant \pi_\Pi(y)$ whence $x \in \Pi[y]$ by convexity of this class: a contradiction. We thus have $\Pi[x] < \Pi[y]$, that is, $x <_\Pi y$. By definition of $\pi_\Pi$, the relation $x =_\Pi y$ implies that $\pi_\Pi(x) = \pi_\Pi(y)$, whereas $\pi_\Pi(x) = \pi_\Pi(y)$ implies that $\Pi[x] \cap \Pi[y] \neq \emptyset$, so $\Pi[x] = \Pi[y]$, so $x =_\Pi y$. □

For any subclass **X** of **S**, we let $\Pi[\mathbf{X}]$ denote the class $\bigcup_{x \in \mathbf{X}} \Pi[x]$.

**Lemma 6.4.** *Let $\mathbf{A}, \mathbf{B}$ be subclasses of **S**. Then the following statements are equivalent:*

*a)* $\mathbf{A} < \Pi[\mathbf{B}]$.

*b)* $\Pi[\mathbf{A}] < \mathbf{B}$.

*c)* $\Pi[\mathbf{A}] < \Pi[\mathbf{B}]$.

**Proof.** All inequalities are vacuously true if $\mathbf{A}=\emptyset$ or $\mathbf{B}=\emptyset$. Assume that $\mathbf{A}$ and $\mathbf{B}$ are non-empty and let $a \in \mathbf{A}$ and $b \in \mathbf{B}$. Assume for contradiction that $\mathbf{A}<\Pi[\mathbf{B}]$, but $\Pi[\mathbf{A}] \not< \Pi[\mathbf{B}]$. Then there exist $a' \in \Pi[a]$ and $b' \in \Pi[b]$ with $a < b' \leqslant a'$. By convexity of $\Pi[a]$, this yields $b' \in \Pi[a]$, whence $a \in \Pi[b]$. This contradiction shows that $\mathbf{A}<\Pi[\mathbf{B}] \Longrightarrow \Pi[\mathbf{A}]<\Pi[\mathbf{B}]$. The inverse implication clearly holds. The equivalence $\Pi[\mathbf{A}]<\mathbf{B} \Longleftrightarrow \Pi[\mathbf{A}]<\Pi[\mathbf{B}]$ holds for similar reasons. $\square$

**Lemma 6.5.** *For $x \in \mathbf{S}$, the three following statements are equivalent:*

*a)* $x$ *is $\Pi$-simple.*

*b)* *There is a cut representation $(L,R)$ of $x$ in $\mathbf{S}$ such that $\Pi[L]<x<\Pi[R]$.*

*c)* $\Pi[x_L^{\mathbf{S}}] < x < \Pi[x_R^{\mathbf{S}}]$.

**Proof.** Since $(x_L^{\mathbf{S}}, x_R^{\mathbf{S}})$ is a cut representation of $x$ in $\mathbf{S}$, the assertion $c)$ implies $b)$.

Conversely, if $(L,R)$ is a cut representation of $x$ in $\mathbf{S}$ with $\Pi[L]<x<\Pi[R]$, then we have $L<\Pi[x]<R$ by the previous lemma. By Proposition 4.11(b), the cut representation $(L,R)$ is cofinal with respect to $(x_L^{\mathbf{S}}, x_R^{\mathbf{S}})$, so $x_L^{\mathbf{S}}<\Pi[x]<x_R^{\mathbf{S}}$. Hence $\Pi[x_L^{\mathbf{S}}]<x<\Pi[x_R^{\mathbf{S}}]$, again by Lemma 6.4. This shows that $b)$ implies $c)$.

Assume now that $x$ is $\Pi$-simple and let us prove $c)$. For $u \in x_L^{\mathbf{S}}$, we have $u \sqsubset x$, so $u \notin \Pi[x]$, whence $u \neq_\Pi x$. We do not have $\Pi[x]<\Pi[u]$ since $x \not< u$, so Lemma 6.3 yields $\Pi[u]<\Pi[x]$, and in particular $\Pi[u]<x$. This proves that $\Pi[x_L^{\mathbf{S}}]<x$, and similar arguments yield $x<\Pi[x_R^{\mathbf{S}}]$.

Assume finally that $c)$ holds and let us prove $a)$. We have $\Pi[x]^\bullet \sqsubseteq x$ so $\Pi[x]^\bullet \in x_L^{\mathbf{S}} \cup \{x\} \cup x_R^{\mathbf{S}}$. Now the class $\Pi[\Pi[x]^\bullet] = \Pi[x]$ is neither strictly greater nor strictly lower than $x$, so our assumption imposes $\Pi[x]^\bullet = x$. We conclude that $x$ is $\Pi$-simple. $\square$

An order $\leqslant$ on a set $S$ is said to be *dense* if for any $a,b \in S$ with $a<b$, there exists a $c \in S$ with $a<c<b$.

**Proposition 6.6.** *Assume that $\mathbf{Smp}_\Pi$ is dense. Then $\Pi$ is the unique convex partition of $\mathbf{S}$ such that $\mathbf{Smp}_\Pi$ is the class of $\Pi$-simple elements of $\mathbf{S}$.*

**Proof.** For $a \in \mathbf{Smp}_\Pi$, let $\mathbf{A}_a$ denote the class of elements $x$ of $\mathbf{S}$ such that no $\Pi$-simple element lies strictly between $a$ and $x$. The definition of the family $(\mathbf{A}_b)_{b \in \mathbf{Smp}_\Pi}$ only depends on the class $\mathbf{Smp}_\Pi$, and not specifically on $\Pi$. For $a \in \mathbf{Smp}_\Pi$, we have $\Pi[a] \subseteq \mathbf{A}_a$.

Conversely, let $x \in \mathbf{A}_a$, and assume for contradiction that $x$ lies outside of $\Pi[a]$, say $a <_\Pi x$. Then $a <_\Pi \pi_\Pi(x)$ and, $\mathbf{Smp}_\Pi$ being dense, there exists a $\Pi$-simple element $b$ between $a$ and $\pi_\Pi(x)$. But $a <_\Pi b <_\Pi \pi_\Pi(x)$ implies $a<b<x$, which contradicts the assumption that there is no simple element between $a$ and $x$. We conclude that $\Pi[a] = \mathbf{A}_a$, which entails in particular that the partition $\Pi$ is uniquely determined by $\mathbf{Smp}_\Pi$. $\square$

If $\mathbf{Smp}_\Pi$ is dense, then we call $\Pi$ the *defining partition* of $\mathbf{Smp}_\Pi$. Notice that this is in particular the case when $\mathbf{Smp}_\Pi$ is a surreal substructure. We next consider a set-theoretic condition under which $\mathbf{Smp}_\Pi$ is always a surreal substructure.

We say that $\Pi$ is *thin* if each member of $\Pi$ has a cofinal and coinitial subset. For instance, the convex partition $\Pi$ of $\mathbf{No}$ where

$$\Pi[x] := \{y \in \mathbf{No} : \exists n \in \mathbb{N}, -n < x - y < n\},$$

is thin. Indeed each class $\Pi[x]$ for $x \in \mathbf{No}$ admits the cofinal and coinitial subset $x + \mathbb{Z}$. See Example 6.15 below for more (counter)examples of thin convex partitions. If $\Pi$ is thin, then we may pick a distinguished family $(\Pi[x])_{x \in \mathbf{S}}$ such that each $\Pi[x]$ for $x \in \mathbf{S}$ is a cofinal and coinitial *subset* of $\Pi[x]$, with $\Pi[x] = \Pi[y] \iff x =_\Pi y$. We write $\Pi[\mathbf{X}] = \bigcup_{x \in \mathbf{X}} \Pi[x]$ for any subclass $\mathbf{X}$ of $\mathbf{S}$.

**Theorem 6.7.** *If $\Pi$ is thin, then $\mathbf{Smp}_\Pi$ is a surreal substructure. If $(L, R)$ is a cut representation in $\mathbf{Smp}_\Pi$, then we have*

$$\{L \mid R\}_{\mathbf{Smp}_\Pi} = \{\Pi[L] \mid \Pi[R]\}_\mathbf{S}.$$

**Proof.** Let $L < R$ be subsets of $\mathbf{Smp}_\Pi$. For $l \in L$ and $r \in R$, we have $\Pi[l] < \Pi[r]$ by Lemma 6.3. Therefore $\Pi[l] < \Pi[r]$ holds as well, which means that $x := \{\Pi[L] \mid \Pi[R]\}_\mathbf{S}$ is well defined. Given $l \in L$ and $l' \in \Pi[l]$, there exists an $l'' \in \Pi[l]$ with $l'' > l'$, since $\Pi[l]$ is cofinal in $\Pi[l]$. It follows that $l' < l'' < x$, whence $\Pi[l] < x$. A similar reasoning shows that $x < \Pi[r]$ for any $r \in R$. By Lemma 6.5, it follows that $x$ is $\Pi$-simple. Let $y \in (L \mid R)_\mathbf{S}$ be $\Pi$-simple. Given $l \in L$ and $r \in R$, the $\Pi$-simplicity of $l, r,$ and $y$ implies that $\Pi[l] < y < \Pi[r]$, and in particular that $\Pi[l] < y < \Pi[r]$. We deduce that $x \sqsubseteq y$, so $x = \{L \mid R\}_{\mathbf{Smp}_\Pi}$. By Proposition 4.7, we conclude that the class $\mathbf{Smp}_\Pi$ is a surreal substructure. □

**Remark 6.8.** The above theorem can be regarded as a strengthening of [34, Theorem 8.4] in a different framework. Indeed, Lurie's result is restricted to the case when $\mathbf{S} = \mathbf{No}$ and requires the additional assumption that

$$\forall a, b, c \in \mathbf{S}, ((a \sqsubseteq b \sqsubseteq c \land \Pi[a] = \Pi[c]) \implies (\Pi[a] = \Pi[b] = \Pi[c])).$$

This condition is equivalent to the condition that $\Pi$ be sharp in our terminology (see below); it fails for the partition $\Pi$ of $\mathbf{No}^{>,>}$ such that

$$\forall a \in \mathbf{No}^{>,>}, \Pi[a] := \{b \in \mathbf{No}^{>,>} : \exists n \in \mathbb{N}, \log_n(b) \asymp \log_n(a)\},$$

which is the defining convex partition of the set $\mathbf{La}$ of log-atomic numbers. Indeed, we have $\omega \sqsubseteq \omega^\omega \sqsubseteq \omega^{\omega^{\omega^{-1}}}$, where $\omega^\omega = \lambda_1 > \Pi[\lambda_0] = \Pi[\omega]$, but

$$\omega^{\omega^{\omega^{-1}}} = \exp\left(\omega^{\frac{2}{\omega}}\right) = \exp_2(2 \log_2 \omega) \in \Pi[\omega].$$

Still, $\mathbf{La}$ is a surreal substructure and even an $\mathbf{No}$-closed one.

When $\Pi$ is thin, the structure $\mathbf{Smp}_\Pi$ is in addition cofinal and coinitial in $\mathbf{S}$, since for $x \in \mathbf{S}$, we have $\mathbf{Smp}_\Pi \ni \{\emptyset \mid \Pi[x]\}_\mathbf{S} \leqslant x \leqslant \{\Pi[x] \mid \emptyset\}_\mathbf{S} \in \mathbf{Smp}_\Pi$. By the previous proposition, we may say that $\mathbf{Smp}_\Pi$ is thin if its defining partition $\Pi$ is thin. If $\Pi$ is not thin, then $\mathbf{Smp}_\Pi$ may fail to be a surreal substructure, but one can prove that there exists a unique $\sqsubseteq$-initial subclass $\mathbf{I}$ of $\mathbf{No}$ and a unique isomorphism between $\mathbf{Smp}_\Pi$ and $\mathbf{I}$.

For instance, we can obtain the ring $\mathbf{Oz} := \mathbf{No}_{\succ} + \mathbb{Z}$ of *omnific integers* of [14, Chapter 5] as $\mathbf{Smp}_{\Pi_{\mathbf{Oz}}}$ where for each number $z \in \mathbf{Oz}$, we set $\Pi_{\mathbf{Oz}}[z] := [z, z+1)$. This is not a surreal substructure since the cut $(0 \mid 1)_{\mathbf{Oz}}$ is empty. Nevertheless, $\mathbf{Oz}$ is $\sqsubseteq$-initial in $\mathbf{No}$. Note that different partitions may yield the same class $\mathbf{Oz}$ (for instance replacing $\Pi_{\mathbf{Oz}}[0]$ and $\Pi_{\mathbf{Oz}}[1]$ with $[0, 1/2)$ and $[1/2, 2)$ respectively and leaving the other classes unchanged), in contrast to the case of dense partitions from Proposition 6.6. The partition $\Pi_2$ in Example 6.15 below is not thin and yet $\mathbf{Smp}_{\Pi_2}$ is a surreal substructure.

**Proposition 6.9.** *Assume that $\Pi$ is thin. Then we have the following uniform cut equation for $\Xi_{\mathbf{Smp}_\Pi}$ and $x \in \mathbf{No}$:*

$$\Xi_{\mathbf{Smp}_\Pi} x = \{\Pi[\Xi_{\mathbf{Smp}_\Pi} x_L] \mid \Pi[\Xi_{\mathbf{Smp}_\Pi} x_R]\}_{\mathbf{S}}.$$

**Proof.** The cut equation follows from Theorem 6.7 and the relation

$$\Xi_{\mathbf{Smp}_\Pi} x = \{\Xi_{\mathbf{Smp}_\Pi} x_L \mid \Xi_{\mathbf{Smp}_\Pi} x_R\}_{\mathbf{Smp}_\Pi}.$$

Now towards uniformity, consider a cut representation $(L, R)$ of a number $y$. We have $\Xi_{\mathbf{Smp}_\Pi} L <_\Pi \Xi_{\mathbf{Smp}_\Pi} R$ so the number $\{\Pi[\Xi_{\mathbf{Smp}_\Pi} L] \mid \Pi[\Xi_{\mathbf{Smp}_\Pi} R]\}_{\mathbf{S}}$ is well defined. Since $(L, R)$ is cofinal with respect to $(y_L, y_R)$ and $\Xi_{\mathbf{Smp}_\Pi}$ is strictly increasing, the number $\{\Pi[\Xi_{\mathbf{Smp}_\Pi} L] \mid \Pi[\Xi_{\mathbf{Smp}_\Pi} R]\}_{\mathbf{S}}$ lies in the cut $(\Pi[\Xi_{\mathbf{Smp}_\Pi} y_L] \mid \Pi[\Xi_{\mathbf{Smp}_\Pi} y_R])_{\mathbf{S}}$, so $\Xi_{\mathbf{Smp}_\Pi} y \sqsubseteq \{\Pi[\Xi_{\mathbf{Smp}_\Pi} L] \mid \Pi[\Xi_{\mathbf{Smp}_\Pi} R]\}_{\mathbf{S}}$. Conversely, we have $L < y < R$, so $\Xi_{\mathbf{Smp}_\Pi} L < \Xi_{\mathbf{Smp}_\Pi} y < \Xi_{\mathbf{Smp}_\Pi} R$. Since $\Xi_{\mathbf{Smp}_\Pi} L \cup \{\Xi_{\mathbf{Smp}_\Pi} y\} \cup \Xi_{\mathbf{Smp}_\Pi} R \subseteq \mathbf{Smp}_\Pi$, we have $\Pi[\Xi_{\mathbf{Smp}_\Pi} L] < \Xi_{\mathbf{Smp}_\Pi} y < \Pi[\Xi_{\mathbf{Smp}_\Pi} R]$, whence $\{\Pi[\Xi_{\mathbf{Smp}_\Pi} L] \mid \Pi[\Xi_{\mathbf{Smp}_\Pi} R]\}_{\mathbf{S}} \sqsubseteq \Xi_{\mathbf{Smp}_\Pi} y$. We conclude that $\Xi_{\mathbf{Smp}_\Pi} y = \{\Pi[\Xi_{\mathbf{Smp}_\Pi} L] \mid \Pi[\Xi_{\mathbf{Smp}_\Pi} R]\}_{\mathbf{S}}$. □

**Corollary 6.10.** *If $\Pi$ is thin and $\mathbf{S}$ is a final segment of $\mathbf{No}$, then $\Xi_{\mathbf{Smp}_\Pi}$ preserves ordinals.*

**Proof.** If $\mu$ is an ordinal, then $(\Pi[\Xi_{\mathbf{Smp}_\Pi} \mu_L] \mid \emptyset)_{\mathbf{S}}$ is a non-empty final segment of $\mathbf{S}$ and thus of $\mathbf{No}$, so by Lemma 4.17, its simplest element $\Xi_{\mathbf{Smp}_\Pi} \mu$ is an ordinal. □

For convex partitions $\Pi, \Pi'$ of $\mathbf{S}$, we write $\Pi \preceq \Pi'$ if we have $\Pi[x] \subseteq \Pi'[x]$ for every $x \in \mathbf{S}$, and say that $\Pi$ is *finer* than $\Pi'$. If $\Pi \preceq \Pi'$, then $\mathbf{Smp}_{\Pi'} \subseteq \mathbf{Smp}_\Pi$.

Recall that a *directed set* is a partial order $(J, \leqslant)$ such that for all $j, j' \in J$, there exists a $j'' \in J$ with $j, j' \leqslant j''$.

**Proposition 6.11.** *Let $\mathbf{S}$ be a surreal substructure. Let $(J, <)$ be a non-empty directed set. If $(\Pi_j)_{j \in J}$ is a $\preceq$-increasing family of thin convex partitions of $\mathbf{S}$, then the intersection $\bigcap_{j \in J} \mathbf{Smp}_{\Pi_j}$ is a surreal substructure with defining thin partition $\Pi_J$ given by*

$$\forall x \in \mathbf{S}, \quad \Pi_J[x] = \bigcup_{j \in J} \Pi_j[x].$$

**Proof.** Given $x \in \mathbf{S}$, the class $\Pi_J[x] := \bigcup_{j \in J} \Pi_j[x]$ is a non-empty convex subclass of $\mathbf{S}$ and $\bigcup_{x \in \mathbf{S}} \Pi_J[x] = \mathbf{S}$. Let $x, y \in \mathbf{S}$ be such that $\Pi_J[x] \cap \Pi_J[y] \neq \emptyset$ and let $i \in J$. Since $J$ is directed, there exists a $j \geqslant i$ in $J$ such that $\Pi_j[x] \cap \Pi_j[y] \neq \emptyset$, whence $\Pi_j[x] = \Pi_j[y]$. In particular, $\Pi_i[x] \subseteq \Pi_J[y]$ and $\Pi_i[y] \subseteq \Pi_J[x]$. Since this is true for any $i \in J$, it follows that $\Pi_J[x] = \Pi_J[y]$, so $\Pi_J$ defines a convex partition of $\mathbf{S}$.

For $x \in \mathbf{S}$, we have $\Pi_J[x_L^{\mathbf{S}}] < x < \Pi_J[x_R^{\mathbf{S}}]$ if and only if $\Pi_j[x_L^{\mathbf{S}}] < x < \Pi_j[x_R^{\mathbf{S}}]$ holds for all $j \in J$, so Lemma 6.5 implies $\bigcap_{j \in J} \mathbf{Smp}_{\Pi_j} = \mathbf{Smp}_{\Pi_J}$. Now for $x \in \mathbf{S}$, the set $\bigcup_{j \in J} \Pi_j x$ is cofinal and coinitial in $\Pi_J[x]$, so $\Pi_J$ is thin. Theorem 6.7 therefore implies that the class $\bigcap_{j \in J} \mathbf{Smp}_{\Pi_j}$ is a surreal substructure. □

**Proposition 6.12.** *Assume that $\mathbf{S}$ is a final segment of $\mathbf{No}$ and that $\Pi \preceq \Pi'$ are thin convex partitions of $\mathbf{S}$. Then for $\lambda \in \mathbf{On}$, we have $\Xi_{\mathbf{Smp}_\Pi} \lambda \leqslant \Xi_{\mathbf{Smp}_{\Pi'}} \lambda$, and in particular $\lambda \leqslant \Xi_{\mathbf{Smp}_\Pi} \lambda$.*

**Proof.** We prove the first inequality by induction on $\lambda \in \mathbf{On}$. Assuming that the inequality holds strictly below $\lambda$, we have

$$\begin{aligned}\Xi_{\mathbf{Smp}_\Pi} \lambda &= \{\Pi[\Xi_{\mathbf{Smp}_\Pi} \lambda_L] \mid \emptyset\}_{\mathbf{S}} \\ \Xi_{\mathbf{Smp}_{\Pi'}} \lambda &= \{\Pi'[\Xi_{\mathbf{Smp}_{\Pi'}} \lambda_L] \mid \emptyset\}_{\mathbf{S}}.\end{aligned}$$

For $\gamma \in \lambda_L$, we have $\Xi_{\mathbf{Smp}_\Pi} \gamma \leqslant \Xi_{\mathbf{Smp}_{\Pi'}} \gamma < \Xi_{\mathbf{Smp}_{\Pi'}} \lambda$ where $\Xi_{\mathbf{Smp}_{\Pi'}} \lambda \in \mathbf{Smp}_{\Pi'} \subseteq \mathbf{Smp}_\Pi$, so $\Pi[\Xi_{\mathbf{Smp}_\Pi} \gamma] < \Xi_{\mathbf{Smp}_{\Pi'}} \lambda$, whence in particular $\Pi[\Xi_{\mathbf{Smp}_\Pi} \lambda_L] < \Xi_{\mathbf{Smp}_{\Pi'}} \lambda$. By Proposition 4.11(a), we have $\Xi_{\mathbf{Smp}_\Pi} \lambda \leqslant \Xi_{\mathbf{Smp}_{\Pi'}} \lambda$, whence the result by induction.

The second inequality is a consequence of the first one in the case when $\Pi$ is the discrete partition of $\mathbf{S}$, which is $\angle$-minimal and for which $\Xi_{\mathbf{Smp}_\Pi} = \Xi_\mathbf{S}$. Since $\mathbf{S}$ is a final segment of $\mathbf{No}$, Proposition 4.17 gives $\Xi_\mathbf{S} 0 = \min(\mathbf{S} \cap \mathbf{On}) \geqslant 0$. Moreover, for all $\lambda \in \mathbf{On}$ with $\lambda > 0$, we have $\Xi_\mathbf{S} \lambda = \{\Xi_\mathbf{S} \lambda_L | \emptyset\}_\mathbf{S} = \{\Xi_\mathbf{S} \lambda_L | \emptyset\}$, which yields $\Xi_\mathbf{S} \lambda \geqslant \lambda$ by induction. □

## 6.2 Sharp convex partitions

We have encountered two different types of projections for surreal substructures. Given an $\mathbf{S}$-closed rooted subclass $\mathbf{X}$ of a surreal substructure $\mathbf{S}$, the topological projection sends every element $x \in \mathbf{S}^{\exists \mathbf{X}^\bullet}$ to the $\sqsubseteq$-maximal initial segment $\mu_\mathbf{X}^\mathbf{S}(x)$ of $x$ lying in $\mathbf{X}$. Given a convex partition $\Pi$ of the surreal substructure $\mathbf{S}$, the $\Pi$-simple projection sends $x \in \mathbf{S}$ to the unique $\Pi$-simple element $\pi_\Pi(x)$ lying in $\Pi[x]$. It is natural to ask whether both types of projections relate to each other.

Given a surreal substructure $\mathbf{S}$ and an $\mathbf{S}$-closed rooted subclass $\mathbf{X}$ with $\mathbf{X}^\bullet = \mathbf{S}^\bullet$, the topological projection $\mu := \mu_\mathbf{X}^\mathbf{S}$ is defined everywhere on $\mathbf{S}$. For each $x \in \mathbf{S}$, we define $\mathbf{M}_\mathbf{X}[x] := \mu^{-1}(\{\mu(x)\})$. It is easy to see that $\mathbf{M}_\mathbf{X}$ defines a partition of $\mathbf{S}$ into non-empty rooted $\sqsubseteq$-convex subclasses, and that $\mathbf{X}$ is the class of roots $\mathbf{M}_\mathbf{X}[x]^\bullet$ where $x$ ranges in $\mathbf{S}$. The members of $\mathbf{M}_\mathbf{X}$ are not necessarily $\leqslant$-convex in $\mathbf{S}$. For instance, one can prove that the structure $\mathbf{S} = \mathbf{Mo}^> + \mathbf{No}^<$ is a $\mathbf{No}^{>,>}$-closed surreal substructure, with $\mathbf{No}^{>,>} = \mathbf{Hull}(\mathbf{S})$, for which $\mathbf{M}_\mathbf{S}[\omega]$ contains $\omega$ and $\omega + 1$ but not $\omega + \omega^{-1}$.

Conversely, given a convex partition $\Pi$ of $\mathbf{S}$, the class $\mathbf{Smp}_\mathbf{S}$ may not be $\mathbf{S}$-closed, and when it is, it may be that $\mu_{\mathbf{Smp}_\mathbf{S}}^\mathbf{S}$ and $\pi_\Pi$ disagree. In some interesting cases, the projections $\mu_{\mathbf{Smp}_\Pi}^\mathbf{S}$ and $\pi_\Pi$ do coincide, and $(\mathbf{Smp}_\Pi, \sqsubseteq, \leqslant)$ has additional properties, as we shall see now.

**Definition 6.13.** *Let $\mathbf{S}$ be a surreal substructure. We say that a convex partition $\Pi$ of $\mathbf{S}$ is* **sharp**, *if the canonical representation in $\mathbf{S}$ of every $\Pi$-simple element $x$ is cofinal with respect to* $(\Pi[x_L \cap \mathbf{Smp}_\Pi], \Pi[x_R \cap \mathbf{Smp}_\Pi])$.

Assume that $\Pi$ is thin and sharp. Then each element $x \in \mathbf{Smp}_\Pi$ admits the cut representation $(\Pi[x_L^{\mathbf{Smp}_\Pi}], \Pi[x_R^{\mathbf{Smp}_\Pi}])$ in $\mathbf{S}$. By Proposition 4.11(b), this cut respresentation is mutually cofinal with $(x_L^\mathbf{S}, x_R^\mathbf{S})$. In view of Remark 4.21, we thus see that the sharpness is equivalent to the fact that the cut $(\Pi[x_L^{\mathbf{Smp}_\Pi}] | \Pi[x_R^{\mathbf{Smp}_\Pi}])_\mathbf{S}$ coincides with the $\sqsubseteq$-final substructure $\mathbf{S}^{\exists x}$ of $\mathbf{S}$ for every $x \in \mathbf{Smp}_\Pi$. This corresponds to the notion of simple representation of [11, Definition 2.2]. We say that $\mathbf{Smp}_\Pi$ is *sharp* in $\mathbf{S}$ if its defining partition is sharp.

The main interest of sharpness lies in the following equivalences:

**Theorem 6.14.** *Let $\Pi$ be a convex partition of the surreal substructure $\mathbf{S}$ such that $\mathbf{Smp}_\Pi$ is a surreal substructure. The following statements are equivalent:*

*a)* $\Pi$ *is sharp.*

*b)* $\mathbf{Smp}_\Pi$ *is $\mathbf{S}$-closed and* $\pi_\Pi = \mu_{\mathbf{Smp}_\Pi}^\mathbf{S}$.

*c)* $\pi_\Pi$ *is $\sqsubseteq$-increasing.*

*d*) $\mathbf{Smp}_\Pi$ is **S**-closed and $\mu^{\mathbf{S}}_{\mathbf{Smp}_\Pi}$ is $\leqslant$-increasing.

**Proof.** Assume that $\Pi$ is sharp. Let us prove *b*), *c*) and *d*). Note that $\mathbf{S}^\bullet$ is $\Pi$-simple, whence $\mathbf{S}^\bullet = (\mathbf{Smp}_\Pi)^\bullet$. We know that $\mu^{\mathbf{S}}_{\mathbf{Smp}_\Pi}$ when it exists is $\sqsubseteq$-increasing, and that $\pi_\Pi$ is $\leqslant$-increasing, so we need only prove that $\mathbf{Smp}_\Pi$ is **S**-closed and $\mu^{\mathbf{S}}_{\mathbf{Smp}_\Pi} = \pi_\Pi$.

Let $a, b \in \mathbf{Smp}_\Pi$ be such that $a \sqsubset b$. We claim that $b$ is simpler than no element of $\Pi[a]$. By symmetry, we may assume without loss of generality that $a < b$. Since $a \in b_L$ and $\Pi$ is sharp, the set $b_L^{\mathbf{S}}$ is cofinal with respect to $\Pi[a]$. Assume for contradiction that we have $b \sqsubseteq x$ for some $x \in \Pi[a]$. Let $y \in \Pi[a]$ be such that $x < y$ and $y \sqsubseteq b$. Then $y \sqsubseteq x$. By Lemma 6.3, we also have $b > \Pi[a]$, whence $y < b$. It follows that $x[\ell(y)] = b[\ell(y)] = 1$, whence $y < x$: a contradiction.

Since $a = \Pi[a]^\bullet$, our claim implies that $a$ is the maximal initial segment of any element of $\Pi[a] = \pi_\Pi^{-1}(\{a\})$ lying in $\mathbf{Smp}_\Pi$, i.e. that $\mu^{\mathbf{S}}_{\mathbf{Smp}_\Pi}$ is defined on $\Pi[a]$ and coincides with $\pi_\Pi$ on this class. Since the classes $\Pi[a]$ cover **S**, we see that $\mu^{\mathbf{S}}_{\mathbf{Smp}_\Pi}$ is defined on **S**, and $\pi_\Pi = \mu^{\mathbf{S}}_{\mathbf{Smp}_\Pi}$. By Proposition 5.8, the structure $\mathbf{Smp}_\Pi$ is **S**-closed.

We next prove that *a*) is a consequence of *b*). Assume for contradiction that $\mathbf{Smp}_\Pi$ is **S**-closed with $\pi_\Pi = \mu_{\mathbf{Smp}_\Pi}$ and that $\Pi$ is not sharp. We treat the case when there are $a, b \in \mathbf{Smp}_\Pi$ such that $a \in b_L$ but $b_L^{\mathbf{S}}$ has a strict upper bound $a''$ in $\Pi[a]$. Then $b_L^{\mathbf{S}} < a'' < b_R^{\mathbf{S}}$, so $b \sqsubseteq a''$, and $b \sqsubseteq \mu^{\mathbf{S}}_{\mathbf{Smp}_\Pi}(a'')$. In particular, $\pi_\Pi(a'') = a \sqsubset \mu^{\mathbf{S}}_{\mathbf{Smp}_\Pi}(a'')$, whence $\pi_\Pi \neq \mu^{\mathbf{S}}_{\mathbf{Smp}_\Pi}$: a contradiction. The other case is similar.

Assume next that $\pi_\Pi$ is $\sqsubseteq$-increasing. For $x \in \mathbf{S}$ and $a \in \mathbf{Smp}_\Pi$ such that $a \sqsubseteq x$, we have $a = \pi_\Pi(a) \sqsubseteq \pi_\Pi(x)$, so $\pi_\Pi(x)$ is the $\sqsubseteq$-maximal $\Pi$-simple initial segment of $x$. This means that $\mathbf{Smp}_\Pi$ is **S**-closed with topological projection $\pi_\Pi$. So *c*) implies *b*).

Assume $\mathbf{Smp}_\Pi$ is **S**-closed and $\mu^{\mathbf{S}}_{\mathbf{Smp}_\Pi}$ is $\leqslant$-increasing. It follows that each fiber $(\mu^{\mathbf{S}}_{\mathbf{Smp}_\Pi})^{-1}(\{\mu^{\mathbf{S}}_{\mathbf{Smp}_\Pi}(x)\})$ of $\mu^{\mathbf{S}}_{\mathbf{Smp}_\Pi}$ where $x \in \mathbf{S}$ is convex for $\leqslant$. As we have seen in the introduction of this section, we can construe $\mathbf{Smp}_\Pi$ as $\mathbf{Smp}_\mathbf{M}$ where for $x \in \mathbf{S}$, we have $\mathbf{M}[x] = \mu^{\mathbf{S}}_{\mathbf{Smp}_\Pi}{}^{-1}(\{\mu^{\mathbf{S}}_{\mathbf{Smp}_\Pi}(x)\})$. By Proposition 6.6, we have $\mu^{\mathbf{S}}_{\mathbf{Smp}_\Pi} = \pi_\mathbf{M} = \pi_\Pi$, so *d*) implies *b*). This concludes the proof. □

**Example 6.15.** Convex partitions of a surreal substructure may or may not be sharp:

- Let $\Pi_>$ denote the partition of **No** where for $x \in \mathbf{No}$, we have

$$\Pi_>[x] = \mathbf{Hull}(x + \mathbb{Z}).$$

This is actually the defining partition of the class $\mathbf{No}_> = \omega \dot{\times} \mathbf{No} = (2 \dot{\times} \mathbf{No})^{\prec \omega}$ of purely infinite surreal numbers, which is sharp, since for $x \in \mathbf{No}_>$, we have $x_L = x_L^{\mathbf{No}_>} + \mathbb{N}$ and $x_R = x_R^{\mathbf{No}_>} - \mathbb{N}$.

- Let $\Pi_1$ denote the partition of **No** where for $x \in \mathbf{No}$, we have

$$\Pi_1[x] = \mathbf{Hull}\left(x + \mathbb{Z}\, \dot{\omega}^{\omega^{-1}}\right).$$

This is a thin convex partition of **No** whose class of $\Pi_2$-simple elements contains $\dot{\omega}^{2^{-\mathbb{N}}}$. However, the number $\dot{\omega}^{\omega^{-1}} = \sup_{\sqsubseteq} \dot{\omega}^{2^{-\mathbb{N}}}$ is not $\Pi_1$-simple since it lies in $\Pi_1[0]$. Thus $\mathbf{Smp}_{\Pi_1}$ is not **No**-closed; *a fortiori* $\Pi_1$ is not sharp.

- Let **C** denote the class $\mathbf{Hull}(^1/_2 \dot{\times} \mathbf{No})$. This is a surreal substructure by Proposition 4.18. Let $\Pi_2$ denote the convex partition of **C** where for $a \in {}^1/_2 \dot{\times} \mathbf{No}$, we have

$$\Pi_2[a] = \mathbf{No}^{\sqsupseteq a \dotplus (-^1/_4)} \sqcup \{a\} \sqcup \mathbf{No}^{\sqsupseteq a \dotplus ^1/_4} \subseteq \mathbf{C}.$$

One can check that each $\mathbf{\Pi}_2[a]$ is a convex subclass of $\mathbf{C}$ and that for $x \in \mathbf{C}$, we have $\mathbf{\Pi}_2[\mu(x)] = \mu^{-1}(\{\mu(x)\})$, where $\mu$ is the topological projection $\mathbf{C} \longrightarrow {}^1\!/_2 \dot{\times} \mathbf{No}$. By Theorem 6.14, $\mathbf{\Pi}_2$ is sharp, but not thin.

We end this subsection with two further properties of sharpness.

**Proposition 6.16.** *Let $(J, <)$ be a non-empty directed set. Let $(\mathbf{\Pi}_j)_{j \in J}$ be a $\trianglelefteq$-increasing family of thin convex partitions of a surreal substructure $\mathbf{S}$. If every $\mathbf{\Pi}_j$ with $j \in J$ is sharp, then the defining thin partition $\mathbf{\Pi}_J$ of $\bigcap_{j \in J} \mathbf{Smp}_{\mathbf{\Pi}_j}$ (defined in Proposition 6.11) is sharp.*

**Proof.** We know by Proposition 6.11 that $\mathbf{\Pi}_J$ is a thin convex partition of $\mathbf{S}$ with $\mathbf{Smp}_{\mathbf{\Pi}_J} = \bigcap_{j \in J} \mathbf{Smp}_{\mathbf{\Pi}_j}$. Let $x \in \mathbf{Smp}_{\mathbf{\Pi}_J}$. For $l \in x_L^{\mathbf{Smp}_{\mathbf{\Pi}_J}}$ and $a \in \mathbf{\Pi}_J[l]$, there is $j \in J$ such that $a \in \mathbf{\Pi}_j[l]$ where $x \in \mathbf{Smp}_{\mathbf{\Pi}_j}$ and $l \in x_L^{\mathbf{Smp}_{\mathbf{\Pi}_j}}$. Since $\mathbf{\Pi}_j$ is sharp, there exists an $x' \in x_L^{\mathbf{S}}$ with $a \leqslant x'$, so $x_L^{\mathbf{S}}$ is cofinal with respect to $\mathbf{\Pi}_J[x_L^{\mathbf{Smp}_{\mathbf{\Pi}}}]$. Likewise $x_R^{\mathbf{S}}$ is coinitial with respect to $\mathbf{\Pi}_J[x_R^{\mathbf{Smp}_{\mathbf{\Pi}}}]$, so $\mathbf{\Pi}_J$ is sharp. □

**Proposition 6.17.** *Let $\mathbf{F}$ be a surreal substructure of $\mathbf{No}$ that is also a final segment. Given a thin and sharp convex partition $\mathbf{\Pi}$ of $\mathbf{F}$, we have $\Xi_{\mathbf{Smp}_{\mathbf{\Pi}}}(\mathbf{On}) = \mathbf{Smp}_{\mathbf{\Pi}} \cap \mathbf{On}$.*

**Proof.** We already know from Corollary 6.10 that $\Xi_{\mathbf{Smp}_{\mathbf{\Pi}}}(\mathbf{On}) \subseteq \mathbf{On}$. Let $a \in \mathbf{No}$ be such that $\Xi_{\mathbf{Smp}_{\mathbf{\Pi}}} a$ is an ordinal. The set $(\Xi_{\mathbf{Smp}_{\mathbf{\Pi}}} a)_R^{\mathbf{F}}$ is both empty and coinitial with respect to $\mathbf{\Pi}[\Xi_{\mathbf{Smp}_{\mathbf{\Pi}}} a_R]$, which implies that $a_R = \emptyset$ and thus that $a$ is an ordinal. □

## 6.3 Group actions

In this subsection, we study one particularly important way in which convex partitions of surreal substructures arise, namely as convex hulls of orbits under a group action.

Let $\mathbf{S}$ be a fixed surreal substructure. We define $\mathcal{F}_{\mathbf{S}}$ to be the (class-sized) group of strictly increasing bijections $g: \mathbf{S} \longrightarrow \mathbf{S}$, with functional composition as the group law. Consider any set-sized subgroup $\mathcal{G}$ of $\mathcal{F}_{\mathbf{S}}$. Then $\mathcal{G}$ naturally acts on $\mathbf{S}$ through function application; we call $\mathcal{G}$ a *function group acting on* $\mathbf{S}$.

**Definition 6.18.** *We define the **halo** $\mathcal{G}[x]$ of an element $x \in \mathbf{S}$ under the action of $\mathcal{G}$ by*

$$\mathcal{G}[x] \;=\; \{y \in \mathbf{S} : \exists g, h \in \mathcal{G}, (gx \leqslant y \leqslant hx)\} = \mathbf{Hull}_{\mathbf{S}}(\mathcal{G}x).$$

**Proposition 6.19.** *The classes $\mathcal{G}[x]$ for $x \in \mathbf{S}$ form a thin convex partition of $\mathbf{S}$.*

**Proof.** Let $x \in \mathbf{S}$. For any $y \in \mathcal{G}[x]$, we have $\mathcal{G}[y] = \mathcal{G}[x]$. Indeed, we have $gx \leqslant y \leqslant hx$ for certain $g, h \in \mathcal{G}$. Given $z \in \mathcal{G}[y]$, we also have $g'y \leqslant z \leqslant h'y$ for certain $g', h' \in \mathcal{G}$, whence $(g'g)x \leqslant g'y \leqslant z \leqslant h'y \leqslant (h'h)x$, so that $z \in \mathcal{G}[x]$. We also have $h^{-1}y \leqslant x \leqslant g^{-1}y$, whence $x \in \mathcal{G}[y]$ and $z \in \mathcal{G}[y]$ for any $z \in \mathcal{G}[x]$. The class $\mathcal{G}[x]$ is convex by definition. For $a \in \mathbf{S}$, we know that $\mathcal{G}[a]$ contains $a$, so the $\mathcal{G}[a]$ for $a \in \mathbf{S}$ form a convex partition of $\mathbf{S}$. For $x \in \mathbf{S}$, the set $\mathcal{G}x$ is cofinal and coinitial in $\mathcal{G}[x]$, so this partition is thin. □

We write $\mathbf{\Pi}_{\mathcal{G}}$ for the partition from Proposition 6.19 and say that an element of $\mathbf{S}$ is $\mathcal{G}$-*simple* if it is $\mathbf{\Pi}_{\mathcal{G}}$-simple. We let $\mathbf{Smp}_{\mathcal{G}}$ denote the class of $\mathcal{G}$-simple elements. Proposition 6.19 implies that every property from Lemmas 6.3, 6.5 and 6.4 applies to the class of $\mathcal{G}$-simple elements. We call $\pi_{\mathcal{G}} := \pi_{\mathbf{\Pi}_{\mathcal{G}}}$ the $\mathcal{G}$-*simple projection* and write $<_{\mathcal{G}}, =_{\mathcal{G}}$, and $\leqslant_{\mathcal{G}}$ instead of $<_{\mathbf{\Pi}_{\mathcal{G}}}, =_{\mathbf{\Pi}_{\mathcal{G}}}$, and $\leqslant_{\mathbf{\Pi}_{\mathcal{G}}}$.

**Proposition 6.20.** $\mathbf{Smp}_{\mathcal{G}}$ *is a surreal substructure with the following uniform cut equation in* **No***:*

$$\forall x \in \mathbf{No}, \Xi_{\mathbf{Smp}_{\mathcal{G}}} x = \{\mathcal{G} \Xi_{\mathbf{Smp}_{\mathcal{G}}} x_L | \mathcal{G} \Xi_{\mathbf{Smp}_{\mathcal{G}}} x_R\}_{\mathbf{S}}$$

**Proof.** This is a direct consequence of Proposition 6.19, Theorem 6.7 and Proposition 6.9, where we take $\mathcal{G}$ ($\mathcal{G}[x]^\bullet$) to be the required cofinal and coinitial subset of $\mathcal{G}[x]$ for each $x \in \mathbf{S}$. □

**Remark 6.21.** If $X$ is a set of strictly increasing bijective functions $\mathbf{S} \longrightarrow \mathbf{S}$, we define $\langle X \rangle$ to be the subgroup of $\mathcal{F}_\mathbf{S}$ generated by $X$, i.e. the smallest subgroup of $\mathcal{F}_\mathbf{S}$ that contains $X$. We say that $X$ is *pointwise cofinal with respect to* $Y$ and we write $Y \leq X$ if

$$\forall x \in \mathbf{S}, \forall f \in \langle Y \rangle, \exists g \in \langle X \rangle, (fx \leqslant gx).$$

This relation is transitive and reflexive. If $Y \leq X$, then $\Pi_{\langle X \rangle} \leq \Pi_{\langle Y \rangle}$, so $\mathbf{Smp}_{\langle Y \rangle} \subseteq \mathbf{Smp}_{\langle X \rangle}$. If $X \leq Y$ and $Y \leq X$, then we say that $X$ and $Y$ are *mutually pointwise cofinal* and we write $X \lessgtr Y$. In that case, we have $\mathbf{Smp}_{\langle X \rangle} = \mathbf{Smp}_{\langle Y \rangle}$.

Let us now specialize Proposition 6.11 to group-induced convex partitions.

**Proposition 6.22.** *Let* $(J, <)$ *be a non-empty directed set. If* $(\mathcal{G}_j)_{j \in J}$ *is a* $\leq$-*increasing family of function groups acting on* $\mathbf{S}$, *then the function group* $\mathcal{G}_J = \langle \mathcal{G}_j : j \in J \rangle$ *generated by* $(\mathcal{G}_j)_{j \in J}$ *satisfies*

$$\mathbf{Smp}_{\mathcal{G}_J} = \bigcap_{j \in J} \mathbf{Smp}_{\mathcal{G}_j}.$$

**Proof.** If $x \in \mathbf{S}$ is $\mathcal{G}_J$-simple, then for $j \in J$, we have $\mathcal{G}_j x_L^\mathbf{S} \subseteq \mathcal{G}_J x_L^\mathbf{S} < x < \mathcal{G}_J x_R^\mathbf{S} \supseteq \mathcal{G}_j x_R^\mathbf{S}$ so $x$ is $\mathcal{G}_j$-simple. Conversely, assume $x \in \mathbf{S}$ is $\mathcal{G}_j$-simple for all $j \in J$. Then let $g = g_{j_1} \cdots g_{j_k} \in \mathcal{G}_J$ where for $1 \leqslant k \leqslant n$, we have $g_{j_k} \in \mathcal{G}_{j_k}$. Since $(J, <)$ is directed and $(\mathcal{G}_j)_{j \in J}$ is $\leq$-increasing, there exists an index $j \in J$ with $j_1, \ldots, j_n \leqslant j$ and an element $g_j \in \mathcal{G}_j$ such that for all $u \in \mathbf{S}$ we have $g_j^{-1} u \leqslant g_{j_i} u \leqslant g_j u$ for all $i \in \{1, \ldots, n\}$, and thus $g_j^{-n} u \leqslant gu \leqslant g_j^n u$. Since $x$ is $\mathcal{G}_j$-simple, we have $g_j^n x_L^\mathbf{S} < x < g_j^{-n} x_R^\mathbf{S}$. This yields $g x_L^\mathbf{S} < x < g x_R^\mathbf{S}$, so $x$ is $\mathcal{G}_J$-simple. This proves that $\bigcap_{j \in J} \mathbf{Smp}_{\mathcal{G}_j} = \mathbf{Smp}_{\mathcal{G}_J}$. □

**Proposition 6.23.** *Let* $I$ *be a non-empty set, and let* $(\mathcal{G}_i)_{i \in I}$ *be a family of function groups acting on* $\mathbf{S}$ *such that each* $\mathbf{Smp}_{\mathcal{G}_i}$ *is sharp in* $\mathbf{S}$. *Then* $\bigcap_{i \in I} \mathbf{Smp}_{\mathcal{G}_i} = \mathbf{Smp}_{\mathcal{G}_I}$ *where* $\mathcal{G}_I = \langle \mathcal{G}_i : i \in I \rangle$.

**Proof.** We have $\bigcap_{i \in I} \mathbf{Smp}_{\mathcal{G}_i} \supseteq \mathbf{Smp}_{\mathcal{G}_I}$ for the same reasons as above. Let $x \in \bigcap_{i \in I} \mathbf{Smp}_{\mathcal{G}_i}$. Let us prove by induction on $n \in \mathbb{N}^>$ that for $g = g_{i_1} \cdots g_{i_n} \in \mathcal{G}_I$, we have $g x_L^\mathbf{S} < x < g x_R^\mathbf{S}$. By Lemma 6.5, this will prove that $x \in \mathbf{Smp}_{\mathcal{G}_I}$. For $n = 1$, the assertion is immediate. Assume therefore that $n \geqslant 2$ and decompose $g = g' g_{i_n}$, where $g' = g_{i_1} \cdots g_{i_{n-1}}$. For every $l \in x_L^\mathbf{S}$, we have $g_{i_n} l \in \mathcal{G}_{i_n}[l]$. Since $x$ is $\mathcal{G}_{i_n}$-simple, the sharpness of $\mathbf{Smp}_{\mathcal{G}_{i_n}}$ implies that there exists an $l' \in x_L^\mathbf{S}$ such that $g_{i_n} l \leqslant l'$. By our inductive hypothesis, we have $g' l' < x$, so $g l < x$. The inequality $x < g x_R^\mathbf{S}$ is proved similarly. □

**Remark 6.24.** The notions of thin convex partitions and function group actions are almost equivalent in the following sense. Let $\mathbf{\Pi}$ be a thin convex partition $\mathbf{\Pi}$ of $\mathbf{S}$, none of whose members has an extremum, and which satisfies the additional condition that there is a regular ordinal $\kappa$ with $\mathrm{cof}(\mathbf{\Pi}[x],<), \mathrm{cof}(\mathbf{\Pi}[x],>) < \kappa$ for all $x \in \mathbf{S}$. Then it can be shown that there is a group $\mathcal{G}$ acting without global fixed points on $\mathbf{S}$ such that $\mathbf{\Pi} = \mathbf{\Pi}_{\mathcal{G}}$. The converse also holds: for any function group $\mathcal{G}$ acting without global fixed points on $\mathbf{S}$, we have $\mathrm{cof}(\mathbf{\Pi}_{\mathcal{G}}[x],<), \mathrm{cof}(\mathbf{\Pi}_{\mathcal{G}}[x],>) < |\mathcal{G}|^+$ for all $x \in \mathbf{S}$.

# 7 Common group actions

## 7.1 Overview of known group actions

We conclude our study of surreal substructures with a closer examination of the action of various common types of function groups. We intentionally introduce these function groups without assigning specific domains; this will allow us to let them act on various surreal substructures.

**Translations**

Given $c \in \mathbf{No}$, we define the *translation by $c$* to be the map

$$T_c : x \longmapsto x + c.$$

The group $\mathcal{T} := \{T_r : r \in \mathbb{R}\}$ acts in particular on $\mathbf{No}$ and $\mathbf{No}^{>,\succ}$. More generally, if $A$ is a set-sized subgroup of $(\mathbf{No}, +)$, then $\mathcal{T}_A := \{T_a : a \in A\}$ acts on $\mathbf{No}$ and $(A|\emptyset)$.

Halos for the action of $\mathcal{T}$ on $\mathbf{No}$ are called *finite halos* $\mathcal{T}[x]$ and $\mathcal{T}$-simple elements correspond to purely infinite numbers. The class $\mathbf{No}_{\succ}$ of purely infinite numbers is sometimes denoted $\mathbb{J}$; see [14, 24].

**Homotheties**

Given $s \in \mathbf{No}^>$, we define the *homothety by the factor $s$* to be the map

$$H_s : x \longmapsto sx.$$

The group $\mathcal{H} := \{H_r : r \in \mathbb{R}^>\}$ acts in particular on $\mathbf{No}, \mathbf{No}^>$, and $\mathbf{No}^{>,\succ}$. More generally, if $M$ is a set-sized subgroup of $(\mathbf{No}^>, \times)$, then $\mathcal{H}_M := \{H_m : m \in M\}$ acts on $\mathbf{No}, \mathbf{No}^>$, and $(M|\emptyset)$.

Halos for the action of $\mathcal{H}$ on $\mathbf{No}^>$ are called *archimedean classes* $\mathcal{H}[x]$ and $\mathcal{H}$-simple elements are called *monomials*. The class of monomials $\mathbf{Mo} = \dot\omega^{\mathbf{No}}$ is parameterized by the $\omega$-map $\Xi_{\mathbf{Mo}}$ and forms a multiplicative cross section that is isomorphic to the value group of $\mathbf{No}$ as a valued field (the valuation being induced by the ordering). The relations $\prec_{\mathcal{H}}, \preccurlyeq_{\mathcal{H}}, =_{\mathcal{H}}$ correspond to the asymptotic relations $\prec, \preccurlyeq$, and $\asymp$ from [29, 1]. Given $x \in \mathbf{No}^{\neq}$, the projection $\pi_{\mathcal{H}}(x)$ coincides with the *dominant monomial* $\mathfrak{d}_x$, when considering $x$ as a generalized series in $\mathbb{R}[[\mathbf{Mo}]]_{\mathbf{On}}$.

**Powers**

Given $s \in \mathbf{No}^>$, we define the *$s$-th power map* by

$$P_s : x \longmapsto x^s = \exp(s \log x).$$

Here exp and log are the exponential and logarithm functions from section 3.1. The group $\mathcal{P} := \{P_r : r \in \mathbb{R}^>\}$ acts in particular on $\mathbf{No}^>$ and $\mathbf{No}^{>,>}$. More generally, if M is a set-sized subgroup of $(\mathbf{No}^>, \times)$, then the group $\mathcal{P}_M := \{P_m : m \in M\}$ acts on $\mathbf{No}^>$ and $(M|\emptyset)$.

Halos $\mathcal{P}[x]$ for the action of $\mathcal{P}$ on $\mathbf{No}^{>,>}$ are sometimes called *multiplicative classes* and $\mathcal{P}$-simple elements *fundamental monomials*. The class $\mathbf{Smp}_\mathcal{P} = \exp(\mathbf{Mo}^>) = \dot{\omega}^{\dot{\omega}^{\mathbf{No}}} = \mathbf{Mo}^{\prec 2}$ of fundamental monomials is parameterized by the $\omega^\omega$-map: see [30, Proposition 2.5].

**Exponentials**

Writing

$$\exp_n := \exp \circ \overset{n\times}{\cdots} \circ \exp$$
$$\log_n := \log \circ \overset{n\times}{\cdots} \circ \log$$

for all $n \in \mathbb{N}$, we define

$$\mathcal{E}^* := \langle \exp \rangle$$
$$\mathcal{E} := \langle \exp_n \circ H_r \circ \log_n : r \in \mathbb{R}^>, n \in \mathbb{N} \rangle.$$

Both $\mathcal{E}^*$ and $\mathcal{E}$ act in particular on $\mathbf{No}^{>,>}$.

Halos $\mathcal{E}[x]$ and $\mathcal{E}^*[x]$ for the actions of $\mathcal{E}$ and $\mathcal{E}^*$ on $\mathbf{No}^{>,>}$ are sometimes called *levels* and *logarithmic-exponential classes* respectively. The $\mathcal{E}$-simple elements are called *log-atomic* numbers and the class $\mathbf{La}$ of such numbers is parameterized by the $\lambda$-map: see [11, Section 5]. The class of $\mathcal{E}^*$-simple elements is denoted by $\mathbf{K}$ and parameterized by the $\kappa$-map: see [30, Section 3].

We notice that each of the above function groups is linearly ordered by

$$f \leqslant g \iff \exists x_0 \in \mathbf{No}, \forall x > x_0, f(x) \leqslant g(x).$$

With the exception of $\mathcal{E}$, all these groups are also abelian. These are both strong properties which need not be imposed for the material of Section 6.3 to apply.

## 7.2 Actions by translations

Throughout this subsection, let $A$ be a fixed set-sized subgroup of $(\mathbf{No}, +)$ and let $\Xi_A := \Xi_{\mathbf{Smp}_{\mathcal{T}_A}}$. If $A \subseteq \mathbf{No}^\prec$, then $\mathcal{T}_A \nleq \mathcal{T}$ so $\mathbf{No}_> = \mathbf{Smp}_\mathcal{T} \subseteq \mathbf{Smp}_{\mathcal{T}_A}$. If $A \nsubseteq \mathbf{No}^\prec$, then given $a \in A \setminus \mathbf{No}^\prec$, the set $\mathbb{Z}a$ is cofinal with respect to $\mathbb{R}$, so $\mathcal{T} \nleq \mathcal{T}_A$, whence $\mathbf{Smp}_{\mathcal{T}_A} \subseteq \mathbf{No}_>$.

**Proposition 7.1.** *If $\mathcal{T}_A$ acts on $\mathbf{No}$, then $\Xi_A : (\mathbf{No}, +, \leqslant, \sqsubseteq) \longrightarrow (\mathbf{Smp}_{\mathcal{T}_A}, +, \leqslant, \sqsubseteq)$ is an isomorphism.*

**Proof.** We already know that $\Xi_A$ is a $(\leqslant, \sqsubseteq)$-isomorphism so we only need to prove that it preserves sums. Let $a, b \in \mathbf{No}$ be such that $\Xi_A$ preserves sums of elements lexicographically strictly simpler than $(a, b)$. Recall that the addition is uniform in the sense that

$$a = \{L_a | R_a\}, \ b = \{L_b | R_b\} \implies a + b = \{a + L_b, L_a + b | a + R_b, R_a + b\}$$

Applying this to the cut equations given by Proposition 6.9 for $\Xi_A$, we obtain

$$\Xi_A a + \Xi_A b = \{\mathcal{T}_A \Xi_A a_L | \mathcal{T}_A \Xi_A a_R\} + \{\mathcal{T}_A \Xi_A b_L | \mathcal{T}_A \Xi_A b_R\}$$
$$= \{\Xi_A a + \mathcal{T}_A \Xi_A b_L, \mathcal{T}_A \Xi_A a_L + \Xi_A b | \Xi_A a + \mathcal{T}_A \Xi_A b_R, \mathcal{T}_A \Xi_A a_R + \Xi_A b\},$$

and by uniformity of the cut equation for $\Xi_A$, we get

$$\begin{aligned}\Xi_A(a+b) &= \{\mathcal{J}_A\,\Xi_A(a+b_L), \mathcal{J}_A\,\Xi_A(a_L+b) \mid \mathcal{J}_A\,\Xi_A(a+b_R), \mathcal{J}_A\,\Xi_A(a_R+b)\} \\ &= \{\Xi_A\,a + \mathcal{J}_A\,\Xi_A\,b_L, \mathcal{J}_A\,\Xi_A\,a_L + \Xi_A\,b \mid \Xi_A\,a + \mathcal{J}_A\,\Xi_A\,b_R, \mathcal{J}_A\,\Xi_A\,a_R + \Xi_A\,b\}.\end{aligned}$$

Thus $\Xi_A(a+b) = \Xi_A\,a + \Xi_A\,b$. By induction, this proves that $\Xi_A$ preserves sums of surreals and consequently that $\mathbf{Smp}_{\mathcal{J}_A}$ is an additive subgroup of $\mathbf{No}$. $\square$

Let us now focus on $\mathcal{J}$. By induction on $\alpha \in \mathbf{On}$, it is easy to see that $\mathbf{No}_{>}^{\prec \alpha} = \dot{\omega}^\alpha \dot{\times} \mathbf{No}$ and $\Xi_{\mathbf{No}_{>}}^{\alpha} x = \dot{\omega}^\alpha \dot{\times} x$ for all $x \in \mathbf{No}$. In particular, this gives a description of $\mathbf{Fix}_{\mathbf{No}_{>}} = \dot{\omega}^\omega \dot{\times} \mathbf{No}$ in terms of sign sequences.

Let us next describe the structures $\mathbf{No}_{>}^{\prec \alpha}$ for $\alpha \in \mathbf{On}$ in terms of Conway normal forms and of $\mathcal{G}$-simplicity for some group $\mathcal{G}$ acting on $\mathbf{No}$. By [15, Corollary 3.1], if $\alpha$ is an ordinal, then the set $\mathbf{No}(\dot{\omega}^\alpha)$ is a subgroup of $(\mathbf{No},+)$, which acts by translations on $\mathbf{No}$. If $\alpha = 1$, then the sets $\{k\,\dot{\omega}^\beta : \beta < \alpha, k \in \mathbb{Z}\}$ and $\mathbf{No}(\dot{\omega}^\alpha)$ are mutually cofinal and coinitial, and $\mathbf{No}_{>} = \mathbf{Smp}_{\mathcal{J}_{\mathbf{No}(\omega)}}$, since $\mathbf{No}(\omega) = \mathbb{R}$. We claim that this generalizes to every ordinal.

**Proposition 7.2.** *For $\alpha \in \mathbf{On}$, we have $\mathbf{No}_{>}^{\prec \alpha} = \dot{\omega}^\alpha \dot{\times} \mathbf{No} = \mathbf{Smp}_{\mathcal{J}_{\mathbf{No}(\dot{\omega}^\alpha)}}$.*

**Proof.** We proceed by induction on $\alpha \in \mathbf{On}$. The result obviously holds for $\alpha = 0$. We saw that it holds for $\alpha = 1$ in Example 5.3. Assume that $\alpha = \beta + 1$ is a successor ordinal. Then the function $\Xi_{\mathbf{No}_{>}^{\prec \beta}}$ is additive by Proposition 7.1, so $\Xi_{\mathbf{No}_{>}^{\prec \beta}}\,\mathbb{Z} = \mathbb{Z}\,\Xi_{\mathbf{No}_{>}^{\prec \beta}}\,1 = \mathbb{Z}\,\dot{\omega}^\beta$ is mutually cofinal and coinitial with $\mathbf{No}(\dot{\omega}^\alpha)$. Let $\theta$ be $\mathcal{J}_{\mathbf{No}(\dot{\omega}^\alpha)}$-simple. Then $\theta$ is $\mathcal{J}_{\mathbf{No}(\dot{\omega}^\beta)}$-simple, so the inductive hypothesis yields $\theta = \Xi_{\mathbf{No}_{>}^{\prec \beta}} x$ for a certain number $x$. Since $\theta \sqsubseteq \theta + \mathbb{Z}\,\dot{\omega}^\beta = \Xi_{\mathbf{No}_{>}^{\prec \beta}}(x + \mathbb{Z})$, we deduce that $x \sqsubseteq x + \mathbb{Z}$. Now for $z \in \mathcal{J}_{\mathbb{Z}}[x]$, there is $n \in \mathbb{N}$ with $x - n < z < x + n$. We cannot have both $x < z$ and $x > z$, so the contrapositive of Lemma 4.4 yields $x \sqsubseteq z$. Thus $x$ is $\mathcal{J}_{\mathbb{Z}}$-simple, so $\theta \in \mathbf{No}_{>}^{\prec \beta} \prec \mathbf{No}_{>} = \mathbf{No}_{>}^{\prec \alpha}$. Conversely, for $\theta \in \mathbf{No}_{>}^{\prec \alpha}$, we have $\theta = \Xi_{\mathbf{No}_{>}^{\prec \beta}} x$ for a certain $x \in \mathbf{No}_{>}$. We have $x \sqsubseteq x + \mathbb{Z}$, so $\theta \sqsubseteq \theta + \mathbb{Z}\,\dot{\omega}^\beta$. Similar arguments as above yield $\theta \sqsubseteq \mathcal{J}_{\mathbf{No}(\dot{\omega}^\alpha)}[\theta]$, whence $\theta \in \mathbf{Smp}_{\mathcal{J}_{\mathbf{No}(\dot{\omega}^\alpha)}}$. This proves that $\mathbf{Smp}_{\mathcal{J}_{\mathbf{No}(\dot{\omega}^\alpha)}} = \mathbf{No}_{>}^{\prec \alpha}$.

If $\alpha$ is a limit ordinal, then Proposition 6.11 yields

$$\begin{aligned}\mathbf{No}_{>}^{\prec \alpha} &= \bigcap_{\beta < \alpha} \mathbf{No}_{>}^{\prec \beta} \\ &= \bigcap_{\beta < \alpha} \mathbf{Smp}_{\mathcal{J}_{\mathbf{No}(\dot{\omega}^\beta)}} \\ &= \mathbf{Smp}_{\mathcal{J}_{\bigcup_{\beta < \alpha} \mathbf{No}(\dot{\omega}^\beta)}} \\ &= \mathbf{Smp}_{\mathcal{J}_{\mathbf{No}(\dot{\omega}^\alpha)}}.\end{aligned}$$

$\square$

A consequence of Propositions 7.1 and 7.2 is that $\Xi_{\mathbf{No}_{>}}^{\alpha}$ is additive for all $\alpha \in \mathbf{On}$. In fact, we even have the following:

**Proposition 7.3.** *For $\alpha \in \mathbf{On}$, the function $\Xi_{\mathbf{No}_{>}}^{\alpha} : \mathbb{R}[[\mathbf{Mo}]]_{\mathbf{On}} \longrightarrow \mathbb{R}[[\mathbf{Mo}]]_{\mathbf{On}}$ is strongly linear, with $\mathbf{No}_{>}^{\prec \alpha} \prec \mathbf{Mo} = \mathbf{Mo} \prec \mathbf{No}^{\sqsupseteq \alpha}$.*

**Proof.** Let $\alpha \in \mathbf{On}^{>}$ and $\Phi := \Xi_{\mathbf{No}_{>}}^{\alpha}$. Let us first show that $\Phi(rx) = r\,\Phi x$ for all $r \in \mathbb{R}$ and $x \in \mathbf{No}$. By Proposition 7.2, the function $\Phi$ is additive, so this holds for any dyadic number $r$. In particular we have $\Phi(0) = \omega^\alpha \dot{\times} 0 = 0$. Let $r$ be a non-dyadic real number. Let $x \in \mathbf{No}$ be such that $\Phi(ry) = r\,\Phi y$ for all $y \in x_{\sqsubset}$. It is well known that $r_{\sqsubset}$ contains only dyadic numbers. By Proposition 7.2 and (3.5), we have

$$\Phi(rx) = \{L_1, L_2 \mid R_1, R_2\}$$

where

$$\begin{aligned}
L_1 &= \Phi(r_L x + r x_L - r_L x_L) + \mathbf{No}(\dot{\omega}^\alpha), \\
L_2 &= \Phi(r_R x + r x_R - r_R x_R) + \mathbf{No}(\dot{\omega}^\alpha), \\
R_1 &= \Phi(r_L x + r x_R - r_L x_R) + \mathbf{No}(\dot{\omega}^\alpha), \text{ and} \\
R_2 &= \Phi(r_R x + r x_L - r_R x_L) + \mathbf{No}(\dot{\omega}^\alpha).
\end{aligned}$$

The cut equation (3.5) for the surreal product by $r$ is uniform [24, Theorem 3.5], so

$$r \Phi x = \{A_1, A_2 | B_1, B_2\},$$

where

$$\begin{aligned}
A_1 &= \{r' \Phi x + r(\Phi x' + \mathbf{No}(\dot{\omega}^\alpha)) - r'(\Phi x' + \mathbf{No}(\dot{\omega}^\alpha))\}, \\
A_2 &= \{r'' \Phi x + r(\Phi x'' + \mathbf{No}(\dot{\omega}^\alpha)) - r''(\Phi x'' + \mathbf{No}(\dot{\omega}^\alpha))\}, \\
B_1 &= \{r' \Phi x + r(\Phi x'' + \mathbf{No}(\dot{\omega}^\alpha)) - r'(\Phi x'' + \mathbf{No}(\dot{\omega}^\alpha))\}, \text{ and} \\
B_2 &= \{r'' \Phi x + r(\Phi x' + \mathbf{No}(\dot{\omega}^\alpha)) - r''(\Phi x' + \mathbf{No}(\dot{\omega}^\alpha))\},
\end{aligned}$$

where $r', r'', x', x''$ respectively range in $r_L, r_R, x_L, x_R$. Let us prove that $L_1$ and $A_1$ are mutually cofinal. Analog relations hold for the other sets so this will yield $r \Phi x = \Phi(rx)$. Since $\Phi$ is additive, for $r' \in r_L$ and $x' \in x_L$, we have

$$\Phi(r' x + r x' - r' x') + \mathbf{No}(\dot{\omega}^\alpha) = \Phi(r' x) + \Phi(r x') - \Phi(r' x') + \mathbf{No}(\dot{\omega}^\alpha).$$

Now $\Phi(r x') = r \Phi x'$ and $\Phi(r' x') = r' \Phi x'$ by our inductive hypothesis. Moreover, we have $\Phi(r' x) = r' \Phi x$, since $r'$ is dyadic. It follows that

$$\Phi(r' x + r x' - r' x') + \mathbf{No}(\dot{\omega}^\alpha) = r' \Phi x + r \Phi x' - r' \Phi x' + \mathbf{No}(\dot{\omega}^\alpha).$$

Since $r$ is non-zero, we have $\{r', r\} \mathbf{No}(\dot{\omega}^\alpha) = \mathbf{No}(\dot{\omega}^\alpha)$, so this set is mutually cofinal with the set $r' \Phi x + r(\Phi x' + \mathbf{No}(\dot{\omega}^\alpha)) - r'(\Phi x' + \mathbf{No}(\dot{\omega}^\alpha))$. Therefore $\Phi$ is $\mathbb{R}$-linear.

Let us next prove by induction that $\mathbf{No}_{\succ}^{\prec \alpha} \prec \mathbf{Mo} = \mathbf{Mo} \prec \mathbf{No}^{\sqsupset \alpha}$. Let $x \in \mathbf{No}$ be such that $\Phi \Xi_{\mathbf{Mo}} y = \Xi_{\mathbf{Mo}} \Xi_{\mathbf{No}^{\sqsupset \alpha}} y$ for all $y \in x_{\sqsubset}$. Let $(L, R)$ be an arbitrary cut representation in $\mathbf{No}$ such that $L$ (resp. $R$) has no maximum (resp. minimum), so that $\Phi L$ (resp $\Phi R$) has no minimum (resp. maximum). Then we note that the cut equation

$$\Phi \{L | R\} = \{\mathcal{J}_{\mathbf{No}(\dot{\omega}^\alpha)} \Phi L | \mathcal{J}_{\mathbf{No}(\dot{\omega}^\alpha)} \Phi R\}$$

simplifies as

$$\Phi \{L | R\} = \{\Phi L | \Phi R\}.$$

Considering the cut representation $(\{0\} \cup \mathbb{R}^{>} \dot{\omega}^{x_L}, \mathbb{R}^{>} \dot{\omega}^{x_R})$ of $\dot{\omega}^x$, we deduce that we have

$$\Phi \dot{\omega}^x = \{\mathbf{No}(\dot{\omega}^\alpha) + \Phi(0), \Phi(\mathbb{R}^{>}) \Xi_{\mathbf{Mo}} x_L | \Phi(\mathbb{R}^{>}) \Xi_{\mathbf{Mo}} x_R\}$$

We have seen that $\Phi$ is $\mathbb{R}$-linear, so the induction hypothesis yields

$$\begin{aligned}
\Phi \dot{\omega}^x &= \{\mathbf{No}(\dot{\omega}^\alpha), \mathbb{R}^{>} \Phi \dot{\omega}^{x_L} | \mathbb{R}^{>} \Phi \dot{\omega}^{x_R}\} \\
&= \{\mathbf{No}(\dot{\omega}^\alpha), \mathbb{R}^{>} \dot{\omega}^{\alpha + x_L} | \mathbb{R}^{>} \dot{\omega}^{\alpha + x_R}\} \\
&= \{\mathbb{R}^{>} \dot{\omega}^{\alpha_L}, \mathbb{R}^{>} \dot{\omega}^{\alpha + x_L} | \mathbb{R}^{>} \dot{\omega}^{\alpha + x_R}\} \quad \text{since } \mathbb{R}^{>} \omega^{\alpha_L} \text{ and } \mathbf{No}(\dot{\omega}^\alpha) \text{ are mutually cofinal} \\
&= \{\mathbb{R}^{>} \dot{\omega}^{(\alpha + x)_L} | \mathbb{R}^{>} \dot{\omega}^{(\alpha + x)_R}\} \\
&= \dot{\omega}^{\alpha + x}.
\end{aligned}$$

We thus have:
$$\mathbf{No}_{>}^{\prec\alpha} \prec \mathbf{Mo} = \mathbf{Mo} \prec \mathbf{No}^{\exists\alpha}.$$

In particular $\Phi$ preserves monomials.

Let $x = \sum x_\mathfrak{m}\, \mathfrak{m}$ be a number considered as a series in $\mathbb{R}[[\mathbf{Mo}]]_{\mathbf{On}}$. By our previous arguments, the number $y = \sum x_\mathfrak{m}\, \Phi\mathfrak{m}$ is well defined. For all $\mathfrak{n} \in \mathbf{Mo}$, we will write $x_{>\mathfrak{n}} = \sum_{\mathfrak{m}>\mathfrak{n}} x_\mathfrak{m}\, \mathfrak{m}$ and $y_{>\Phi\mathfrak{n}} = \sum_{\mathfrak{m}>\mathfrak{n}} x_\mathfrak{m}\, \Phi\mathfrak{m}$. Let us prove by induction on the order type $\ell_{\mathbf{Mo}}(x)$ of $(\operatorname{supp} x, >)$ that $y = \Phi x$; this will conclude the proof. The additivity and $\mathbb{R}$-linearity of $\Phi$ yield the result for $\ell_{\mathbf{Mo}}(x) < \omega$. If $\ell_{\mathbf{Mo}}(x)$ is successor and infinite, then $\operatorname{supp} x$ has a minimum $\mathfrak{m}_x$ and $x = x_{>\mathfrak{m}_x} + x_{\mathfrak{m}_x}\mathfrak{m}_x$, so

$$\begin{aligned}
\Phi x &= \Phi x_{>\mathfrak{m}_x} + \Phi(x_{\mathfrak{m}_x}\mathfrak{m}_x) \\
&= \left(\sum_{\mathfrak{m}>\mathfrak{m}_x} x_\mathfrak{m}\, \Phi\mathfrak{m}\right) + x_{\mathfrak{m}_x}\, \Phi\mathfrak{m}_x \\
&= y.
\end{aligned}$$

Assume now that $\ell_{\mathbf{Mo}}(x)$ is an infinite limit. Since $\Phi$ is strictly increasing and monomial preserving, [24, Lemma 5.3] yields

$$\begin{aligned}
x &= \{x_{>\mathfrak{n}} + (x_\mathfrak{n} - 2^{-\mathbb{N}})\,\mathfrak{n} \mid x_{>\mathfrak{n}} + (x_\mathfrak{n} + 2^{-\mathbb{N}})\,\mathfrak{n}\}, \\
y &= \{y_{>\Phi\mathfrak{n}} + (x_\mathfrak{n} - 2^{-\mathbb{N}})\,\Phi\mathfrak{n} \mid y_{>\Phi\mathfrak{n}} + (x_\mathfrak{n} + 2^{-\mathbb{N}})\,\Phi\mathfrak{n}\},
\end{aligned}$$

where $\mathfrak{n}$ ranges over $\operatorname{supp} x$. Notice that the left (resp. right) options in the above representation of $x$ have no maximum (resp. minimum), so

$$\Phi x = \{\Phi x_{>\mathfrak{n}} + \Phi((x_\mathfrak{n} - 2^{-\mathbb{N}})\,\mathfrak{n}) \mid \Phi x_{>\mathfrak{n}} + \Phi((x_\mathfrak{n} + 2^{-\mathbb{N}})\,\mathfrak{n})\}.$$

Our inductive hypothesis yields

$$\begin{aligned}
\Phi x &= \{y_{>\Phi\mathfrak{n}} + (x_\mathfrak{n} - 2^{-\mathbb{N}})\,\Phi\mathfrak{n} \mid y_{>\Phi\mathfrak{n}} + (x_\mathfrak{n} + 2^{-\mathbb{N}})\,\Phi\mathfrak{n}\} \\
&= y.
\end{aligned}$$

This concludes the proof. $\square$

**Proposition 7.4.** *For $\alpha \in \mathbf{On}$, we have*

$$\mathbf{No}_{>}^{\prec\alpha} = \mathbb{R}[[\dot{\omega}^{\mathbf{No}^{\exists\alpha}}]]_{\mathbf{On}}.$$

*In particular $\mathbf{No}_{>}^{\prec\alpha}$ is a non-unitary subring of $\mathbf{No}$, and*

$$\mathbf{Fix}_{\mathbf{No}_{>}} = \mathbb{R}[[\dot{\omega}^{\mathbf{No}^{>,\succ}}]]_{\mathbf{On}}.$$

**Proof.** The strong linearity of $\Xi_{\mathbf{No}_{>}^{\prec\alpha}}$ and the relation $\mathbf{No}_{>}^{\prec\alpha} \prec \mathbf{Mo} = \mathbf{Mo} \prec \mathbf{No}^{\exists\alpha}$ give

$$\begin{aligned}
\mathbf{No}_{>}^{\prec\alpha} &= \Xi_{\mathbf{No}_{>}^{\prec\alpha}}\, \mathbb{R}[[\mathbf{Mo}]]_{\mathbf{On}} \\
&= \mathbb{R}[[\Xi_{\mathbf{No}_{>}^{\prec\alpha}}\, \mathbf{Mo}]]_{\mathbf{On}} \\
&= \mathbb{R}[[\mathbf{Mo} \prec \mathbf{No}^{\exists\alpha}]]_{\mathbf{On}} \\
&= \mathbb{R}[[\dot{\omega}^{\mathbf{No}^{\exists\alpha}}]]_{\mathbf{On}}.
\end{aligned}$$

That this forms a (non-unitary) subring follows from the fact that $\mathbf{No}_{>}^{\exists\alpha} = (\alpha_L | \emptyset)$ is closed under addition, whence $\dot{\omega}^{\mathbf{No}^{\exists\alpha}}$ is closed under multiplication. $\square$

## 7.3 Actions by homotheties

In this subsection, M is a set-sized subgroup of $(\mathbf{No}^>, \times)$ and $\Xi_M$ the defining isomorphism of $\mathbf{Smp}_{\mathcal{H}_M}$. We will distinguish between confined and ample subgroups. We say that M is *confined* if it is a subgroup of $1+\mathbf{No}^<$ and *ample* if not. If M is ample, then given $m \in M \setminus (1+\mathbf{No}^<)$, the maximum $a = \max(m, m^{-1}) > 1$ satisfies $a - 1 \succcurlyeq 1$, which implies that $a^{\mathbb{N}}$ is cofinal with respect to $\mathbb{R}$. Thus $\mathcal{H} \not\leq \mathcal{H}_M$ on $\mathbf{No}^>$, so $\mathbf{Smp}_{\mathcal{H}_M} \subseteq \mathbf{Mo}$. If M is confined, then $\mathcal{H}_M \not\leq \mathcal{H}$, so $\mathbf{Mo} \subseteq \mathbf{Smp}_{\mathcal{H}_M}$. For $\alpha \in \mathbf{On}$, natural examples of ample multiplicative subgroups include $\mathbf{No}(\varepsilon_\alpha)^>$ for $\alpha \in \mathbf{On}$, whereas natural examples of confined multiplicative subgroups include $1+\mathbf{No}(\varepsilon_\alpha)^<$.

**Remark 7.5.** If M is confined, then $1/2 \notin 2 \dot{\times} \mathbf{No}$ is $\mathcal{H}_M$-simple but $(1/2)_R^{\geqslant} = \{1\}$ is not coinitial with respect to $\mathcal{H}_M[1]$ which contains elements strictly below 1. So $\mathbf{Smp}_{\mathcal{H}_M}$ is not sharp in $\mathbf{No}^>$. The standard monomial group $\mathbf{Mo}$ is sharp both in $\mathbf{No}^>$ and in $\mathbf{No}^{>,>}$ by [11, Corollary 4.17], but this observation does not generalize to arbitrary ample multiplicative subgroups M of $\mathbf{No}^{>,>}$. For instance, if $M = \mathbb{Q}^> \dot{\omega}^{\mathbb{Z}\omega^{1/3}}$, then $\dot{\omega}^{\omega^{1/2}}$ is $\mathcal{H}_M$-simple and $1 \in (\dot{\omega}^{\omega^{1/2}})_L^{\mathbf{Smp}_{\mathcal{H}_M}}$ but the number $\dot{\omega}^{\omega^{1/3}} \in \mathcal{H}_M[1]$ lies strictly above $(\dot{\omega}^{\omega^{1/2}})_L^{\mathbf{No}^>}$.

**Proposition 7.6.** *Assume that* M *is ample and let* $\mathcal{H}_M$ *act on* $\mathbf{No}^>$. *Then the parameterization* $\Xi_M$ *of* $\mathbf{Smp}_{\mathcal{H}_M}$ *is an isomorphism* $(\mathbf{No}^>, +, \leqslant, \sqsubseteq) \longrightarrow (\mathbf{Smp}_{\mathcal{H}_M}, \times, \leqslant, \sqsubseteq)$.

**Proof.** We only need to prove that $\Xi_M$ is a morphism $(\mathbf{No}, +) \longrightarrow (\mathbf{Smp}_{\mathcal{H}_M}, \times)$. Consider monomials $\mathfrak{m}, \mathfrak{n} \in \mathbf{Mo}$ with cut representations $(L_\mathfrak{m}, R_\mathfrak{m})$ and $(L_\mathfrak{n}, R_\mathfrak{n})$ such that $\mathbb{R} L_\mathfrak{m} \subseteq \mathbf{Hull}(L_\mathfrak{m})$, $\mathbb{R} R_m \subseteq \mathbf{Hull}(R_\mathfrak{m})$, and likewise for $\mathfrak{n}$. Then [11, Proposition 4.19] yields

$$\mathfrak{m}\mathfrak{n} = \{L_\mathfrak{m} \mathfrak{n} + \mathfrak{m} L_\mathfrak{n} \mid R_\mathfrak{m} \mathfrak{n}, \mathfrak{m} R_\mathfrak{n}\}.$$

Given $x, y \in \mathbf{No}$, this applies in particular to the cut representation $(\{0\} \cup M \Xi_M x_L, M \Xi_M x_R)$ of $\Xi_M x$ (and likewise for $\Xi_M y$) since M is ample. We thus have

$$\begin{aligned}\Xi_M x \Xi_M y &= \{0, M\Xi_M x_L \Xi_M y + \Xi_M x M \Xi_M y_L \mid M \Xi_M x_R \Xi_M y, \Xi_M x M \Xi_M y_R\} \\ &= \{0, M\Xi_M x_L \Xi_M y + M \Xi_M x \Xi_M y_L \mid M\Xi_M x_R \Xi_M y, M \Xi_M x \Xi_M y_R\}\end{aligned}$$

Note that $\Xi_M x \Xi_M y > M \Xi_M x_L \Xi_M y, M \Xi_M x \Xi_M y_L$. Assume $x_L \neq \emptyset$ and $y_L \neq \emptyset$. Since M is ample, there exists a $c \in M$ such that $c \geqslant 2$. For $l_x \in x_L$, $l_y \in y_L$, and $m, m' \in M$, we have $m \Xi_M l_x \Xi_M y + m' \Xi_M x \Xi_M l_y \leqslant \max(cm, cm') \max(\Xi_M l_x \Xi_M y, \Xi_M x \Xi_M l_y)$. This proves the following relation (which also holds when $x_L = \emptyset$ or $y_L = \emptyset$, by what precedes):

$$\Xi_M x \Xi_M y = \{0, M\Xi_M x_L \Xi_M y, M \Xi_M x \Xi_M y_L \mid M\Xi_M x_R \Xi_M y, M\Xi_M x \Xi_M y_R\}.$$

Now let $x, y$ be numbers such that for any $a, b \in \mathbf{No}$ with $a \sqsubseteq x$, $b \sqsubseteq y$, and $(a, b) \neq (x, y)$, we have $\Xi_M(a+b) = \Xi_M a \Xi_M b$. Then

$$\begin{aligned}\Xi_M(x+y) &= \{0, M\Xi_M(x_L + y), M\Xi_M(x+y_L) \mid M\Xi_M(x_R + y), M\Xi_M(x+y_R)\} \\ &= \{0, M\Xi_M x_L \Xi_M y, M\Xi_M x \Xi_M y_L \mid M\Xi_M x_R \Xi_M y, M\Xi_M x \Xi_M y_R\} \\ &= \Xi_M x \Xi_M y.\end{aligned}$$

We conclude by induction. □

The above proof fails if M is confined, since then $\Xi_M 1 = 2$ and $\Xi_M 2 = 3 \neq 2 \times 2$.

**Corollary 7.7.** *If* M *is ample, then the* $\mathcal{H}_M$-*simple projection* $\pi_{\mathcal{H}_M}$ *is a surjective morphism* $(\mathbf{No}^>, \times, \leqslant) \longrightarrow (\mathbf{Smp}_{\mathcal{H}_M}, \times, \leqslant)$.

**Proof.** We only need to prove that $\pi_{\mathcal{H}_\mathrm{M}}$ preserves products. Given $x, y \in \mathbf{No}^>$, the relation $\mathrm{M} xy = x \mathrm{M} y$ implies $xy =_{\mathcal{H}_\mathrm{M}} \pi_{\mathcal{H}_\mathrm{M}}(x)\, \pi_{\mathcal{H}_\mathrm{M}}(y)$. Proposition 7.6 implies that $\pi_{\mathcal{H}_\mathrm{M}}(x)\, \pi_{\mathcal{H}_\mathrm{M}}(y)$ is $\mathcal{H}_\mathrm{M}$-simple, whence $\pi_{\mathcal{H}_\mathrm{M}}(xy) = \pi_{\mathcal{H}_\mathrm{M}}(x)\, \pi_{\mathcal{H}_\mathrm{M}}(y)$. □

**Proposition 7.8.** *Let $\mathbf{R} \supseteq \mathbb{R}$ be a proper convex subring of $(\mathbf{No}, +, \times, \leqslant)$ with a cofinal subset. Let $\mathbf{M} = (\mathbf{R}^\times)^>$ and write $\mathfrak{M}_\mathbf{R}$ for the convex subgroup $\mathbf{M} \cap \mathbf{Mo}$ of $\mathbf{Mo}$. Define $\mathfrak{N}_\mathbf{R}$ to be the group $\mathbf{Smp}_{\mathcal{H}_\mathrm{M}}$. Then there is a canonical strongly linear isomorphism of ordered valued fields*

$$\mathbf{No} \longrightarrow \mathbb{R}[[\mathfrak{M}_\mathbf{R}]]_{\mathbf{On}}[[\mathfrak{N}_\mathbf{R}]]_{\mathbf{On}}.$$

**Proof.** By [1, Page 713], we only need to prove that $\mathfrak{M}_\mathbf{R}$ is a convex subgroup of $\mathbf{Mo}$ with $\mathfrak{M}_\mathbf{R} \cap \mathfrak{N}_\mathbf{R} = \{1\}$ and $\mathbf{Mo} = \mathfrak{M}_\mathbf{R} \mathfrak{N}_\mathbf{R}$. Since $\mathbf{R}$ has a cofinal subset, the group $\mathbf{M}$ has an ample cofinal and coinitial subgroup M and we may apply the two previous results to $\mathbf{M}$.

Intersections and convex hulls of subgroups are again subgroups, so $\mathfrak{M}_\mathbf{R}$ is a convex subgroup of $\mathbf{Mo}$. We claim that for $\mathfrak{a} \in \mathbf{Mo}$, we have $\mathfrak{a} = \mathfrak{m}\mathfrak{n}$ where $\mathfrak{n} := \pi_{\mathcal{H}_\mathrm{M}}(\mathfrak{a}) \in \mathfrak{N}_\mathbf{R}$ and $\mathfrak{m} := \frac{\mathfrak{a}}{\mathfrak{n}} \in \mathfrak{M}_\mathbf{R}$. Indeed, as a product of monomials, $\mathfrak{m}$ is a monomial. Furthermore, Corollary 7.7 yields

$$\pi_{\mathcal{H}_\mathrm{M}}(\mathfrak{m}) = \frac{\pi_{\mathcal{H}_\mathrm{M}}(\mathfrak{a})}{\pi_{\mathcal{H}_\mathrm{M}}(\mathfrak{n})} = \frac{\pi_{\mathcal{H}_\mathrm{M}}(\mathfrak{a})}{\pi_{\mathcal{H}_\mathrm{M}}(\mathfrak{a})} = 1$$

whence $\mathfrak{m} \in \mathcal{H}_\mathbf{M}[1]$. This means that there exist $r \leqslant r' \in (\mathbf{R}^\times)^>$ with $r \leqslant \mathfrak{m} \leqslant r'$. In other words, we have $\mathfrak{m} \in \mathfrak{M}_\mathbf{R}$. This concludes the proof. □

**Remark 7.9.** Assume now that $\mathbf{Hull}(\mathrm{M})$ is closed under exp. In [10], an alternative to Gonshor's definition of the exponential function has been proposed in terms of Conway's $\omega$-map. This definition can be generalized [10, Proposition 2.12] by replacing the $\omega$-map by $\Xi_\mathrm{M}$. This yields an alternative exponential function $\exp_\mathrm{M}$ on $\mathbf{No}$ for which $(\mathbf{No}, +, \times, \exp_\mathrm{M})$ is an elementary extension of $(\mathbb{R}, +, \times, \exp)$. The exponentials $\exp$ and $\exp_\mathrm{M}$ coincide on $\mathbb{R}_\mathrm{M} := \mathbb{R}[[\mathfrak{M}_\mathbf{R}]]_{\mathbf{On}}$, but $\exp_\mathrm{M}$ grows faster than $\exp$ on $\mathbb{R}_\mathrm{M}[[\mathfrak{N}_\mathbf{R}]]_{\mathbf{On}}$. It would be interesting to see if the properties of $\mathbf{No}$ as the exponential field of generalized series $(\mathbb{R}_\mathrm{M}[[\mathfrak{N}_\mathbf{R}]]_{\mathbf{On}}, +, \times, \exp_\mathrm{M})$ over $\mathbb{R}_\mathrm{M}$ are similar to those of $(\mathbb{R}[[\mathbf{Mo}]]_{\mathbf{On}}, +, \times, \exp)$ over $\mathbb{R}$.

## 7.4 Exponential groups

Let us now study the action of $\mathcal{E}$ and $\mathcal{E}^*$ on $\mathbf{No}^{>,\succ}$. Given $x \in \mathbf{No}$, recall that one traditionally writes $\lambda_x := \Xi_\mathcal{E} x$ and $\kappa_x := \Xi_{\mathcal{E}^*} x$.

The parameterization $\lambda = \Xi_\mathcal{E}$ of the class $\mathbf{La} := \mathbf{Smp}_\mathcal{E}$ was first given in [11]. It was also shown there that $\mathbf{La}$ coincides with the class of *log-atomic* surreal numbers, which consists of those numbers $x \in \mathbf{No}^{>,\succ}$ such that $\log_n x \in \mathbf{Mo}$ for all $n \in \mathbb{N}$. Such numbers were essential for the definition of well-behaved formal derivations on $\mathbf{No}$. This was first achieved in [11], while building on analogue results in the context of transseries [38, 26].

The structure $\mathbf{K} := \mathbf{Smp}_{\mathcal{E}^*}$ of $\kappa$-numbers was introduced and studied in detail in [30], as an intermediate subclass between fundamental monomials and the log-atomic numbers. It turns out that the structure $\mathbf{K}$ is not big enough to describe all log-atomic numbers. Indeed, it was noticed in [35] that $\mathbf{K} = \mathbf{La} \prec \mathbf{No}_{\succ}$, as a corollary of [3, Proposition 2.5].

**Proposition 7.10.** [3, Proposition 2.5] *For all $x \in \mathbf{No}$, we have*

$$\exp(\lambda_x) = \lambda_{x+1}$$

**Proof.** We rely on the following uniform version of [11, Theorem 3.8(1)] from [3, Lemma 2.4]: if $\mathfrak{m} = \{L|R\}$ is a monomial, where $\mathbb{R} L \subseteq \mathbf{Hull}(L)$ and $\mathbb{R} R \subseteq \mathbf{Hull}(R)$, then

$$\exp(\mathfrak{m}) = \{\mathfrak{m}^{\mathbb{N}}, \exp(L) | \exp(R)\}.$$

In fact, we have $\mathcal{P} \subseteq \mathcal{E} < \{\exp\}$ on $\mathbf{No}^{>,>}$, so $\exp(\mathfrak{m}) > \mathcal{E}\mathfrak{m} \supseteq \mathfrak{m}^{\mathbb{N}}$, and

$$\exp(\mathfrak{m}) = \{\mathcal{E}\mathfrak{m}, \exp(L) | \exp(R)\}. \tag{7.1}$$

Now let $x$ be a number with $\lambda_{u+1} = \exp(\lambda_u)$ for all $u \in x_{\sqsubset}$. Then $x + 1 = \{x, x_L + 1 | x_R + 1\}$. The uniformity of the cut equation for the $\lambda$-map thus yields

$$\begin{aligned}
\lambda_{x+1} &= \{\mathcal{E}\lambda_x, \mathcal{E}\lambda_{x_L+1} | \mathcal{E}\lambda_{x_R+1}\} \\
&= \{\mathcal{E}\lambda_x, \mathcal{E}\exp(\lambda_{x_L}) | \mathcal{E}\exp(\lambda_{x_R})\} \\
&= \{\mathcal{E}\lambda_x, \exp \circ \mathcal{E}\lambda_{x_L} | \exp \circ \mathcal{E}\lambda_{x_R}\} \quad &\text{(since } \exp \circ \mathcal{E} = \mathcal{E} \circ \exp\text{)} \\
&= \exp \lambda_x &\text{(by (7.1))}
\end{aligned}$$

The result follows by induction. $\square$

**Corollary 7.11.** [11] $\mathbf{Smp}_{\mathcal{E}}$ *coincides with the class of log-atomic surreal numbers.*

**Proof.** We have $\log_n \lambda_x = \lambda_{x-n} \in \mathbf{La}$ for all $n \in \mathbb{N}$, whence $\log_n \mathbf{La} \subseteq \mathbf{La} \subseteq \mathbf{Mo}$. This shows that every element of $\mathbf{La}$ is log-atomic.

Conversely, let $\lambda$ be a log-atomic number and assume $\lambda \notin \mathbf{La}$. Note that $\pi_{\mathcal{E}}(\lambda)$ is log-atomic by our previous argument. Assume for instance that $\pi_{\mathcal{E}}(\lambda) < \lambda$. For $n \in \mathbb{N}$, we have $\log_n \pi_{\mathcal{E}}(\lambda) \neq \log_n \lambda$. Since both $\log_n \lambda$ and $\log_n \pi_{\mathcal{E}}(\lambda)$ are monomials, it follows that $\log_n \pi_{\mathcal{E}}(\lambda) \prec \log_n \lambda$. We deduce that $(\exp_n \circ \mathcal{H} \circ \log_n)(\pi_{\mathcal{E}}(\lambda)) < \lambda$, whence $\mathcal{E}\pi_{\mathcal{E}}(\lambda) < \lambda$, which contradicts the defining relation $\pi_{\mathcal{E}}(\lambda) =_{\mathcal{E}} \lambda$. Likewise, $\pi_{\mathcal{E}}(\lambda) > \lambda$ is impossible. We conclude that $\lambda = \pi_{\mathcal{E}}(\lambda) \in \mathbf{Smp}_{\mathcal{E}}$. $\square$

**Proposition 7.12.** [35] *We have* $\mathbf{K} = \mathbf{La} \prec \mathbf{No}_>$.

**Proof.** Following Mantova-Matusinski, we have the following equivalences for any number $x \in \mathbf{No}$:

$$\begin{aligned}
x \in \mathbf{No}_> &\iff x_L + \mathbb{N} < x < x_R - \mathbb{N} \\
&\iff \exp_{\mathbb{N}}(\lambda_{x_L}) < \lambda_x < \log_{\mathbb{N}}(\lambda_{x_R}) \\
&\iff \exp_{\mathbb{N}}(\mathcal{E}\lambda_{x_L}) < \lambda_x < \log_{\mathbb{N}}(\mathcal{E}\lambda_{x_R}) \\
&\iff \mathcal{E}^*(\lambda_{x_L}) < \lambda_x < \mathcal{E}^*(\lambda_{x_R}) \\
&\iff \lambda_x \in \mathbf{K}.
\end{aligned}$$
$\square$

**Corollary 7.13.** $\mathbf{K}$ *is sharp in* $\mathbf{No}^{>,>}$.

**Proof.** Let $\kappa \in \mathbf{K}$, $\kappa' \in \kappa_L^{\mathbf{K}}$, and $\kappa'' \in \kappa_R^{\mathbf{K}}$. There are unique numbers $\theta, \theta', \theta'' \in \mathbf{No}_>$ with $\kappa = \lambda_{\theta}$, $\kappa' = \lambda_{\theta'}$, and $\kappa'' = \lambda_{\theta''}$. Let $n \in \mathbb{N}$. We have $\theta > \theta' + \mathbb{N}$ and $\theta' + n = \{\theta'_L + n, \theta' + n_L | \theta'_R + n\}$, where

$$\theta_R + n > \theta > \theta' + n > \theta'_L + n \cup \theta' + n_L,$$

so $\theta \sqsupseteq \theta' + n$. We deduce that $\exp_{\mathbb{N}}(\kappa') = \lambda_{\theta' + \mathbb{N}} \sqsubseteq \kappa$. Symmetric arguments yield $\log_{\mathbb{N}}(\kappa'') = \lambda_{\theta'' - \mathbb{N}} \sqsubseteq \kappa$. Since $\exp_{\mathbb{N}}(\kappa')$ is cofinal in $\mathcal{E}^*[\kappa']$ and $\log_{\mathbb{N}}(\kappa'')$ is coinitial in $\mathcal{E}^*[\kappa'']$, this proves that $\mathbf{K} = \mathbf{Smp}_{\mathcal{E}^*}$ is sharp. $\square$

On the other hand, the class **La** is not sharp:

**Proposition 7.14.** *The structure **La** is not sharp in* $\mathbf{No}^{>,\succ}$.

**Proof.** Given $r \in \mathbb{R}^>$ and $x \in \mathbf{No}^{>,\succ}$, we have
$$x^{e^r} = (\exp_2 \circ T_r \circ \log_2)(x) < (\exp_3 \circ T_2 \circ \log_3)(x).$$

We deduce that the element $(\exp_3 \circ T_2 \circ \log_3)(\omega)$ of $\mathcal{E}[\omega]$ is a strict upper bound for $\dot\omega^{e^{\mathbb{R}^>}}$ and hence for $\dot\omega^{\mathbb{N}}$. Note that $\lambda_1 = \dot\omega^\omega$, so $\dot\omega^{\mathbb{N}}$ is cofinal in $(\lambda_1)_L$. We have $(\lambda_1)_L^{\mathbf{La}} = \{\lambda_0\} = \{\omega\}$, so $(\lambda_1)_L$ is not cofinal in $\mathcal{E}[(\lambda_1)_L^{\mathbf{La}}]$. This means that the defining partition of **La** is not sharp. □

# 8　Nested surreal numbers

## 8.1　Nested transseries and surreal numbers

The study of generalized transseries solutions to functional equations was started in [17, 26]. It is well known that non-trivial solutions of the functional equation $E(x+1) = \exp E(x)$ grow faster than any iterated exponential. This motivates the introduction of "hyperseries" [17, 38, 2, 16] as a generalization of transseries that allows for transfinite iterates of exponentiation and logarithm. In [26, section 2.7.1], it was pointed out that functional equations of the kind
$$f(x) = \sqrt{x} + e^{f(\log x)} \tag{8.1}$$

admit natural symbolic solutions of the form
$$f(x) = \sqrt{x} + e^{\sqrt{\log x} + e^{\sqrt{\log\log x} + e^{\cdot^{\cdot^{\cdot}}}}}. \tag{8.2}$$

The formal calculus with this kind of expressions requires a second extension of Écalle's original theory from [17] with so-called "nested transseries". In our context, it is also natural to study those surreal numbers
$$y = f(\omega) = \sqrt{\omega} + e^{\sqrt{\log \omega} + e^{\sqrt{\log\log \omega} + e^{\cdot^{\cdot^{\cdot}}}}} \tag{8.3}$$

that are obtained by substituting $\omega$ for $x$ in such a generalized transseries. More specifically, one may wonder whether there exist sequences $(y_i)_{i \in \mathbb{N}} \in \mathbf{No}_>^{\mathbb{N}}$ with
$$y_i = \sqrt{\log_i \omega} + e^{y_{i+1}},$$

for all $i \in \mathbb{N}$. In this section, we will show that the class of such numbers actually forms a surreal substructure. This shows in particular that expressions of the form (8.2) or (8.3) are highly ambiguous and therefore somewhat misleading.

In order to develop a sound calculus for nested transseries and surreal numbers such as (8.2) and (8.3) it is crucial to decide which expressions of the form (8.2) should be considered to be well-formed. For instance, the functional equation
$$g(x) = \sqrt{x} + e^{g(\log x)} + \log x \tag{8.4}$$

admits a "natural" solution
$$g(x) = \sqrt{x} + e^{\sqrt{\log x} + e^{\sqrt{\log_2 x} + e^{\cdot^{\cdot^{\cdot}} + \log_3 x} + \log_2 x}} + \log x. \tag{8.5}$$

However, such expressions do not behave well for basic calculus operations. For instance, the syntactic derivative of (8.5) is given by

$$g'(x) = \frac{1}{2\sqrt{x}} + \frac{1}{x} + \left(\frac{1}{2x\sqrt{\log x}} + \frac{1}{x \log x}\right) e^{g(\log x)} +$$
$$\left(\frac{1}{2x \log x \sqrt{\log_2 x}} + \frac{1}{x \log x \log_2 x}\right) e^{g(\log x)} e^{g(\log_2 x)} + \cdots.$$

However, the sum

$$\frac{1}{x} + \frac{1}{x \log x} e^{g(\log x)} + \frac{1}{x \log x \log_2 x} e^{g(\log x)} e^{g(\log_2 x)} + \cdots$$

does not converge in the sense of section 2.3. Fortunately, as pointed out in [26, section 2.7.1], the equation (8.4) is a perturbation of (8.1) and its solutions can naturally be expressed in terms of $f$.

The above counterexample led the second author to introduce the abstract notion of so-called *fields of transseries* [27] which excludes transseries such as (8.5). Generalizing the combinatorial ideas from [26], this enabled him and his student Schmeling to construct derivations and right compositions on fields of transseries [38]. This theory reappeared crucially in Berarducci and Mantova's construction of a well-behaved derivation $\partial_{\mathrm{BM}}$ on **No** [11]. Indeed, one of the main ingredients of their construction is the proof [11, Theorem 8.10] that **No** is a field of transseries in the sense of [27, 38]. In particular, it satisfies the following condition:

**T4.** Let $(\mathfrak{m}_i)_{i \in \mathbb{N}} \in \mathbf{Mo}^{\mathbb{N}}$ be a sequence of monomials with $\mathfrak{m}_{i+1} \in \operatorname{supp} \log \mathfrak{m}_i$ for all $i$. Then there exists an $i_0 \in \mathbb{N}$ with

$$\forall i \geqslant i_0, \quad \mathfrak{m}_{i+1} \preccurlyeq \operatorname{supp} \log \mathfrak{m}_i \wedge (\log \mathfrak{m}_i)_{\mathfrak{m}_{i+1}} \in \{-1, 1\}.$$

This condition can be regarded as a formal translation of the idea that all surreal numbers should be "well nested". In particular, it rules out the existence of surreal numbers of the form

$$\sqrt{\omega} + e^{\sqrt{\log \omega} + e^{\sqrt{\log_2 \omega} + e^{\cdot^{\cdot^{\cdot} + \log_3 \omega}} + \log_2 \omega}} + \log \omega.$$

## 8.2 Admissible sequences

Given sequences $(\varphi_i)_{i \in \mathbb{N}} \in \mathbf{No}^{\mathbb{N}}$ and $(\epsilon_i)_{i \in \mathbb{N}} \in \{-1, 1\}^{\mathbb{N}}$, let us study how to give a meaning to expressions of the type

$$\varphi_0 + \epsilon_0 e^{\varphi_1 + \epsilon_1 e^{\varphi_2 + \epsilon_2 e^{\cdot^{\cdot^{\cdot}}}}}. \tag{8.6}$$

In this subsection, we start with the determination of lower and upper bounds for (8.6). We say that $(\varphi, \epsilon)$ is a *signed sequence* if

**SS1.** $\varphi_i \geqslant 0$ for all $i \geqslant 2$.

**SS2.** $\varphi_i = 0 \Longrightarrow \epsilon_i = 1$ for all $i \geqslant 2$.

**SS3.** $\varphi_i > 0$ for infinitely many $i$.

**SS4.** $\varphi_i \in \mathbf{No}_{\succ}$ for all $i \geqslant 1$.

In that case, we may define a signed sequence $(\varphi_{\nearrow k}, \epsilon_{\nearrow k})$ for every $k \in \mathbb{N}$ by taking $(\varphi_{\nearrow k})_i := \varphi_{k+i}$ and $(\epsilon_{\nearrow k})_i := \epsilon_{k+i}$ for all $i \in \mathbb{N}$.

Assume that $(\varphi, \epsilon)$ is a fixed signed sequence. For all $i, j \in \mathbb{N}$ with $i \leqslant j$, we define functions $\Phi_i, \Phi_{i;}, \Phi_{j;i}: \mathbf{No} \longrightarrow \mathbf{No}$ by

$$\Phi_i(x) := \varphi_i + \epsilon_i e^x.$$

$$\Phi_{i;}(x) := (\Phi_0 \circ \cdots \circ \Phi_{i-1})(x) = \varphi_0 + \epsilon_0 e^{\varphi_1 + \epsilon_1 e^{\varphi_2 + \epsilon_2 e^{\cdot^{\cdot^{\varphi_{i-1} + \epsilon_{i-1} e^x}}}}}$$

$$\Phi_{j;i}(x) := (\Phi_i \circ \cdots \circ \Phi_{j-1})(x) = \varphi_i + \epsilon_i e^{\varphi_{i+1} + \epsilon_{i+1} e^{\varphi_{i+2} + \epsilon_{i+2} e^{\cdot^{\cdot^{\varphi_{j-1} + \epsilon_{j-1} e^x}}}}}.$$

By convention, we understand that $\Phi_{0;}(x) = x$ and $\Phi_{j;i}(x) = x$ whenever $i = j$.

Writing $\epsilon_{i;} := \epsilon_{;i} := \epsilon_0 \cdots \epsilon_{i-1}$ and $\epsilon_{j;i} := \epsilon_{i;j} := \epsilon_i \cdots \epsilon_{j-1}$, we notice that $\Phi_i, \Phi_{i;}$, and $\Phi_{j;i}$ are strictly increasing if $\epsilon_i = 1$, $\epsilon_{i;} = 1$, and $\epsilon_{j;i} = 1$, respectively, and strictly decreasing in the contrary case. We will write $\Phi_{;i}$ and $\Phi_{i;j}$ for the partial inverses of $\Phi_{i;}$ and $\Phi_{j;i}$. We will also use the abbreviations

$$\begin{aligned} x_{i;} &:= \Phi_{i;}(x) & x_{j;i} &:= \Phi_{j;i}(x) \\ x_{;i} &:= \Phi_{;i}(x) & x_{i;j} &:= \Phi_{i;j}(x). \end{aligned}$$

For instance, we have

$$x_{1;3} = \varphi_1 + \epsilon_1 e^{\varphi_2 + \epsilon_2 e^x}$$

for all $x$ and

$$x_{;1} = \log \frac{x - \varphi_0}{\epsilon_0},$$

whenever $\frac{x - \varphi_0}{\epsilon_0} > 0$. For all $i \in \mathbb{N}$, we next define

$$\begin{aligned} L_{i;} &:= (\varphi_i - \epsilon_{;i} \mathbb{R}^{>} \operatorname{supp} \varphi_i)_{i;} & L &:= \bigcup_{i \in \mathbb{N}} L_{i;} \\ R_{i;} &:= (\varphi_i + \epsilon_{;i} \mathbb{R}^{>} \operatorname{supp} \varphi_i)_{i;} & R &:= \bigcup_{i \in \mathbb{N}} R_{i;} \end{aligned}$$

We finally define

$$\mathbf{S} := \{x \in \mathbf{No} : \forall i \in \mathbb{N}, x_{;i} - \varphi_i \prec \operatorname{supp} \varphi_i\}.$$

In the remainder of this section, the signed sequence $(\varphi, \epsilon)$ will mostly remain fixed. In the rare cases when $(\varphi, \epsilon)$ needs to be varied, we will use subscripts, e.g. by writing $\mathbf{S}_{\varphi, \epsilon}$ instead of $\mathbf{S}$. For each $k \in \mathbb{N}$ we also write $\mathbf{S}_{\nearrow k} := \mathbf{S}_{\varphi_{\nearrow k}, \epsilon_{\nearrow k}}$.

**Lemma 8.1.** *If $x \in \mathbf{S}$ or $x \in (L | R)$, then $x_{;i}$ is well defined for all $i \in \mathbb{N}$.*

**Proof.** If $x \in \mathbf{S}$, then the definition of $\mathbf{S}$ implicitly assumes that $x_{;i}$ is well defined for all $i \in \mathbb{N}$. If $x \in (L | R)$, so in particular $L < R$, then let us prove the lemma by induction on $i$. The result clearly holds for $i = 0$. Assuming that $x_{;i}$ is well defined, let $j > i$ be minimal such that $\varphi_j \neq 0$. Applying $\Phi_{;i}$ to the inequality

$$L_{j;} < x < R_{j;},$$

we obtain

$$\epsilon_{;i} (L_{j;})_{;i} < \epsilon_{;i} x_{;i} < \epsilon_{;i} (R_{j;})_{;i}.$$

By definition, we have

$$\begin{aligned} (L_{j;})_{;i} &= \varphi_i + \epsilon_i \exp_{j-i}(\varphi_j - \epsilon_{;j} \mathbb{R}^{>} \operatorname{supp} \varphi_j) \\ (R_{j;})_{;i} &= \varphi_i + \epsilon_i \exp_{j-i}(\varphi_j + \epsilon_{;j} \mathbb{R}^{>} \operatorname{supp} \varphi_j), \end{aligned}$$

whence

$$\epsilon_{;i+1} \exp_{j-i}(\varphi_j - \epsilon_{;j} \mathbb{R}^{>} \operatorname{supp} \varphi_j) < \epsilon_{;i+1} \frac{x_{;i} - \varphi_i}{\epsilon_i} < \epsilon_{;i+1} \exp_{j-i}(\varphi_j + \epsilon_{;j} \mathbb{R}^{>} \operatorname{supp} \varphi_j).$$

Both in the cases when $\epsilon_{;i+1}=1$ and when $\epsilon_{;i+1}=-1$, it follows that $(x_{;i}-\varphi_i)/\epsilon_i$ is bounded from below by the exponential of a surreal number, whence $(x_{;i}-\varphi_i)/\epsilon_i>0$. In particular, $x_{;i+1}=\log((x_{;i}-\varphi_i)/\epsilon_i)$ is well defined. This completes the induction. □

**Proposition 8.2.** *We have* $\mathbf{S}=(L|R)$.

**Proof.** Let $x\in\mathbf{S}$ and $i\in\mathbb{N}$. If $\epsilon_{;i}=1$, then $\Phi_{i;}$ is strictly increasing, whence

$$L_{i;}<x<R_{i;} \iff \varphi_i-\mathbb{R}^{>}\operatorname{supp}\varphi_i<x_{;i}<\varphi_i+\mathbb{R}^{>}\operatorname{supp}\varphi_i$$
$$\iff -\mathbb{R}^{>}\operatorname{supp}\varphi_i<x_{;i}-\varphi_i<\mathbb{R}^{>}\operatorname{supp}\varphi_i$$
$$\iff x_{;i}-\varphi_i\prec\operatorname{supp}\varphi_i.$$

Otherwise $\Phi_{i;}$ is strictly decreasing, whence

$$L_{i;}<x<R_{i;} \iff \varphi_i+\mathbb{R}^{>}\operatorname{supp}\varphi_i>x_{;i}>\varphi_i-\mathbb{R}^{>}\operatorname{supp}\varphi_i$$
$$\iff \mathbb{R}^{>}\operatorname{supp}\varphi_i>x_{;i}-\varphi_i>-\mathbb{R}^{>}\operatorname{supp}\varphi_i$$
$$\iff x_{;i}-\varphi_i\prec\operatorname{supp}\varphi_i.$$

In both cases, we conclude that $L_{i;}<x<R_{i;}$ if and only if $x_{;i}-\varphi_i\prec\operatorname{supp}\varphi_i$. Since this equivalence holds for all $i\in\mathbb{N}$, the result follows. □

We say that the signed sequence $(\varphi,\epsilon)$ is *admissible* if

**AS.** $L<R$.

**Proposition 8.3.** *The following statements are equivalent.*

a) $(\varphi,\epsilon)$ *is admissible.*

b) $\mathbf{S}$ *is a surreal substructure.*

c) $\forall i\in\mathbb{N},\forall\mathfrak{m}\in\operatorname{supp}\varphi_i,\forall j>i,\exists\psi\in\mathbf{No}^{<\operatorname{supp}\varphi_j},\mathfrak{m}\succ(\varphi_j\mathbin{+\mkern-10mu+}\psi)_{j;i}-\varphi_i.$

**Proof.** We have b) $\Longrightarrow$ a) by the previous proposition. If $(\varphi,\epsilon)$ is admissible, then $\mathbf{S}=(L|R)$ is a surreal substructure by Proposition 4.18(b). We also obtain c) by taking $\psi\in\mathbf{S}_{\nearrow j}-\varphi_j$. Indeed, we have $(\varphi_j\mathbin{+\mkern-10mu+}\psi)_{j;i}-\varphi_i\in(\mathbf{S}_{\nearrow j})_{j;i}-\varphi_i\subseteq\mathbf{S}_{;i}-\varphi_i$, whence $(\varphi_j\mathbin{+\mkern-10mu+}\psi)_{j;i}-\varphi_i\prec\operatorname{supp}\varphi_i$, by the definition of $\mathbf{S}$. The definition of $\mathbf{S}$ also yields $\psi\prec\operatorname{supp}\varphi_j$. Assume finally that c) is satisfied and let us prove a).

Let $i,j\in\mathbb{N}$. If $i=j$, then $L_{i;}<R_{i;}$ follows by definition and strict monotonicity of the function $\Phi_{i;}$. Assume that $i<j$. Let $\mathfrak{m}\in\operatorname{supp}\varphi_i$ and consider a $\psi\in\mathbf{No}^{<\operatorname{supp}\varphi_j}$ with $(\varphi_j\mathbin{+\mkern-10mu+}\psi)_{j;i}-\varphi_i\prec\mathfrak{m}$. Such a $\psi$ exists by c) and the class $\mathbf{C}_{\mathfrak{m}}$ of such numbers $\psi$ is a convex surreal substructure by Proposition 4.18(d). Moreover the family $(\mathbf{C}_{\mathfrak{m}})_{\mathfrak{m}\in\operatorname{supp}\varphi_i}$ is decreasing on $(\operatorname{supp}\varphi_i,\succcurlyeq)$ so by Proposition 4.18(e), its intersection is non-empty. Given $y$ in this intersection, we have $L_{i;}<((\varphi_j\mathbin{+\mkern-10mu+}y)_{j;i})_{i;}=(\varphi_j\mathbin{+\mkern-10mu+}y)_{j;}$, since $(\varphi_j\mathbin{+\mkern-10mu+}y)_{j;i}-\varphi_i\prec\operatorname{supp}\varphi_i$. Similarly, $((\varphi_j\mathbin{+\mkern-10mu+}y)_{j;i})_{i;}=(\varphi_j\mathbin{+\mkern-10mu+}y)_{j;}<R_j$, since $y=(\varphi_j\mathbin{+\mkern-10mu+}y)-\varphi_j\prec\operatorname{supp}\varphi_j$. This shows that $L_{i;}<(\varphi_j\mathbin{+\mkern-10mu+}y)_{j;}<R_{j;}$. By symmetry, we obtain the same conclusion if $i>j$, i.e. $(\varphi,\epsilon)$ is admissible. □

## 8.3 Nested sequences

Let $(\varphi,\epsilon)$ be a fixed admissible sequence. Now that we have described lower and upper bounds $L$ and $R$ for expressions of the form (8.6), our next goal is to determine those elements $y\in\mathbf{S}=(L|R)$ such that

$$\operatorname{supp}\varphi_i\succ\frac{y_{;i}-\varphi_i}{\epsilon_i}\in\mathbf{Mo}$$

for all $i \in \mathbb{N}$. Such elements are called *nested surreal numbers* and we denote by $\mathbf{Ne} = \mathbf{Ne}_{\varphi,\epsilon}$ the class of nested surreal numbers with respect to our fixed admissible sequence $(\varphi, \epsilon)$.

It turns out that not all admissible sequences $(\varphi, \epsilon)$ give rise to nested surreal numbers (see Example 8.14 below). We say that $(\varphi, \epsilon)$ is *nested* if

**NS.** $\operatorname{supp} \varphi_i \succ e^{\mathbf{S}_{\nearrow(i+1)}}$, for all $i \in \mathbb{N}$.

The main objective of this subsection is to show that $\mathbf{Ne}$ is a surreal substructure whenever $(\varphi, \epsilon)$ is nested (in particular, $\mathbf{Ne}$ is non-empty). In the next subsection, we will give various examples and sufficient conditions for **NS** to be satisfied.

We will say that $(\varphi, \epsilon)$ is *large* if we have $\varphi_1 > 0$ or $(\varphi_1, \epsilon_1) = (0, 1)$. Notice that the admissible sequences $(\varphi_{\nearrow i}, \epsilon_{\nearrow i})$ for $i > 0$ are always large. Let us first show how to reduce the general case to the case when $(\varphi, \epsilon)$ is large. Assuming that $(\varphi, \epsilon)$ is not large, let $(\varphi', \epsilon')$ be the large nested sequence with $(\varphi'_0, \epsilon'_0) = (0, 1)$, $(\varphi'_1, \epsilon'_1) = (-\varphi_1, -\epsilon_1)$, and $\varphi'_i = \varphi_i$ for $i \geqslant 2$. Assume that we know how to show that $\mathbf{Ne}_{\varphi', \epsilon'}$ is a surreal substructure of $\mathbf{Mo}$. Writing $S(a) := -a$ and $I(a) = a^{-1}$, we have $\Xi_{\mathbf{Mo}} \circ S = I \circ \Xi_{\mathbf{Mo}}$, whence $I$ induces a strictly decreases self-$\sqsubseteq$-embedding on $\mathbf{Mo}$. It follows that the function $x \longmapsto I \circ \Xi_{\mathbf{Ne}_{\varphi', \epsilon'}}(-x)$ is an embedding of $\mathbf{Mo}$ into itself. Hence the range $(\mathbf{Ne}_{\varphi', \epsilon'})^{-1}$ of this mapping is a surreal substructure, and so is $\mathbf{Ne} = \varphi_0 + \epsilon_0 (\mathbf{Ne}_{\varphi', \epsilon'})^{-1}$.

In the remainder of this section, let $(\varphi, \epsilon)$ be a fixed large nested sequence.

**Lemma 8.4.** *For $x, y \in \mathbf{S}$, we have*
$$(x_{;1} - \varphi_1) / \epsilon_1 =_{\mathcal{E}} (y_{;1} - \varphi_1) / \epsilon_1.$$

**Proof.** Choose $i \in \mathbb{N}^>$ minimal with $\varphi_i \neq 0$. We have $x_{;i} - \varphi_i, y_{;i} - \varphi_i \prec \operatorname{supp} \varphi_i$, whence
$$\tfrac{1}{2} y_{;i} < x_{;i} < 2 y_{;i}$$
and
$$\exp_{i-1} y_{;i} =_{\mathcal{E}} \exp_i(\tfrac{1}{2} y_{;i}) < \exp_i x_{;i} < \exp_i(2 y_{;i}) =_{\mathcal{E}} \exp_i y_{;i}$$

We observe that $(x_{;1} - \varphi_1) / \epsilon_1 = \exp_{i-1} x_{;i}$ and $(y_{;1} - \varphi_1) / \epsilon_1 = \exp_i y_{;i}$. By convexity of $\mathcal{E}[(y_{;1} - \varphi_1)/\epsilon_1]$, we have $(x_{;1} - \varphi_1)/\epsilon_1 \in \mathcal{E}[(y_{;1} - \varphi_1)/\epsilon_1]$, whence the result. □

**Lemma 8.5.** *We have a $\sqsubseteq$-embedding*
$$\Phi_{1;} : \mathbf{S}_{\nearrow 1} \cap (\varphi_1 + \epsilon_1 \mathbf{Mo}) \longrightarrow \mathbf{S}.$$

**Proof.** Recall that $\Phi_{1;}(x) = \varphi_0 + \epsilon_0 e^x$ for all $x \in \mathbf{No}$. Let us first show that $\mathbf{U} := \mathbf{S}_{\nearrow 1} \cap (\varphi_1 + \epsilon_1 \mathbf{Mo})$ is a surreal substructure. By **NS**, we have $\mathbf{S}_{\nearrow 1} = \varphi_1 + \epsilon_1 e^{\mathbf{S}_{\nearrow 2}}$. Writing $\mathbf{S}_{\nearrow 2} =: (L_{\nearrow 2} | R_{\nearrow 2})$, as for $\mathbf{S}$, we observe that $L_{\nearrow 2}$ and $R_{\nearrow 2}$ are sets of purely infinite numbers, respectively without maximum and minimum. By Proposition 4.18(b), it follows that $\mathbf{S}_{\nearrow 2} \cap \mathbf{No}_{\succ} = (L_{\nearrow 2} | R_{\nearrow 2})_{\mathbf{No}_{\succ}}$ is a convex surreal substructure of $\mathbf{No}_{\succ}$. By Proposition 4.18(d), we deduce that $\mathbf{U} = \varphi_1 + \epsilon_1 e^{\mathbf{S}_{\nearrow 2} \cap \mathbf{No}_{\succ}}$ is a convex surreal substructure of $\varphi_1 + \epsilon_1 \mathbf{Mo}^{<\operatorname{supp} \varphi_1}$.

By Proposition 4.28 and **NS**, the function $x \longmapsto \varphi_0 + \epsilon_0 x$ is a $\sqsubseteq$-embedding on $e^{\mathbf{U}}$, so it remains to be shown that $\exp$ is a $\sqsubseteq$-embedding on $\mathbf{U}$. Towards this, consider numbers $u, v \in \mathbf{U}$ with $u \sqsubseteq v$. Since $u, v \in \varphi_1 + \epsilon_1 \mathbf{Mo}^{<\operatorname{supp} \varphi_1}$, Proposition 4.28 implies that $u = \varphi_1 + \epsilon_1 \mathfrak{u}$ and $v = \varphi_1 + \epsilon_1 \mathfrak{v}$ for certain infinite monomials $\mathfrak{u}$ and $\mathfrak{v}$ with $\mathfrak{u} \sqsubseteq \mathfrak{v}$.

Consider $\mathfrak{m} \in \mathbf{Mo}^>$. The cuts $(\mathbb{R}^> \mathfrak{m}_L^{\mathbf{Mo}} | \mathbb{R}^> \mathfrak{m}_R^{\mathbf{Mo}})$ and $(\mathfrak{m}_L | \mathfrak{m}_R)$ are mutually cofinal. Given (7.1), it follows that
$$\forall \mathfrak{m} \in \mathbf{Mo}^>, e^{\mathfrak{m}} = \{\mathcal{E}\mathfrak{m}, \mathcal{D}e^{\mathfrak{m}_L^{\mathbf{Mo}}} \mid \mathcal{D}e^{\mathfrak{m}_R^{\mathbf{Mo}}}\}.$$

Proposition 4.36 therefore implies that exp is a $\sqsubseteq$-embedding on $\mathcal{E}[\mathfrak{m}] \cap \mathbf{Mo}^{>}$ for every $\mathfrak{m} \in \mathbf{Mo}^{>}$. Using Lemma 8.4, we deduce that $e^{\mathfrak{u}} \sqsubseteq e^{\mathfrak{v}}$.

Just before Lemma 8.4, we already noticed that $I: \mathbf{Mo} \longrightarrow \mathbf{Mo}; \mathfrak{m} \longmapsto \mathfrak{m}^{-1}$ is a $\sqsubseteq$-embedding. Since $\mathfrak{u}$ and $\mathfrak{v}$ are monomials, it follows that $e^{\epsilon_1 \mathfrak{u}} = (e^{\mathfrak{u}})^{\epsilon_1} \sqsubseteq (e^{\mathfrak{v}})^{\epsilon_1} = e^{\epsilon_1 \mathfrak{v}}$. By [11, Proposition 4.23], we conclude that $e^u = e^{\varphi_1 + \epsilon_1 \mathfrak{u}} \sqsubseteq e^{\varphi_1 + \epsilon_1 \mathfrak{v}} = e^v$. $\square$

In order to show that $\mathbf{Ne}$ is a surreal substructure, let us now introduce a suitable function group $\mathcal{G}$ acting on $\mathbf{S}$. At a second stage, we will show that $\mathbf{Ne} = \mathbf{Smp}_{\mathcal{G}}$. Theorem 6.20 then implies that $\mathbf{Ne}$ is a surreal substructure.

**Lemma 8.6.** *Given $x \in \mathbf{S}$ and $r \in \mathbb{R}^{>}$, we have $\varphi_0 + r(x - \varphi_0) \in \mathbf{S}$.*

**Proof.** Let $y = \varphi_0 + r(x - \varphi_0)$. Let us show by induction on $i \in \mathbb{N}$ that
$$L_{i;} < y < R_{i;}$$
and $y_{;i} - x_{;i} \preccurlyeq 1$ whenever $i \geqslant 1$. This is clear for $i = 0$, so assume $i > 0$. If $i = 1$, then $y_{;i} - x_{;i} = \log r \preccurlyeq 1$. If $i > 1$, then the induction hypothesis yields
$$y_{;i} - x_{;i} = \log \frac{y_{;i-1} - \varphi_{i-1}}{x_{;i-1} - \varphi_{i-1}} = \log\left(1 + \frac{y_{;i-1} - x_{;i-1}}{x_{;i-1} - \varphi_{i-1}}\right) = \log(1 + o(1)) = o(1) \preccurlyeq 1.$$

By **NS**, we also have $\operatorname{supp} \varphi_i \succ e^{x_{;i+1}}$, whence
$$\varphi_i - \mathbb{R}^{>} \operatorname{supp} \varphi_i < \varphi_i + \epsilon_i e^{x_{;i+1}} < \varphi_i + \mathbb{R}^{>} \operatorname{supp} \varphi_i.$$

We have $\operatorname{supp} \varphi_i \succ 1$ by **SS4**. Since $y_{;i} = x_{;i} + O(1) = \varphi_i + \epsilon_i e^{x_{;i+1}} + O(1)$, this yields
$$\varphi_i - \mathbb{R}^{>} \operatorname{supp} \varphi_i < y_{;i} < \varphi_i + \mathbb{R}^{>} \operatorname{supp} \varphi_i.$$

Applying $\Phi_{i;}$, we conclude that $L_{i;} < y < R_{i;}$, which completes our proof by induction. $\square$

The lemma implies that $\mathbf{S}_{\nearrow i} - \varphi_i$ is closed under the action of $\mathcal{H}$ for all $i \in \mathbb{N}$. This allows us to define a strictly increasing bijection
$$\Psi_{i,r}: \mathbf{S} \longrightarrow \mathbf{S}; x \longmapsto (\varphi_i + r(x_{;i} - \varphi_i))_{i;}$$
for all $i \in \mathbb{N}$ and $r \in \mathbb{R}^{>}$. We take
$$\mathcal{G} := \langle \Psi_{i,r} : r \in \mathbb{R}^{>}, i \in \mathbb{N} \rangle$$
to be the function group generated by these functions. As usual, we will write $\mathcal{G}_{\nearrow i}$ for the function group obtained by applying this definition for $(\varphi_{\nearrow i}, \epsilon_{\nearrow i})$ instead of $(\varphi, \epsilon)$.

**Lemma 8.7.** *Given $x \in \mathbf{S}$, we have:*

a) *For each $i > 0$, the set $\Psi_{i, \mathbb{R}^{>}}(x)$ contains strict upper and lower bounds for $\Psi_{i-1, \mathbb{R}^{>}}(x)$.*

b) *The set $\{\Psi_{i,r}(x) : r \in \mathbb{R}^{>}, i \in \mathbb{N}, i > j\}$ is cofinal and coinitial in $\mathcal{G}[x]$ for all $j \in \mathbb{N}$.*

c) *For $y \in \mathcal{G}[x]$, we have $\varphi_0 + \epsilon_0 \mathfrak{d}_{y - \varphi_0} \in \mathcal{G}[x]$, whence $\mathcal{G}[x]^{\bullet} \in \varphi_0 + \epsilon_0 \mathbf{Mo}$.*

d) *$(\mathcal{G}_{\nearrow 1}[x_{;1}])_{1;} = \mathcal{G}[x]$.*

e) *$\mathcal{G}[x]^{\bullet}_{;1} = \mathcal{G}_{\nearrow 1}[x_{;1}]^{\bullet}$.*

**Proof.**

a) The number $x_{;i+1}$ is positive infinite, so we have
$$\varphi_i + 2^{-\epsilon_i}(x_{;i} - \varphi_i) + \mathbb{Z} < x_{;i} < \varphi_i + 2^{\epsilon_i}(x_{;i} - \varphi_i) + \mathbb{Z},$$

whence
$$e^{\varphi_i + 2^{-\epsilon_i}(x_{;i} - \varphi_i)} \prec e^{x_{;i}} \prec e^{\varphi_i + 2^{\epsilon_i}(x_{;i} - \varphi_i)}.$$

If $\epsilon_{i-1} = 1$, then it follows that
$$\varphi_{i-1} + e^{\varphi_i + 2^{-\epsilon_i}(x_{;i} - \varphi_i)} < \varphi_{i-1} + \mathbb{R}^> e^{x_{;i}} < \varphi_{i-1} + e^{\varphi_i + 2^{\epsilon_i}(x_{;i} - \varphi_i)}.$$

Applying $\Phi_{i-1;}$, we obtain
$$\Psi_{i,2^{-\epsilon_i}}(x) < \Psi_{i-1,\mathbb{R}^>}(x) < \Psi_{i,2^{\epsilon_i}}(x).$$

If $\epsilon_{i-1} = 1$, then a similar reasoning yields
$$\Psi_{i,2^{\epsilon_i}}(x) < \Psi_{i-1,\mathbb{R}^>}(x) < \Psi_{i,2^{-\epsilon_i}}(x).$$

In both cases, this shows that $\Psi_{i,\mathbb{R}^>}(x)$ contains strict upper and lower bounds for $\Psi_{i,\mathbb{R}^>}(x)$.

b) By induction on $j \in \mathbb{N}$, let us show that $\Psi_{j,\mathbb{R}^>}$ is strictly cofinal and coinitial with respect to $\mathcal{G}_{<j} := \langle \Psi_{i,r} : i < j, r \in \mathbb{R}^> \rangle \subseteq \mathcal{G}$. Note that $\mathcal{G}_{<0} = \{\mathrm{id}_{\mathbf{S}}\}$. In view of a), this clearly holds for $j = 0$.

Assuming that this assertion holds for a given $j \in \mathbb{N}$, let us first show that $\Psi_{j,\mathbb{R}^>}$ is cofinal with respect to $\mathcal{G}_{\leqslant j} := \mathcal{G}_{<j+1}$. Given $x' = (\Psi_{j,r_1} \circ \gamma_1 \circ \cdots \circ \Psi_{j,r_l} \circ \gamma_l)(x)$ with $\gamma_1, \ldots, \gamma_l \in \mathcal{G}_{<j}$, we must show that $x' < \Psi_{j,s}(x)$ for some $s \in \mathbb{R}^>$. Using a second induction on $l$, we may find an $s' \in \mathbb{R}^>$ with $y := (\Psi_{j,r_2} \circ \gamma_2 \circ \cdots \circ \Psi_{j,r_l} \circ \gamma_l)(x) < \Psi_{j,s'}(x)$. Using the induction hypothesis on $j$, it follows that $\gamma_1(y) < \Psi_{j,t}(y) < \Psi_{j,s't}(x)$ for some $t \in \mathbb{R}^>$, whence $x' = \Psi_{j,r_1}(\gamma_1(y)) < \Psi_{j,r_1 s' t}(x)$.

In a similar way, one shows that $\Psi_{j,\mathbb{R}^>}$ is strictly coinitial with respect to $\mathcal{G}_{\leqslant j}$. Applying a) for $i = j+1$, it also follows that $\Psi_{j+1,\mathbb{R}^>}$ is strictly coinitial with respect to $\mathcal{G}_{<j+1}$. We conclude by induction.

c) We have $\varphi_0 + \mathbb{R}^> (y - \varphi_0) \subseteq \mathcal{G}[x]$, whence $\varphi_0 + \epsilon_0 \mathfrak{d}_{y-\varphi_0} \in \mathbf{Hull}(\varphi_0 + \mathbb{R}^> (y - \varphi_0)) \subseteq \mathcal{G}[x]$.

d) Applying b) to $j = 0$ yields $\mathcal{G}[x] := \mathbf{Hull}(\Psi_{i,\mathbb{R}^>}(x) : i \geqslant 1)$. Consequently,
$$\begin{aligned}(\mathcal{G}_{\nearrow 1}[x_{;1}])_{1;} &= \mathbf{Hull}(((\varphi_{i+1} + \mathbb{R}^>(x_{;(i+1)} - \varphi_{i+1}))_{i+1;1})_{1;} : i \geqslant 0) \\ &= \mathbf{Hull}(\Psi_{i,\mathbb{R}^>}(x) : i \geqslant 1) \\ &= \mathcal{G}[x].\end{aligned}$$

e) Let $a = \mathcal{G}[x]^\bullet$ and $b = \mathcal{G}_{\nearrow 1}[x_{;1}]^\bullet$. By d), we have $b \sqsubseteq a_{;1}$, whence Lemma 8.5 implies $b_{1;} \sqsubseteq (a_{;1})_{1;} = a$. Since $b_{1;} \in \mathcal{G}[x]$, it follows that $a \sqsubseteq b_{1;} \sqsubseteq a$, whence $a = b_{1;}$. □

**Theorem 8.8.** *The class* **Ne** *is a surreal substructure.*

**Proof.** Let us first show that the root $a = \mathcal{G}[x]^\bullet$ of each halo with $x \in \mathbf{S}$ is a nested monomial. Indeed, Lemma 8.7(e) implies that $a_{;i} = \mathcal{G}_{\nearrow i}[x_{;i}]^\bullet$ for all $i \in \mathbb{N}$, by induction on $i$. In combination with Lemma 8.7(c), this yields $(a_{;i} - \varphi_i)/\epsilon_i \in \mathbf{Mo}$ for all $i \in \mathbb{N}$, as required.

In order to conclude that **Ne** coincides with the surreal substructure $\mathbf{Smp}_\mathcal{G}$, it remains to be shown that each halo contains at most one nested monomial. Given $a < b$ in **Ne**, it suffices to show that $\mathcal{G}[a] < b$. Let $i \in \mathbb{N}$ and $r \in \mathbb{R}^>$. If $\epsilon_{;i+1} = 1$, then $\epsilon_{;i} = \epsilon_i$ and $(a_{;i} - \varphi_i)/\epsilon_i < (b_{;i} - \varphi_i)/\epsilon_i$. Those are monomials, so $\mathbb{R}^>(a_{;i} - \varphi_i)/\epsilon_i < (b_{;i} - \varphi_i)/\epsilon_i$, whence $\Psi_{i,\mathbb{R}^>}(a) < b$. Similarly, if $\epsilon_{;i+1} = -1$, then $(b_{;i} - \varphi_i)/\epsilon_i < \mathbb{R}^>(a_{;i} - \varphi_i)/\epsilon_i$, whence again $\Psi_{i,\mathbb{R}^>}(a) < b$. Using Lemma 8.7(b), we conclude that $\mathcal{G}[a] < b$. □

## 8.4 Sufficient conditions for nestedness

Let $(\varphi, \epsilon)$ be a signed sequence. The conditions **AS** and **NS** may not be so easy to check for $(\varphi, \epsilon)$. Let us mention a few stronger sufficient conditions that imply **AS** and **NS**.

**Proposition 8.9.** *Let $(\varphi, \epsilon)$ be a signed sequence such that*

$$\forall i > 0, \forall j > i, \forall \psi \in \mathbf{No}^{\prec \mathrm{supp}\, \varphi_j}, \quad (\varphi_j \mathbin{+\mkern-8mu+} \psi)_{j;i} - \varphi_i \prec \mathrm{supp}\, \varphi_i.$$

*Then $(\varphi, \epsilon)$ is a nested sequence.*

**Proof.** The condition clearly implies the one from Proposition 8.3(*c*), which is equivalent to **AS**. Given $i \in \mathbb{N}$, let us next show that $\mathrm{supp}\, \varphi_i \succ \mathrm{e}^{\mathbf{S}_{\nearrow (i+1)}}$. Let $j > i$ be minimal with $\varphi_j \neq 0$. Given $\xi \in \mathbf{S}_{\nearrow (i+1)}$ and $\psi := \log_{j-(i+1)} \xi \in \mathbf{S}_{\nearrow j}$, we obtain $\psi \prec \mathrm{supp}\, \varphi_j$, whence $\mathrm{e}^\xi = (\psi_{j;i} - \varphi_i)/\epsilon_i \prec \mathrm{supp}\, \varphi_i$. □

**Example 8.10.** This proposition is in particular satisfied for the signed sequence $(\varphi, \epsilon)$ from the introduction with $\varphi_i = \sqrt{\log_i \omega}$ and $\epsilon_i = 1$ for all $i \in \mathbb{N}$.

**Example 8.11.** The proposition is also satisfied for any signed sequence $(\varphi, \epsilon)$ with $\varphi_{2i} = 0$ and $\epsilon_{2i} = 1$ for $i \in \mathbb{N}$ and $\varphi_{2i-1} = \log_{3i} \omega$ for $i > 0$.

Given a signed sequence $(\varphi, \epsilon)$ that satisfies a suitable condition **NS**\* (see below), Schmeling constructs a field of transseries that contains the corresponding nested transseries [38, Section 2.5]. Following [29, p. 6] and [12, p. 14], we conjecture that every field of transseries embeds into **No**. As part of our program to prove this conjecture, let us mention two more specific conjectures that concern nested transseries.

**Conjecture 8.12.** *Let $(\varphi, \epsilon)$ be a signed sequence such that the following holds:*

**NS**\*. $\forall i > 0, \forall \mathfrak{m} \in \mathrm{supp}\, \varphi_i, \exists j > i, \forall \psi \in \mathbf{No}_{\succ}^{\prec \mathrm{supp}\, \varphi_j}, \quad (\varphi_j \mathbin{+\mkern-8mu+} \psi)_{j;i} - \varphi_i \prec \mathfrak{m}.$

*Then $(\varphi, \epsilon)$ is a nested sequence.*

**Example 8.13.** The condition **NS**\* is satisfied for the sequence $((\log_i \omega)_{i \in \mathbb{N}^>}, ((-1)^i)_{i \in \mathbb{N}^>})$, which does not satisfy the condition from Proposition 8.9. It is also satisfied for

$$\varphi_0 = \sum_{k \in \mathbb{N}} \mathrm{e}^{\sqrt{\omega} - \mathrm{e}^{\sqrt{\log \omega} + \mathrm{e}^{\cdot^{\cdot^{\sqrt{\log_k \omega}}}}}},$$

$\epsilon_0 = -1$ and $\epsilon_i = 1$, $\varphi_i = \sqrt{\log_i \omega}$ for all $i > 0$. This sequence also does not satisfy the requirement of Proposition 8.9.

Let us finish with a counterexample of a signed sequence $(\varphi, \epsilon)$ that satisfies **AS** but not **NS**.

**Example 8.14.** Consider the nested sequence $((\sqrt{\log_i \omega})_{i \in \mathbb{N}}, (1)_{i \in \mathbb{N}})$ that gives rise to nested numbers of the form

$$x = \sqrt{\omega} \mathbin{+\mkern-8mu+} \mathrm{e}^{\sqrt{\log \omega} + \mathrm{e}^{\sqrt{\log_2 \omega} + \mathrm{e}^{\cdot^{\cdot^{\cdot}}}}}.$$

Given such a number $x$, we define $(\varphi_0, \epsilon_0) := (x - \sqrt{\omega}, 1)$ as well as $(\varphi_i, 1) := (\sqrt{\log_i \omega}, 1)$ for all $i \in \mathbb{N}$. By definition, $(\varphi_{\nearrow 1}, \epsilon_{\nearrow 1})$ is nested so there is $u \in \mathbf{Ne}_{\nearrow 1}$ with $u < \log(x - \sqrt{\omega})$. The number $\varphi_0 + \mathrm{e}^u$ lies in **Ne**, so the sequence $(\varphi, \epsilon)$ is admissible. However we have $\mathrm{e}^v = \varphi_0$, so $\mathrm{e}^v \not\prec \mathrm{supp}\, \varphi_0$. This means that $\varphi_0 + \mathrm{e}^v$ does not lie in **S** and thus that $(\varphi, \epsilon)$ is not nested.

There even exist admissible sequences $(\varphi, \epsilon)$ with $\mathbf{Ne} = \emptyset$. However, we conjecture that

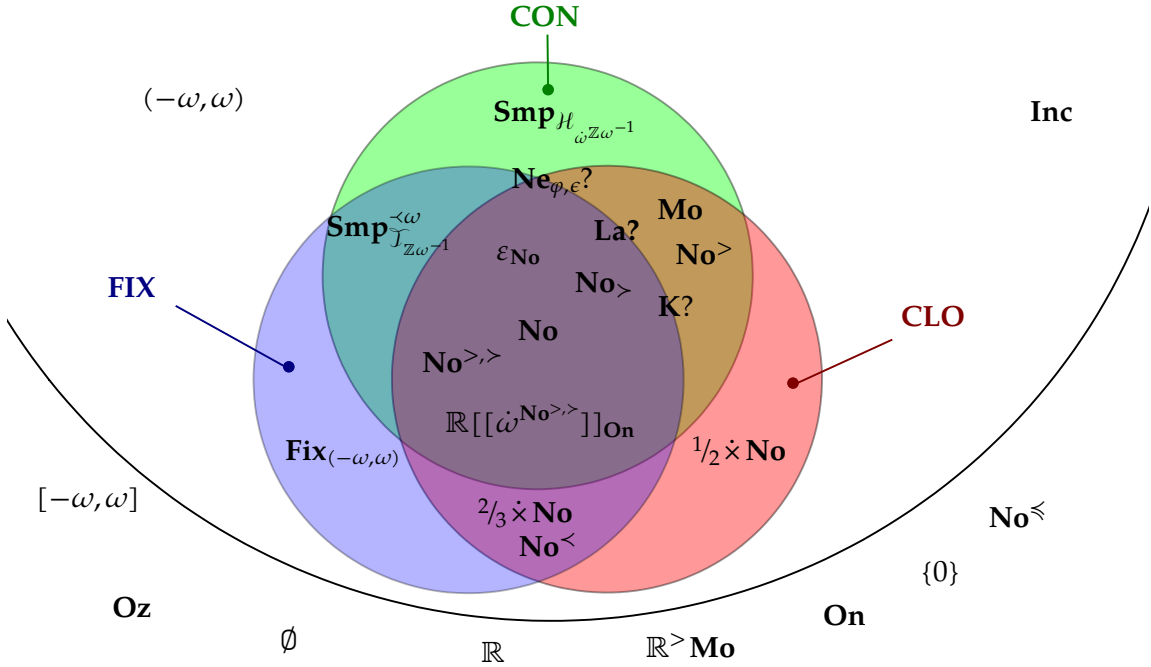

**Figure A.1.** A tiny glimpse of the landscape of surreal substructures.

**Conjecture 8.15.** *For every admissible sequence $(\varphi, \epsilon)$, there exists a $k \in \mathbb{N}$ such that $(\varphi_{\nearrow k}, \epsilon_{\nearrow k})$ is nested.*

We have made good progress on Conjectures 8.12 and 8.15 in the more general setting of hyperseries. We plan to report on this in a forthcoming paper.

# Appendix A  An atlas of surreal substructures

We have encountered several types of surreal substructures: intervals and convex surreal substructures, $\sqsubseteq$-final substructures, structures of fixed points, and structures obtained through convex partitions or group actions. Those different families of surreal substructures have non-trivial intersections. Figure A.1 gives a glimpse of the resulting landscape. We have used the following criteria for our classification:

- Surreal substructures lie in the great circle.
- **No**-closed surreal substructures lie in the rightmost smaller circle (**CLO**).
- Structures obtained through convex partitions of convex subclasses of **No** lie in the middle-upper smaller circle (**CON**).
- Structures of fixed points lie in the leftmost smaller circle (**FIX**).

All the represented classes in Figure A.1 satisfy the property that their non-empty cuts are rooted, which is not the case for other simple classes such as $\mathbf{No} + 1$. Equivalently, they are uniquely $(\leqslant, \sqsubseteq)$-isomorphic to a $\sqsubseteq$-initial subclass of **No**.

Question marks indicate that we do not know whether **La** and **K** may be construed as structures of fixed points. The nature of $\mathbf{Ne}_{\varphi, \epsilon}$ may change as a function of $(\varphi, \epsilon)$; we assume that $(\varphi, \epsilon)$ is nested. The class **La** is **No**-closed, but this result is not entirely trivial. We derived it from a computation of sign sequences of log-atomic numbers which is too long to produce here.

| S | No | No$^>$ | No$^{>,>}$ | S$_{\varphi,\epsilon}$ |
|---|---|---|---|---|
| $\Pi[x], x \in \mathbf{S}$ | **Hull**$(x+\mathbb{Z})$ | **Hull**$(\mathbb{R}^> x)$ | **Hull**$(\exp_\mathbb{Z}(x))$ | **Hull**$((\mathbb{R}^> x_{;i})_{i;}: i \in \mathbb{N}^>)$ |
| **Smp**$_\Pi$ | **No**$_>$ | **Mo** | **K** | **Ne**$_{\varphi,\epsilon}$ |

Table A.1. Examples of surreal substructures that correspond to classes of $\Pi$-simplest elements.

| S | **No** | $a \dotplus \mathbf{No}$ | $a \dottimes \mathbf{No}, a > 0$ | **Mo** | **No**$_>$ | **Mo**$\prec$**No**$^{\exists \omega^{-1}}$ |
|---|---|---|---|---|---|---|
| **Fix**$_\mathbf{S}$ | **No** | $(a \dotttimes \omega) \dotplus \mathbf{No}$ | $(\sup_\sqsubseteq a \dotttimes a \dotttimes \cdots) \dotttimes \mathbf{No}$ | $\varepsilon_{\mathbf{No}}$ | $\mathbb{R}[[\dot{\omega}^{\mathbf{No}^{>,>}}]]_{\mathbf{On}}$ | $\mathbf{La} \cap (\dot{\omega}^{\omega^{-1}}, \omega)$ |

Table A.2. Examples of surreal substructures obtained as classes of fixed points.

| U$\prec$V | **No**$^{\exists \beta}$ | **No**$^>$ | **No**$_>$ | **Mo** | **La** | $\varepsilon_{\mathbf{No}}$ |
|---|---|---|---|---|---|---|
| **No**$^{\exists x}$ | **No**$^{\exists x + \beta}$ | **No**$^{\exists x + 1}$ | $x \dotplus \mathbf{No}_>$ | $x \dotplus \mathbf{Mo}$ | $x \dotplus \mathbf{La}$ | $x \dotplus \varepsilon_{\mathbf{No}}$ |
| **No**$_>$ | $(\omega \dotttimes \beta) \dotplus \mathbf{No}_>$ | **No**$^{>0}_>$ | $\mathbb{R}[[\dot{\omega}^{\mathbf{No}^{>1}}]]_{\mathbf{On}}$ | $\dot{\omega}^{\mathbf{No}^>}$ | ? | $\varepsilon_{\mathbf{No}}$ |
| **Mo** | $\dot{\omega}^\beta \dotttimes \mathbf{Mo}$ | **Mo**$^{>\omega}$ | ? | **Smp**$_\mathcal{D}$ | ? | $\varepsilon_{\mathbf{No}}$ |
| **La** | ? | **La**$^{>\omega}$ | **K** | ? | ? | ? |

Table A.3. Imbrications of various common surreal substructures. The symbols ? signify that we were unable so far to determine an intelligible description of the corresponding imbrication.

Next we give a few examples of surreal substructures that were obtained as $\Pi$-simplest elements for convex partitions, through fixed points, and as imbrications of other surreal substructures.

**Remark A.1.** The identity $\mathbf{La} \cap (\dot{\omega}^{\omega^{-1}}, \omega) = (\mathbf{Mo} \prec \mathbf{No}^{\exists \omega^{-1}})^{\prec \omega}$ is given as an illustration; we refer to [4] for a proof. This is also an intermediate step in our computation of sign sequences of log-atomic numbers. There is, for every purely infinite number $\theta$ and integer $n \in \mathbb{Z}$, a similar description of $\mathbf{La} \cap (\lambda_{\theta+n}, \lambda_{\theta+n+1})$ in terms of fixed points of certains simple surreal substructures.

# Appendix B  Set-theoretic issues

**Proper classes as sets**

Strictly speaking, statements such as "**No** forms a real closed field" *de facto* do not make sense. Indeed, **No** is a proper class and not a set, whereas the definition of real closed fields relies on set theory. The most common standard for set theory is ZFC, i.e. Zermelo–Fraenkel's axioms with the axiom of choice. From a foundational point of view, it is more convenient to base the theory of surreal numbers on Neumann-Bernays-Gödel's set theory with the axiom of global choice (NBG set theory for short), which is a conservative extension of ZFC [13, 21].

**Set-sized relativations**

In the other direction, many of the results from this paper that were derived for class-sized surreal substructures admit set-sized analogues. More precisely, given a regular infinite ordinal $\kappa$, then many statements about $(\mathbf{No}, \leqslant, \sqsubseteq)$ can be relativized to $(\mathbf{No}(\kappa), \leqslant, \sqsubseteq)$, in which case "sets of cardinality $<\kappa$" play a similar role with respect to "sets of cardinality $\kappa$" as general "sets" with respect to "proper classes".

For instance, a *surreal substructure of* $\mathbf{No}(\kappa)$ is a subset $S \subseteq \mathbf{No}(\kappa)$ such that the set $(L|R) \cap \mathbf{No}(\kappa)$ is rooted for any two subsets $L < R$ in $S$ with $|L|, |R| < \kappa$. In other words, the surreal substructures of $\mathbf{No}(\kappa)$ are the isomorphic copies of $(\mathbf{No}(\kappa), \leqslant, \sqsubseteq)$ inside itself, and they behave similarly to usual surreal substructures in many respects. In particular, if $\kappa$ is the cardinality of $\mathbf{No}$ in $\mathrm{ZFC}_{\kappa'}$ with $\kappa' > \kappa$ as above, then surreal substructures can actually be considered as set-sized relativations of this kind.

**Cofinality**

In ZFC, the *cofinality* $\mathrm{cof}(X, \leqslant)$ of a linearly ordered set $(X, \leqslant)$ is equivalently

- the least order type of a cofinal well-ordered subset of $(X, \leqslant)$,
- the least cardinal of a cofinal subset of $(X, \leqslant)$,
- the unique regular ordinal which embeds in a cofinal way in $(X, \leqslant)$.

Assuming NBG set theory and regarding $\mathbf{On}$ as an initial, regular ordinal, this definition naturally extends to proper classes. In particular, every convex subclass $\mathbf{X}$ of a surreal substructure $\mathbf{S}$ has a cofinality $\mathrm{cof}(\mathbf{X}, \leqslant)$ in $\mathbf{On} \cup \{\mathbf{On}\}$, and elementary properties of the cofinality apply in our case. For instance, mutually cofinal convex subclasses of $\mathbf{No}$ have the same cofinality.

# Glossary

# Index